\documentclass[final,1p,times]{elsarticle}

\usepackage{amssymb}
\usepackage{amsmath} 
\usepackage{graphicx}
\usepackage{hyperref}
\usepackage{booktabs}
\usepackage{listings}
\usepackage{xcolor}
\usepackage{ragged2e}
\usepackage{float}
\usepackage{algorithm} 
\usepackage{algpseudocode} 
\usepackage{caption}
\usepackage{setspace}
\usepackage[utf8]{inputenc}
\usepackage{lineno} 
\usepackage{xcolor}
\usepackage{siunitx}
\usepackage{graphicx}
\usepackage{subcaption}
\usepackage{placeins}
\usepackage{array}
\usepackage{booktabs}
\usepackage{caption}
\usepackage{algorithm}
\usepackage{algpseudocode}
\journal{}

\begin{document}

\begin{frontmatter}



\title{Steady-State and Transient Heat Conduction Analysis Using a Polygonal Cell-Based Smoothed Finite Element Method}

\author[1]{Miao Zhang}
\author[2]{Yang Yang}
\author[1]{Chao Su}
\author[1]{Yijia Dong}
\author[1]{Mingjiao Yan}
\affiliation[1]{
    organization={College of Water Conservancy and Hydropower Engineering, Hohai University},
    city={Nanjing},
    postcode={210098},
    state={Jiangsu},
    country={China}
}
\affiliation[2]{
    organization={PowerChina Kunming Engineering Corporation Limited},
    city={Kunming},
    postcode={650051},
    state={Yunnan},
    country={China}
}

\begin{abstract}
This paper presents a polygonal cell-based smoothed finite element method (CS-FEM) for two-dimensional steady-state and transient heat-conduction analysis. In the proposed formulation, Wachspress shape functions are employed to construct the temperature approximation over general polygonal elements, and the smoothed temperature gradient is evaluated through boundary integration over cell-based smoothing domains. The resulting formulation is implemented in ABAQUS through the user-defined element (UEL) interface, enabling heat-conduction analysis using polygonal meshes within a commercial finite element environment. Several numerical examples, including a linear patch test, steady-state benchmark problems, and transient heat-conduction problems with different boundary conditions, are investigated to verify the accuracy, convergence behavior, and robustness of the proposed method. The numerical results show good agreement with analytical or reference solutions. Compared with conventional FEM using triangular and quadrilateral elements, the proposed polygonal CS-FEM exhibits favorable accuracy and convergence performance, while providing greater flexibility in mesh generation for complex geometries. The proposed framework therefore offers an accurate and robust numerical approach for steady-state and transient heat-conduction analysis.
\end{abstract}



\begin{highlights}
\item A polygonal CS-FEM is developed for heat-conduction analysis.
\item Smoothed temperature gradients are evaluated by boundary integration.
\item The proposed formulation is implemented in ABAQUS through UEL.
\end{highlights}

\begin{keyword}
Heat conduction \sep Polygonal element \sep Cell-based smoothed finite element method \sep Wachspress shape function \sep ABAQUS UEL \sep Transient analysis
\end{keyword}

\end{frontmatter}




\section{Introduction}

Heat-conduction analysis is an important topic in many engineering and scientific fields, such as civil engineering~\cite{hung_anh_overview_2021}, aerospace engineering~\cite{chang_influence_2025,wang_experimental_2026}, mechanics~\cite{liu_time-discontinuous_2024,vila_real_thermo-elastoviscoplastic_2004}, and electrical engineering~\cite{zhang_review_2021}. In practical applications, heat conduction problems often involve complex geometries, heterogeneous materials, localized heat sources, and transient boundary conditions. Although analytical solutions are available for some idealized problems~\cite{belghazi_analytical_2010,pophillat_analytical_2020}, they are generally difficult to obtain when the geometry, material parameters, and boundary conditions become complicated. Therefore, robust numerical methods with high accuracy and flexible mesh adaptability are required for the analysis of steady-state and transient heat-conduction problems.

A variety of numerical methods have been developed for heat-conduction analysis, including the finite element method (FEM)~\cite{bilal_finite_2020,yuan_energy_2025}, finite difference method (FDM)~\cite{gu_generalized_2019,zhang_local_2025}, finite volume method (FVM)~\cite{li_efficient_2023,xu_numerical_2026}, meshless method (MLM)~\cite{meng_fast_2022,zhang_improved_2013,singh_numerical_2004}, boundary element method (BEM)~\cite{simoes_experimental_2012,yao_precise_2014,yu_precise_2015}, and singular boundary method (SBM)~\cite{gu_singular_2012,wei_aca-sbm_2016}. Among these approaches, the FEM has been widely used because of its rigorous theoretical basis, numerical stability, and general applicability to engineering problems~\cite{chaquet_solving_2025,chaquet_temperature_2025,malek_three-dimensional_2020}. However, the accuracy and convergence of the FEM are strongly affected by mesh quality. When conventional triangular or quadrilateral elements are severely distorted, the numerical accuracy may deteriorate significantly~\cite{feng_element_2018,gao_element_2017,kazemzadeh-parsi_solution_2009}. This limitation becomes more pronounced for domains with holes, curved boundaries, or local refinement requirements, where repeated remeshing may be required to obtain high-quality elements.

To reduce the dependence of numerical accuracy on mesh quality, the smoothed finite element method (SFEM) has attracted increasing attention in heat-conduction analysis~\cite{chong_update-grid_2018,li_time_2026,zhao_novel_2025,li_smoothed_2014,wu_novel_2025}. By introducing a smoothing operation over element-based, node-based, edge-based, or cell-based smoothing domains, the compatible gradient field can be replaced by a smoothed gradient field evaluated through boundary integration. For heat-conduction problems, this treatment is particularly useful because the heat flux is directly determined by the temperature gradient through Fourier's law. Therefore, an accurate and stable approximation of the temperature gradient is essential for both temperature prediction and heat-flux recovery. Compared with the conventional FEM, the S-FEM can reduce the sensitivity to element distortion and improve the robustness of numerical integration.

In addition to strain or gradient smoothing, polygonal and polyhedral discretization techniques provide another effective way to enhance mesh flexibility. Unlike triangular and quadrilateral elements, polygonal elements can have an arbitrary number of sides and can naturally accommodate complex boundaries, hanging nodes, and local mesh refinement. Therefore, polygonal elements have been increasingly used in numerical simulations involving complex geometries~\cite{cui_new_2018,nguyen_polygonal_2017,nguyen-hoang_new_2017,nguyen-xuan_polytopal_2019,zhang_modeling_2017}. For heat-transfer problems, Yan et al.~\cite{yan_polyhedral_2025} developed a polyhedral scaled boundary finite element method for three-dimensional heat-conduction analysis. Jiang et al.~\cite{jiang_n-side_2023} proposed an $n$-side cell-based smoothed finite element method for heat-transfer problems. Zhou et al.~\cite{zhou_new_2023} presented a general analytical polygonal boundary element method for three-dimensional transient nonlinear heat-conduction problems. These studies demonstrate the potential of polygonal and polyhedral methods in heat-transfer analysis.

Although significant progress has been achieved, a systematic polygonal cell-based smoothed finite element method for both steady-state and transient heat-conduction problems still deserves further investigation. In particular, the combination of polygonal Wachspress interpolation, cell-based smoothing domains, boundary-integral evaluation of the smoothed temperature gradient, transient heat-capacity formulation, and ABAQUS user-defined element implementation has not been sufficiently discussed in a unified framework. Moreover, because ABAQUS/CAE cannot directly post-process field variables associated with user-defined elements, a complete computational procedure including preprocessing, UEL-based analysis, and post-processing is necessary for practical applications. Therefore, it is meaningful to develop a polygonal CS-FEM framework that can be conveniently implemented in a commercial finite element platform and verified through both steady-state and transient benchmark problems.

In this study, a polygonal cell-based smoothed finite element method is developed for two-dimensional steady-state and transient heat-conduction analysis. In the proposed formulation, Wachspress shape functions are used to construct the temperature approximation over general polygonal elements, and cell-based smoothing domains are introduced to evaluate the smoothed temperature gradient through boundary integration. This treatment avoids the direct evaluation of shape-function derivatives inside polygonal elements and improves the robustness of heat-flux recovery. The formulation is further implemented in ABAQUS through the user-defined element (UEL) interface, enabling both steady-state and transient heat-conduction analyses within a commercial finite element environment. Several benchmark examples involving regular domains, complex geometries, transient responses, and heat-flux boundary conditions are presented to verify the accuracy, convergence behavior, and computational performance of the proposed method.

The remainder of this paper is organized as follows. Section~\ref{sec:2} presents the theoretical formulation of the polygonal CS-FEM for heat-conduction problems. Section~\ref{sec:3} describes the ABAQUS UEL implementation and the corresponding computational framework. Section~\ref{sec:4} verifies the proposed method through numerical examples. Finally, Section~\ref{sec:5} summarizes the main conclusions.

\section{Theoretical formulation of the polygonal CS-FEM for heat conduction problems}
\label{sec:2}

We consider a heat-conduction body occupying the domain $\Omega \subset \mathbb{R}^2$, bounded by $\Gamma=\Gamma_T\cup\Gamma_q$ with $\Gamma_T\cap\Gamma_q=\varnothing$, where $\Gamma_T$ and $\Gamma_q$ denote the portions of the boundary with prescribed temperature and prescribed outward heat flux, respectively. For a heterogeneous and possibly anisotropic medium, Fourier's law is written as
\begin{equation}
\mathbf{q} = -\mathbf{k}\nabla T ,
\label{eq:fourier_heat_CSFEM}
\end{equation}
where $T$ is the temperature field, $\mathbf{q}$ is the heat-flux vector, and $\mathbf{k}$ is the symmetric positive-definite thermal conductivity tensor. Together with energy conservation, the transient heat-conduction equation can be expressed as
\begin{equation}
\nabla \cdot (\mathbf{k}\nabla T) + r - \rho c \frac{\partial T}{\partial t} = 0 
\quad \text{in } \Omega ,
\label{eq:gov-heat-CSFEM}
\end{equation}
where $\rho$ is the density, $c$ is the specific heat capacity, and $r$ denotes the volumetric heat source term.

The boundary conditions are given by
\begin{equation}
T = \bar{T}
\quad \text{on } \Gamma_T ,
\label{eq:bc-temp-CSFEM}
\end{equation}
\begin{equation}
\mathbf{k}\nabla T \cdot \mathbf{n} = -\bar{q}_n
\quad \text{on } \Gamma_q ,
\label{eq:bc-flux-heat-CSFEM}
\end{equation}
where $\bar{T}$ is the prescribed temperature, $\bar{q}_n$ is the prescribed outward normal heat flux, and $\mathbf{n}$ is the outward unit normal vector. For transient analyses, the initial condition is prescribed as
\begin{equation}
T(\mathbf{x},0)=T_0(\mathbf{x}) \quad \text{in } \Omega .
\label{eq:ic-temp-CSFEM}
\end{equation}

\subsection{Galerkin weak form and spatial discretization}

Let $w$ be an admissible weighting function satisfying $w=0$ on $\Gamma_T$. Multiplying Eq.~\eqref{eq:gov-heat-CSFEM} by $w$ and integrating over $\Omega$ gives
\begin{equation}
\int_{\Omega} w \, \nabla \cdot (\mathbf{k}\nabla T)\, d\Omega
+ \int_{\Omega} w\, r \, d\Omega
- \int_{\Omega} w\, \rho c\, \frac{\partial T}{\partial t}\, d\Omega = 0 .
\label{eq:weak-raw-heat-CSFEM}
\end{equation}

Applying the divergence theorem to the first term yields
\begin{equation}
\int_{\Omega} w \, \nabla \cdot (\mathbf{k}\nabla T)\, d\Omega
= \int_{\Gamma} w\, \mathbf{k}\nabla T \cdot \mathbf{n} \, d\Gamma
- \int_{\Omega} (\nabla w)^T \mathbf{k}\nabla T \, d\Omega .
\label{eq:div-heat-CSFEM}
\end{equation}

With $w=0$ on $\Gamma_T$ and using the Neumann boundary condition on $\Gamma_q$, the Galerkin weak form is obtained as
\begin{equation}
\int_{\Omega} (\nabla w)^T \mathbf{k}\nabla T \, d\Omega
+ \int_{\Omega} w\, \rho c\, \frac{\partial T}{\partial t}\, d\Omega
= \int_{\Omega} w\, r \, d\Omega
+ \int_{\Gamma_q} w\, \bar{q}_n \, d\Gamma .
\label{eq:weak-form-heat-CSFEM}
\end{equation}

The domain $\Omega$ is discretized into polygonal elements $\{\Omega_e\}$. Within an element $\Omega_e$, the temperature and weighting functions are approximated as
\begin{equation}
T^h(\mathbf{x},t)=\sum_{I=1}^{n_e} N_I(\mathbf{x})\, T_I(t), 
\qquad 
w^h(\mathbf{x})=\sum_{J=1}^{n_e} N_J(\mathbf{x})\, \delta T_J ,
\label{eq:trial-test-heat-CSFEM}
\end{equation}
where $N_I$ are polygonal shape functions, $n_e$ is the number of nodes of the polygonal element, and $\delta T_J$ are arbitrary nodal variations satisfying the essential boundary conditions. In the present formulation, Wachspress shape functions are adopted for convex polygonal elements because they satisfy the partition of unity and linear completeness conditions. Substituting Eq.~\eqref{eq:trial-test-heat-CSFEM} into Eq.~\eqref{eq:weak-form-heat-CSFEM} and assembling all element contributions lead to the semi-discrete system
\begin{equation}
\mathbf{K}\mathbf{T}
+ \mathbf{C}\dot{\mathbf{T}}
= \mathbf{F} ,
\label{eq:semi-heat-CSFEM}
\end{equation}
where $\mathbf{T}$ collects the global nodal temperatures, $\dot{\mathbf{T}}$ is its time derivative, $\mathbf{K}$ is the global conductivity matrix, $\mathbf{C}$ is the global thermal capacity matrix, and $\mathbf{F}$ contains the source and Neumann boundary contributions.

\subsection{Construction of polygonal smoothing cells}

For each polygonal element $\Omega_e$, a cell-based smoothing strategy is employed. A polygonal element with $n_e$ vertices is divided into $n_c=n_e$ non-overlapping triangular smoothing cells. The representative point of the element is chosen as the polygon centroid $\mathbf{x}_o$. For a polygon with vertices $\mathbf{x}_1,\mathbf{x}_2,\ldots,\mathbf{x}_{n_e}$ arranged counterclockwise, the $C$-th smoothing cell is constructed as
\begin{equation}
\Omega_C = \triangle(\mathbf{x}_o,\mathbf{x}_C,\mathbf{x}_{C+1}),
\qquad C=1,2,\ldots,n_e ,
\label{eq:smoothing-cell-construction}
\end{equation}
where $\mathbf{x}_{n_e+1}=\mathbf{x}_1$. The element domain is therefore decomposed as
\begin{equation}
\Omega_e=\bigcup_{C=1}^{n_c}\Omega_C, 
\qquad 
\Omega_C \cap \Omega_D = \varnothing \quad (C\ne D),
\label{eq:smoothing-cell-union}
\end{equation}
where $\Omega_C$ has boundary $\Gamma_C$ and area $A_C$. This construction is convenient for arbitrary convex polygonal elements because only line integrations along the boundaries of the smoothing cells are required.

A piecewise-constant smoothing function is defined over each smoothing cell as
\begin{equation}
\Phi(\mathbf{x}-\mathbf{x}_C)=
\begin{cases}
1/A_C, & \mathbf{x}\in \Omega_C,\\
0, & \mathbf{x}\notin \Omega_C,
\end{cases}
\label{eq:smoothing-func-heat-CSFEM}
\end{equation}
where $\mathbf{x}_C$ is a representative point of the smoothing cell. The smoothed temperature gradient is defined as the weighted average of the compatible gradient over the smoothing cell:
\begin{equation}
\widetilde{\nabla} T (\mathbf{x}_C)
= \int_{\Omega_C} \nabla T(\mathbf{x}) \, \Phi(\mathbf{x}-\mathbf{x}_C)\, d\Omega .
\label{eq:smooth-grad-heat-CSFEM}
\end{equation}

Applying integration by parts to Eq.~\eqref{eq:smooth-grad-heat-CSFEM} gives
\begin{equation}
\widetilde{\nabla} T (\mathbf{x}_C)
= \int_{\Gamma_C} T(\mathbf{x})\, \mathbf{n}(\mathbf{x})\,\Phi(\mathbf{x}-\mathbf{x}_C)\, d\Gamma
- \int_{\Omega_C} T(\mathbf{x})\, \nabla\Phi(\mathbf{x}-\mathbf{x}_C)\, d\Omega .
\label{eq:smooth-parts-heat-CSFEM}
\end{equation}

Since $\Phi$ is constant inside $\Omega_C$, one has $\nabla\Phi=\mathbf{0}$ within the smoothing cell. Therefore, the smoothed gradient can be evaluated by the following boundary integral:
\begin{equation}
\widetilde{\nabla} T (\mathbf{x}_C)
= \frac{1}{A_C} \int_{\Gamma_C} T(\mathbf{x})\, \mathbf{n}(\mathbf{x}) \, d\Gamma .
\label{eq:smooth-boundary-heat-CSFEM}
\end{equation}
This expression avoids the direct evaluation of shape-function derivatives inside the polygonal element, which is one of the main advantages of the cell-based smoothed formulation.

Using the interpolation $T^h(\mathbf{x})=\sum_{I=1}^{n_e} N_I(\mathbf{x})T_I$, Eq.~\eqref{eq:smooth-boundary-heat-CSFEM} becomes
\begin{equation}
\widetilde{\nabla} T^h (\mathbf{x}_C)
= \sum_{I=1}^{n_e} \widetilde{\mathbf{B}}_I^{(C)}\, T_I ,
\label{eq:smooth-grad-shape-heat-CSFEM}
\end{equation}
where the smoothed gradient operator associated with node $I$ is defined as
\begin{equation}
\widetilde{\mathbf{B}}_I^{(C)}
= \frac{1}{A_C} \int_{\Gamma_C} N_I(\mathbf{x})\, \mathbf{n}(\mathbf{x}) \, d\Gamma .
\label{eq:BI-heat-CSFEM}
\end{equation}

For straight boundary segments $e=1,\dots,m_C$ of $\Gamma_C$, with length $L_e$ and outward unit normal vector $\mathbf{n}_e$, Eq.~\eqref{eq:BI-heat-CSFEM} can be rewritten as
\begin{equation}
\widetilde{\mathbf{B}}_I^{(C)}
= \frac{1}{A_C} \sum_{e=1}^{m_C} 
\left( \int_{e} N_I(\mathbf{x})\, d\Gamma \right)\mathbf{n}_e .
\label{eq:BI-edge-sum-heat-CSFEM}
\end{equation}

The line integral along each segment is evaluated by one-dimensional Gauss quadrature:
\begin{equation}
\int_{e} N_I(\mathbf{x}) \, d\Gamma
\approx 
\sum_{g=1}^{n_g} N_I(\mathbf{x}_g)\, w_g \, L_e ,
\label{eq:edge-gauss-heat-CSFEM}
\end{equation}
where $\mathbf{x}_g$ and $w_g$ are the Gauss point and weight, respectively.

\subsection{Element matrices and heat-flux recovery}

In the polygonal CS-FEM, both the trial and test function gradients in the conductivity term are replaced by their smoothed counterparts over each smoothing cell. Defining the cell-wise smoothed gradient matrix as
\begin{equation}
\widetilde{\mathbf{B}}^{(C)}
= 
\begin{bmatrix}
\widetilde{\mathbf{B}}_1^{(C)} &
\widetilde{\mathbf{B}}_2^{(C)} &
\cdots &
\widetilde{\mathbf{B}}_{n_e}^{(C)}
\end{bmatrix},
\label{eq:B-heat-CSFEM}
\end{equation}
the smoothed element conductivity matrix is obtained as
\begin{equation}
\mathbf{K}_e
= \sum_{C=1}^{n_c} 
A_C \, 
\widetilde{\mathbf{B}}^{(C)T} 
\mathbf{k} \,
\widetilde{\mathbf{B}}^{(C)} .
\label{eq:Ke-heat-CSFEM}
\end{equation}

The thermal capacity matrix is computed by standard area integration:
\begin{equation}
\mathbf{C}_e
= \int_{\Omega_e} N^T \rho c\, N \, d\Omega
= \sum_{C=1}^{n_c}
\int_{\Omega_C} N^T \rho c\, N \, d\Omega .
\label{eq:Ce-heat-CSFEM}
\end{equation}
It should be noted that the conductivity matrix is evaluated using the smoothed gradient operator, whereas the capacity matrix is evaluated using the shape functions themselves. This treatment preserves the consistency of the transient heat-storage term and avoids introducing artificial smoothing into the capacity contribution.

The elemental source vector and Neumann boundary vector are given by
\begin{equation}
\mathbf{f}_{r,e} 
= \int_{\Omega_e} N^T r \, d\Omega ,
\label{eq:fr-heat-CSFEM}
\end{equation}
\begin{equation}
\mathbf{f}_{q,e} 
= \int_{\Gamma_{q,e}} N^T \bar{q}_n \, d\Gamma .
\label{eq:fq-heat-CSFEM}
\end{equation}
Therefore, the elemental external vector can be written as
\begin{equation}
\mathbf{F}_e=\mathbf{f}_{r,e}+\mathbf{f}_{q,e}.
\label{eq:element-force-heat-CSFEM}
\end{equation}

After element assembly and enforcement of the essential boundary conditions on $\Gamma_T$, the global system takes the form of Eq.~\eqref{eq:semi-heat-CSFEM}.

The heat flux can be recovered directly from the smoothed temperature gradient. For the $C$-th smoothing cell, the recovered heat flux is computed as
\begin{equation}
\widetilde{\mathbf{q}}^{(C)}
=
-\mathbf{k}\widetilde{\nabla}T^h(\mathbf{x}_C)
=
-\mathbf{k}\widetilde{\mathbf{B}}^{(C)}\mathbf{T}_e ,
\label{eq:cell-heat-flux-CSFEM}
\end{equation}
where $\mathbf{T}_e$ is the vector of nodal temperatures of element $\Omega_e$. The element-averaged heat flux can be obtained by area-weighted averaging over all smoothing cells:
\begin{equation}
\widetilde{\mathbf{q}}_e
=
\frac{1}{A_e}
\sum_{C=1}^{n_c}
A_C \widetilde{\mathbf{q}}^{(C)} ,
\label{eq:element-heat-flux-CSFEM}
\end{equation}
where $A_e$ is the area of the polygonal element.

\subsection{Consistency and stability of the smoothed approximation}

The consistency of the proposed polygonal CS-FEM can be examined from the properties of the smoothed gradient operator. For a constant temperature field $T=T_0$, the exact temperature gradient is zero. From Eq.~\eqref{eq:smooth-boundary-heat-CSFEM}, one obtains
\begin{equation}
\widetilde{\nabla}T
=
\frac{T_0}{A_C}
\int_{\Gamma_C}\mathbf{n}\,d\Gamma
=
\mathbf{0},
\label{eq:constant-consistency-CSFEM}
\end{equation}
because the integral of the outward normal vector over a closed boundary is zero. Hence, the proposed formulation does not produce spurious heat flux for a constant temperature field.

For a linear temperature field,
\begin{equation}
T(\mathbf{x})=a+b x+c y ,
\label{eq:linear-temperature-CSFEM}
\end{equation}
the exact gradient is constant:
\begin{equation}
\nabla T =
\begin{bmatrix}
b\\
c
\end{bmatrix}.
\label{eq:exact-linear-gradient-CSFEM}
\end{equation}
Since Wachspress shape functions satisfy linear completeness on convex polygonal elements, the linear temperature field can be exactly represented by the nodal interpolation. Moreover, the boundary-integral form of the smoothed gradient gives the average gradient over the smoothing cell. Therefore, the exact constant gradient can be recovered for a linear temperature field, indicating that the proposed element satisfies the first-order consistency requirement and can pass the linear patch test.

The stability of the steady-state problem follows from the positive definiteness of the thermal conductivity tensor. For any admissible element temperature vector $\mathbf{T}_e$, the quadratic form of the element conductivity matrix is
\begin{equation}
\mathbf{T}_e^T\mathbf{K}_e\mathbf{T}_e
=
\sum_{C=1}^{n_c}
A_C
\left(
\widetilde{\mathbf{B}}^{(C)}\mathbf{T}_e
\right)^T
\mathbf{k}
\left(
\widetilde{\mathbf{B}}^{(C)}\mathbf{T}_e
\right)
\geq 0 .
\label{eq:positive-semi-definite-CSFEM}
\end{equation}
Thus, the element conductivity matrix is symmetric positive semi-definite. After applying sufficient essential boundary conditions, the assembled global conductivity matrix becomes positive definite, which ensures the stability and uniqueness of the steady-state solution.

\subsection{Time integration}

Let the time step be $\Delta t$ and denote the solution at time $t^n$ by $\mathbf{T}^n$. Using the generalized $\theta$-method, the time derivative and intermediate temperature are approximated as
\begin{equation}
\dot{\mathbf{T}}^{\,n+\theta} \approx \frac{\mathbf{T}^{n+1}-\mathbf{T}^{n}}{\Delta t},
\qquad 
\mathbf{T}^{n+\theta} = (1-\theta)\mathbf{T}^n+\theta \mathbf{T}^{n+1},
\label{eq:theta-method-heat-CSFEM}
\end{equation}
where $\theta\in(0,1]$. Applying the $\theta$-method to Eq.~\eqref{eq:semi-heat-CSFEM} gives
\begin{equation}
\left(\frac{1}{\Delta t}\mathbf{C}^{n+\theta}+\theta\,\mathbf{K}^{n+\theta}\right)\mathbf{T}^{n+1}
=
\mathbf{F}^{n+\theta}
+
\left(\frac{1}{\Delta t}\mathbf{C}^{n+\theta}-(1-\theta)\,\mathbf{K}^{n+\theta}\right)\mathbf{T}^{n}.
\label{eq:theta-scheme-heat-CSFEM}
\end{equation}

For the fully implicit backward Euler scheme with $\theta=1$, Eq.~\eqref{eq:theta-scheme-heat-CSFEM} reduces to
\begin{equation}
\left( \mathbf{K}^{n+1} + \frac{1}{\Delta t}\mathbf{C}^{n+1} \right)\mathbf{T}^{n+1}
=
\mathbf{F}^{n+1} + \frac{1}{\Delta t}\mathbf{C}^{n+1}\mathbf{T}^{n}.
\label{eq:backward-euler-heat-CSFEM}
\end{equation}
Since the backward Euler scheme is unconditionally stable for linear heat-conduction problems, it is adopted in the present transient analysis unless otherwise specified.

\section{Implementation}
\label{sec:3}

\subsection{Overview of the computational framework}

\begin{figure}[H]
  \centering
  \includegraphics[width=0.9\textwidth]{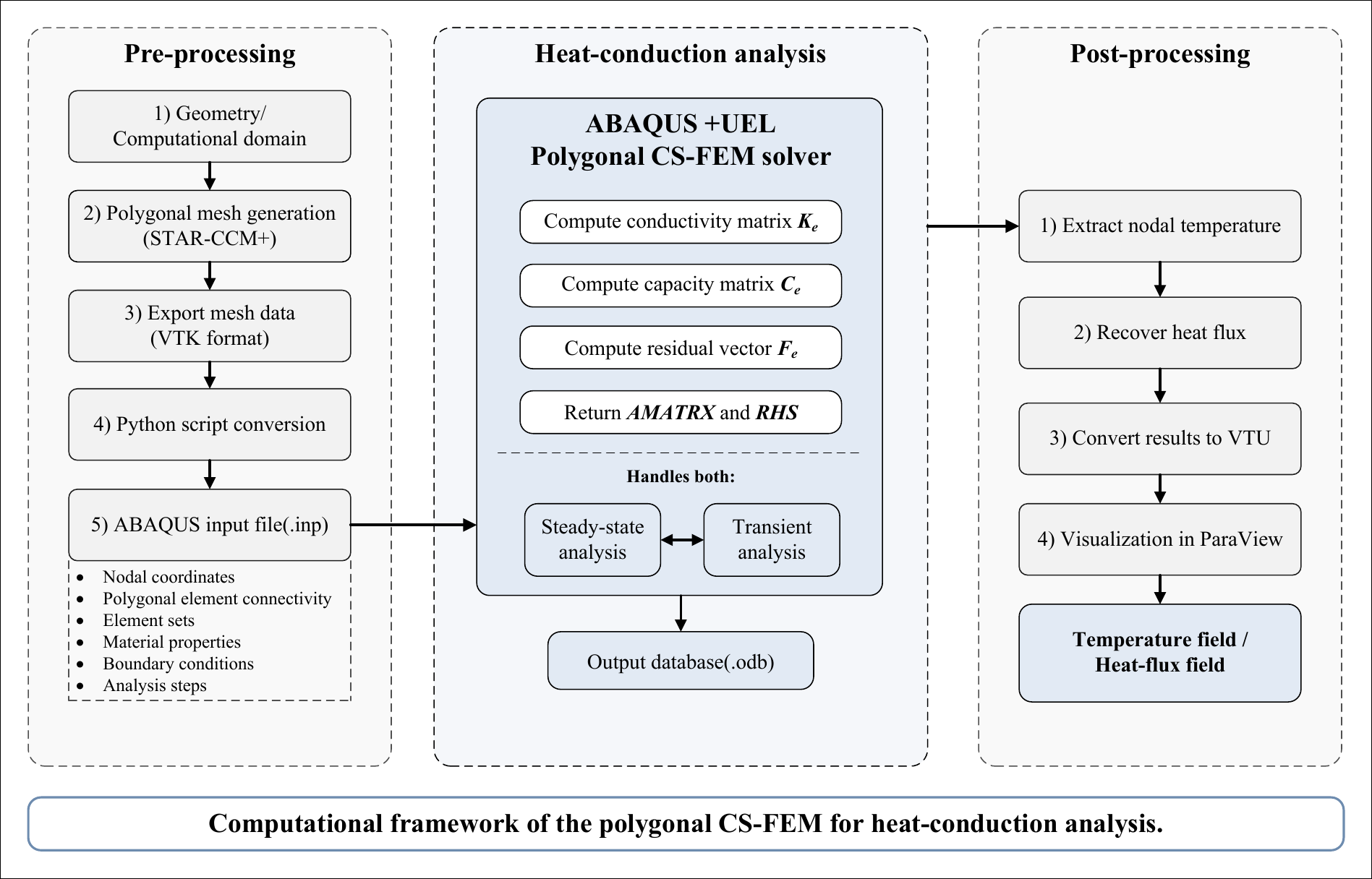}
  \caption{Overall workflow of the polygonal CS-FEM for heat-conduction analysis.}
  \label{fig:workflow}  
\end{figure}

The proposed polygonal CS-FEM is implemented in ABAQUS through the user-defined element (UEL) interface for steady-state and transient heat-conduction analyses. As shown in Fig.~\ref{fig:workflow}, the overall computational procedure consists of three main stages: pre-processing, heat-conduction analysis, and post-processing.

In the pre-processing stage, the computational domain is discretized using polygonal meshes generated by STAR-CCM+. The mesh data are then exported in VTK format and converted into an ABAQUS input file using an in-house Python script. The generated input file contains nodal coordinates, polygonal element connectivity, element sets, material parameters, boundary conditions, and analysis-step information. In the analysis stage, the element conductivity matrix, thermal capacity matrix, and residual vector are computed in the UEL subroutine and passed to ABAQUS through \texttt{AMATRX} and \texttt{RHS}. Since ABAQUS/CAE cannot directly visualize field results associated with UEL elements, the nodal temperature and recovered heat-flux quantities are extracted from the output database and converted into VTU files for visualization in ParaView.

\subsection{ABAQUS UEL implementation}

In the ABAQUS UEL framework, the nodal temperature is treated as the primary degree of freedom of each node. Although ABAQUS uses the array variable \texttt{U} to store the current nodal unknowns, it represents the nodal temperature vector in the present heat-conduction formulation. For each polygonal element, the UEL subroutine computes the elemental conductivity matrix $\mathbf{K}_e$, thermal capacity matrix $\mathbf{C}_e$, and external load vector $\mathbf{F}_e$, where
\begin{equation}
\mathbf{F}_e=\mathbf{f}_{r,e}+\mathbf{f}_{q,e}.
\end{equation}

For steady-state heat-conduction analysis, the elemental residual equation is written as
\begin{equation}
\mathbf{R}_e=\mathbf{F}_e-\mathbf{K}_e\mathbf{T}_e ,
\end{equation}
where $\mathbf{T}_e$ is the elemental nodal temperature vector. Therefore, the quantities passed to ABAQUS are
\begin{equation}
\mathrm{AMATRX}=\mathbf{K}_e ,
\label{eq:uel-steady-amatrx}
\end{equation}
\begin{equation}
\mathrm{RHS}=\mathbf{F}_e-\mathbf{K}_e\mathbf{T}_e .
\label{eq:uel-steady-rhs}
\end{equation}

For transient heat-conduction analysis, the backward Euler scheme is adopted. The elemental residual equation at time $t^{n+1}$ is given by
\begin{equation}
\mathbf{R}_e^{n+1}
=
\mathbf{F}_e^{n+1}
-
\mathbf{K}_e^{n+1}\mathbf{T}_e^{n+1}
-
\frac{1}{\Delta t}\mathbf{C}_e^{n+1}
\left(
\mathbf{T}_e^{n+1}-\mathbf{T}_e^n
\right).
\end{equation}
Accordingly, the tangent matrix and residual vector returned by the UEL are
\begin{equation}
\mathrm{AMATRX}
=
\mathbf{K}_e^{n+1}
+
\frac{1}{\Delta t}\mathbf{C}_e^{n+1},
\label{eq:uel-transient-amatrx}
\end{equation}
\begin{equation}
\mathrm{RHS}
=
\mathbf{F}_e^{n+1}
-
\mathbf{K}_e^{n+1}\mathbf{T}_e^{n+1}
-
\frac{1}{\Delta t}\mathbf{C}_e^{n+1}
\left(
\mathbf{T}_e^{n+1}-\mathbf{T}_e^n
\right).
\label{eq:uel-transient-rhs}
\end{equation}
Compared with the previous expression using a mass-like matrix $\mathbf{M}$, the notation $\mathbf{C}_e$ is used here to emphasize that it represents the thermal capacity matrix in the heat-conduction problem.

The main variables used in the ABAQUS UEL implementation are summarized in Tab.~\ref{tab:uel_variables}.

\begin{table}[H]
\centering
\caption{Main variables used in the ABAQUS UEL implementation.}
\label{tab:uel_variables}
\begin{tabular}{lll}
\toprule
Quantity & UEL variable & Description \\
\midrule
Current nodal temperature & \texttt{U} & Elemental temperature vector $\mathbf{T}_e^{n+1}$ \\
Temperature increment & \texttt{DU} & Increment of nodal temperature \\
Element coordinates & \texttt{COORDS} & Coordinates of polygonal element nodes \\
Material parameters & \texttt{PROPS} & Thermal conductivity, density, and specific heat \\
Element matrix & \texttt{AMATRX} & Conductivity matrix or transient tangent matrix \\
Residual vector & \texttt{RHS} & Elemental heat-balance residual vector \\
State variables & \texttt{SVARS} & User-defined output variables, if required \\
Time increment & \texttt{DTIME} & Time step size $\Delta t$ \\
\bottomrule
\end{tabular}
\end{table}

The implementation procedure of the polygonal CS-FEM in the UEL subroutine is summarized in Algorithm~\ref{alg:uel_csfem}.

\begin{algorithm}[H]
\caption{UEL implementation of the polygonal CS-FEM for heat-conduction analysis}
\label{alg:uel_csfem}
\begin{algorithmic}[1]
\State Read nodal coordinates, element connectivity, material parameters, and current nodal temperature.
\State Initialize $\mathbf{K}_e=\mathbf{0}$, $\mathbf{C}_e=\mathbf{0}$, and $\mathbf{F}_e=\mathbf{0}$.
\State Construct triangular smoothing cells for the polygonal element.
\For{each smoothing cell $\Omega_C$}
    \State Compute the area $A_C$ and the outward unit normal vector of each cell edge.
    \For{each boundary segment of $\Gamma_C$}
        \State Evaluate the Wachspress shape functions at the Gauss points.
        \State Compute the smoothed gradient operator $\widetilde{\mathbf{B}}^{(C)}$ by boundary integration.
    \EndFor
    \State Assemble the cell contribution to the conductivity matrix $\mathbf{K}_e$.
    \State Compute and assemble the contribution to the thermal capacity matrix $\mathbf{C}_e$.
\EndFor
\State Compute the heat-source vector $\mathbf{f}_{r,e}$ and the Neumann boundary vector $\mathbf{f}_{q,e}$, if present.
\State Set $\mathbf{F}_e=\mathbf{f}_{r,e}+\mathbf{f}_{q,e}$.
\If{steady-state analysis}
    \State Set $\mathrm{AMATRX}=\mathbf{K}_e$.
    \State Set $\mathrm{RHS}=\mathbf{F}_e-\mathbf{K}_e\mathbf{T}_e$.
\Else
    \State Set $\mathrm{AMATRX}=\mathbf{K}_e+\mathbf{C}_e/\Delta t$.
    \State Set $\mathrm{RHS}=\mathbf{F}_e-\mathbf{K}_e\mathbf{T}_e^{n+1}-\mathbf{C}_e(\mathbf{T}_e^{n+1}-\mathbf{T}_e^n)/\Delta t$.
\EndIf
\State Recover the smoothed heat flux for post-processing, if required.
\end{algorithmic}
\end{algorithm}

\subsection{Implementation of polygonal CS-FEM elements}

In the present implementation, a polygonal element with $n_e$ nodes is subdivided into $n_e$ triangular smoothing cells by connecting the polygon centroid to its vertices. This subdivision is only used for numerical integration and smoothing-domain construction; it does not introduce additional degrees of freedom. The unknown temperature field is interpolated by Wachspress shape functions defined over the original polygonal element.

The generated polygonal elements are assumed to be convex when Wachspress interpolation is adopted. For each polygonal element, the vertices are ordered counterclockwise to ensure a consistent definition of the outward normal vectors. The element centroid is first computed from the nodal coordinates. Then, each smoothing cell is constructed using the centroid and two adjacent polygon vertices. For every smoothing cell, the smoothed gradient matrix is evaluated by one-dimensional Gauss integration along the cell boundary. This procedure avoids the direct computation of compatible shape-function derivatives inside the polygon and improves the robustness of the element formulation for irregular polygonal meshes.

The recovered heat flux in each smoothing cell is calculated as
\begin{equation}
\widetilde{\mathbf{q}}^{(C)}
=
-\mathbf{k}\widetilde{\mathbf{B}}^{(C)}\mathbf{T}_e .
\end{equation}
For post-processing, the cell-wise heat flux can be further averaged over the polygonal element or projected to nodes. In the present work, the nodal temperature is directly obtained from the ABAQUS solution vector, whereas the heat-flux-related quantities are recovered from the smoothed temperature gradient.

Essential temperature boundary conditions are imposed using the standard ABAQUS \texttt{*Boundary} option. Prescribed heat fluxes are assembled into the elemental residual vector through the Neumann boundary contribution $\mathbf{f}_{q,e}$. Internal heat sources, if considered, are included through the source vector $\mathbf{f}_{r,e}$. In this way, both Dirichlet and Neumann boundary conditions can be consistently handled within the ABAQUS UEL framework.

\begin{figure}[H]
  \centering
  \includegraphics[width=0.9\textwidth]{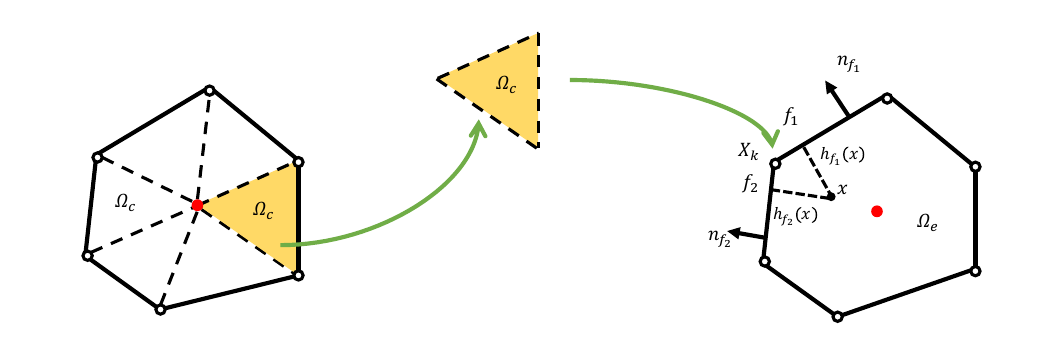}
  \caption{Construction of smoothing cells for a polygonal CS-FEM element.}
  \label{fig:polyelement}  
\end{figure}

\subsection{Definition of polygonal elements in ABAQUS}

The ABAQUS input file for the proposed UEL contains nodal coordinates, element connectivity, degrees of freedom, material parameters, boundary conditions, and analysis-step information. Since the polygonal elements are implemented as user-defined elements, different UEL types can be defined according to the number of vertices of the polygonal element, such as U3, U4, U5, $\ldots$, and U8. Taking an eight-node polygonal element U8 as an example, the corresponding input format is shown in Tab.~\ref{tab:1}.

\begin{table}[H]
\centering
\caption{Example of the ABAQUS input format for an eight-node polygonal UEL.}
\label{tab:1} 
\begin{tabular}{p{0.8cm}p{12cm}}
\toprule
Line & Input command \\
\midrule
1 & \texttt{*User element, nodes=8, type=U8, properties=4, coordinates=2, variables=8} \\
2 & \texttt{11} \\
3 & \texttt{*Element, type=U8, elset=E1} \\
4 & \texttt{1, 1, 2, 3, 5, 4, 9, 10, 12} \\
5 & \texttt{*Uel property, elset=E1} \\
6 & \texttt{1.0, 1.0, 1.0, 1.0} \\
\bottomrule
\end{tabular}
\end{table}

In Tab.~\ref{tab:1}, Line 1 defines the user element type, including the number of nodes, the number of material parameters, the spatial dimension, and the number of state variables. Line 2 specifies the active degree of freedom, where degree of freedom 11 is used as the temperature degree of freedom in ABAQUS. Lines 3 and 4 define the element connectivity. Lines 5 and 6 define the material parameters of the element set \texttt{E1}. In the present implementation, the material parameters are arranged as
\begin{equation}
\texttt{PROPS}=
\left[
k_x,\; k_y,\; \rho,\; c
\right],
\end{equation}
where $k_x$ and $k_y$ are the thermal conductivities in the $x$ and $y$ directions, respectively. For isotropic heat conduction, $k_x=k_y=k$.

\section{Numerical examples}
\label{sec:4}

In this section, five two-dimensional heat-conduction problems are investigated to examine the accuracy, convergence behavior, and applicability of the proposed polygonal CS-FEM. Unless otherwise specified, the thermal conductivity coefficient $k$, specific heat capacity $c$, and density $\rho$ are taken as unity. For comparison, conventional FEM analyses are performed in ABAQUS using DC2D3 triangular elements and DC2D4 quadrilateral elements.

To quantitatively evaluate the numerical accuracy, the relative temperature error is defined as
\begin{equation}
e=
\sqrt{
\frac{
\sum_{i=1}^{N_{\mathrm{node}}}
\left(T_i-T_i^{\mathrm{Ref}}\right)^2
}{
\sum_{i=1}^{N_{\mathrm{node}}}
\left(T_i^{\mathrm{Ref}}\right)^2
}
},
\label{eq:relative_temperature_error}
\end{equation}
where $T_i$ is the numerical temperature at node $i$, $T_i^{\mathrm{Ref}}$ is the corresponding analytical or reference solution, and $N_{\mathrm{node}}$ is the total number of nodes used for error evaluation.

\subsection{Linear patch test}

A linear patch test is first conducted to examine the consistency of the proposed polygonal CS-FEM. The computational domain is a unit square, $\Omega=[0,1]\times[0,1]~\mathrm{m}^2$. The thermal conductivity is taken as a constant, and no internal heat source is considered. The exact temperature field is prescribed as
\begin{equation}
T(x,y)=1+2x+3y .
\label{eq:patch_exact_temperature}
\end{equation}

The exact temperature given by Eq.~\eqref{eq:patch_exact_temperature} is imposed on all boundaries. Therefore, the nodal temperature on the boundary is assigned as
\begin{equation}
T_i=1+2x_i+3y_i .
\label{eq:patch_boundary_temperature}
\end{equation}

\begin{figure}[htbp]
\centering
\begin{subfigure}[b]{0.34\textwidth}
\centering
\includegraphics[width=1\textwidth]{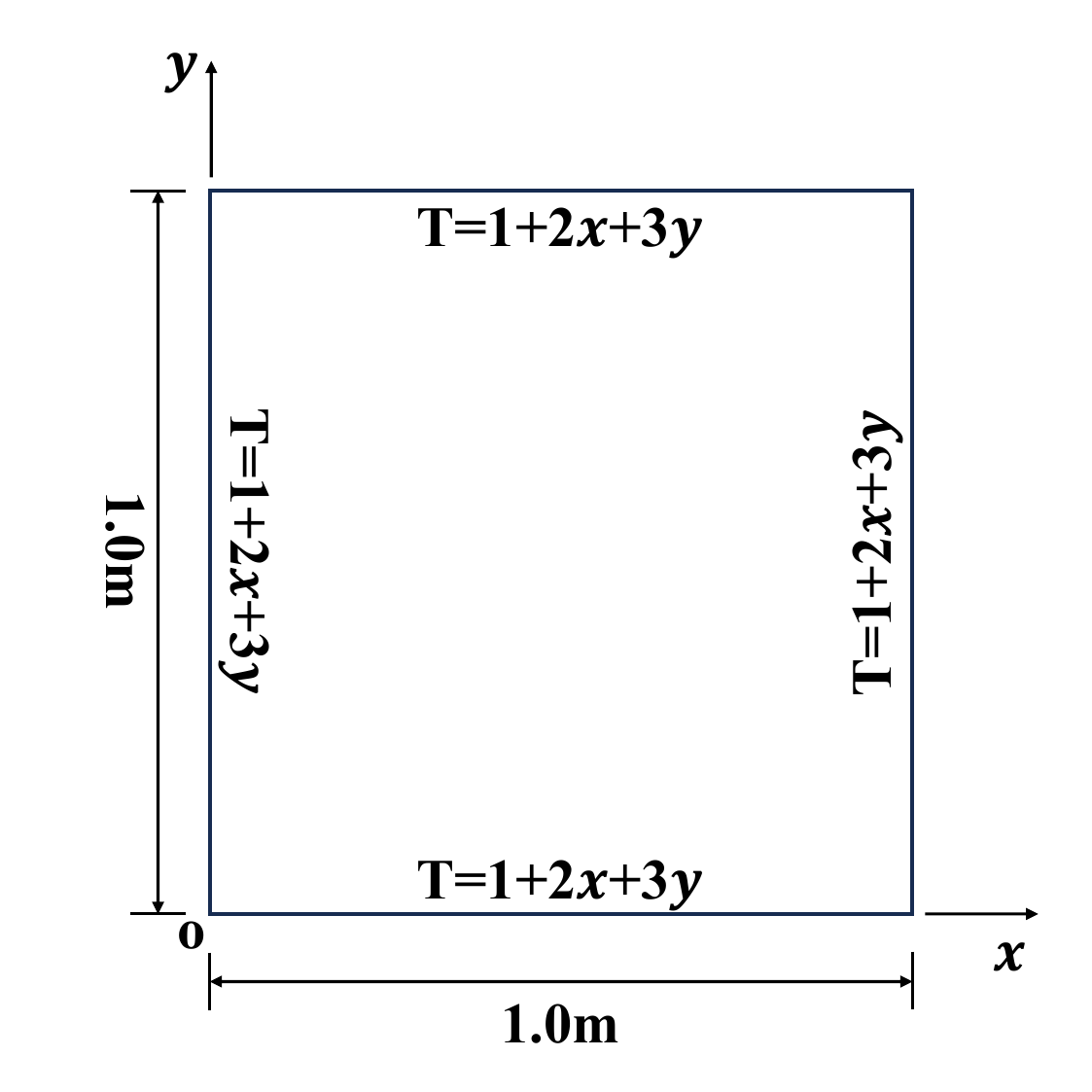}
\caption{}
\label{fig:geo_patch}
\end{subfigure}
\hfill
\begin{subfigure}[b]{0.3\textwidth}
\centering
\includegraphics[width=1\textwidth]{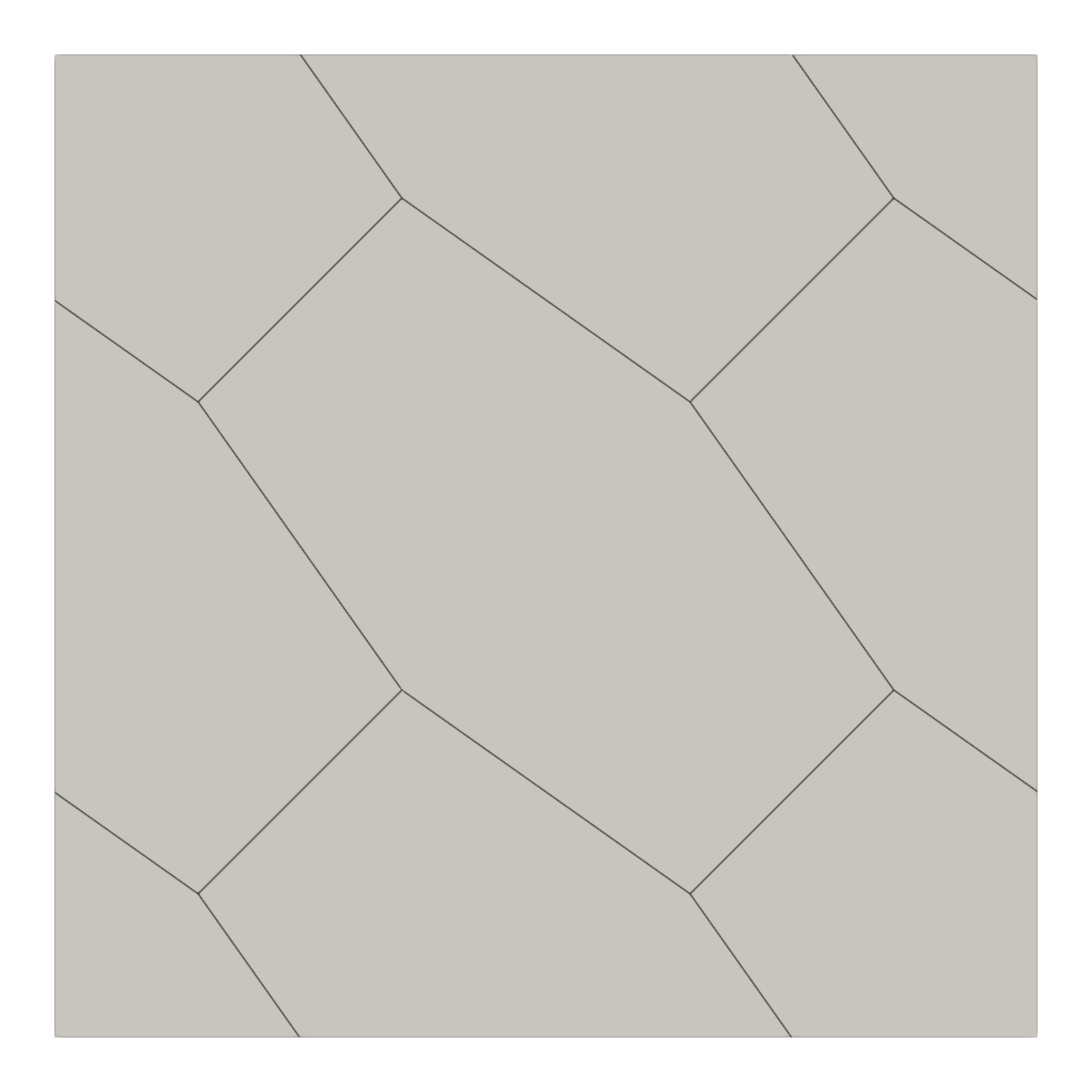}
\caption{}
\label{fig:patch_polyM}
\end{subfigure}
\hfill
\begin{subfigure}[b]{0.3\textwidth}
\centering
\includegraphics[width=1\textwidth]{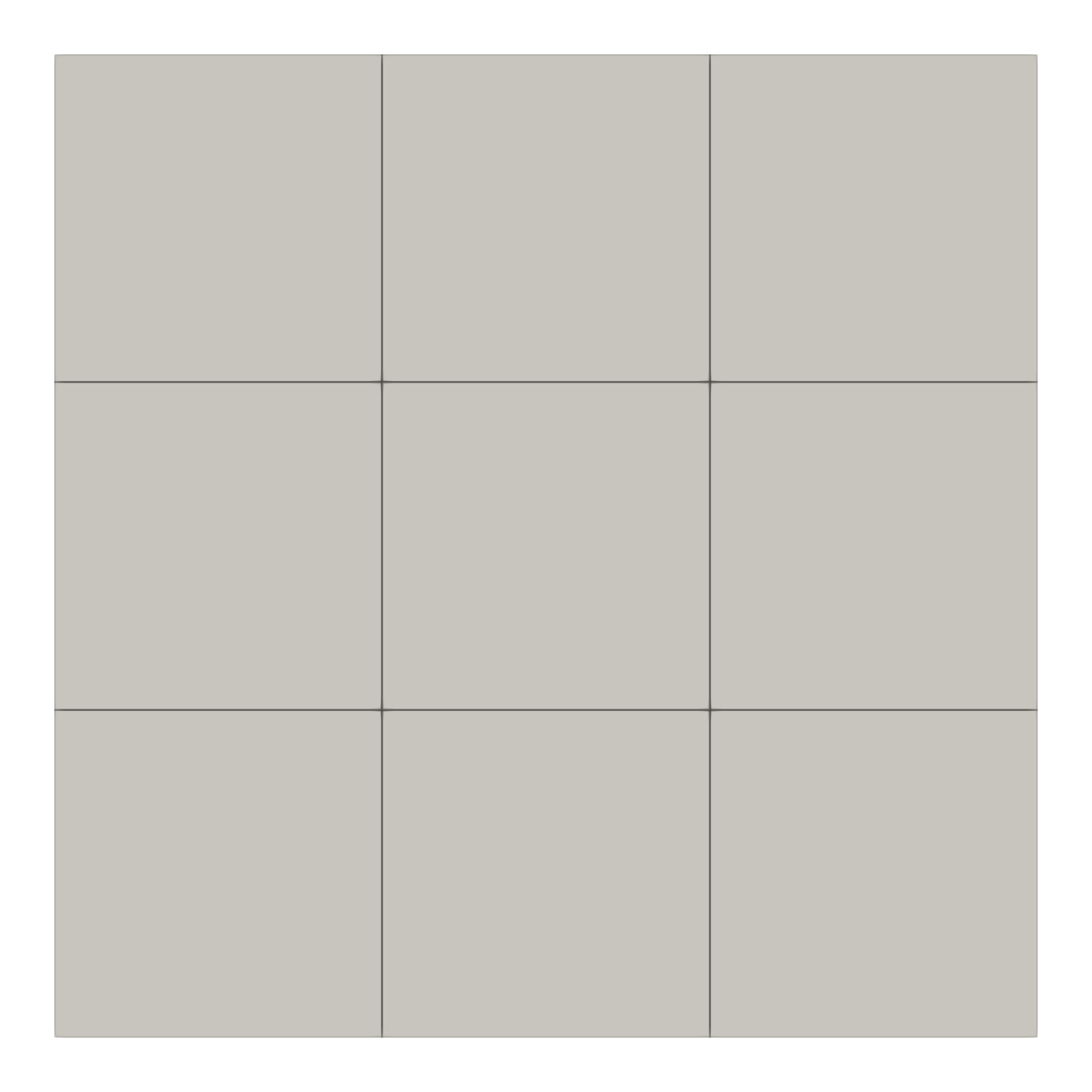}
\caption{}
\label{fig:patch_quadM}
\end{subfigure}
\caption{Computational model for the linear patch test: (a) geometry and boundary conditions; (b) polygonal mesh; and (c) quadrilateral mesh.}
\label{fig:linear_patch_model}
\end{figure}

Fig.~\ref{fig:linear_patch_model} shows the square patch with a side length of 1 m, where the boundary condition $T = 1 + 2x + 3y$ is applied on all boundaries. The patch is discretized using 9 polygonal elements and 9 quadrilateral elements, respectively.

\FloatBarrier

\begin{figure}[htbp]
\centering
\begin{subfigure}[b]{0.3\textwidth}
\centering
\includegraphics[width=1\textwidth]{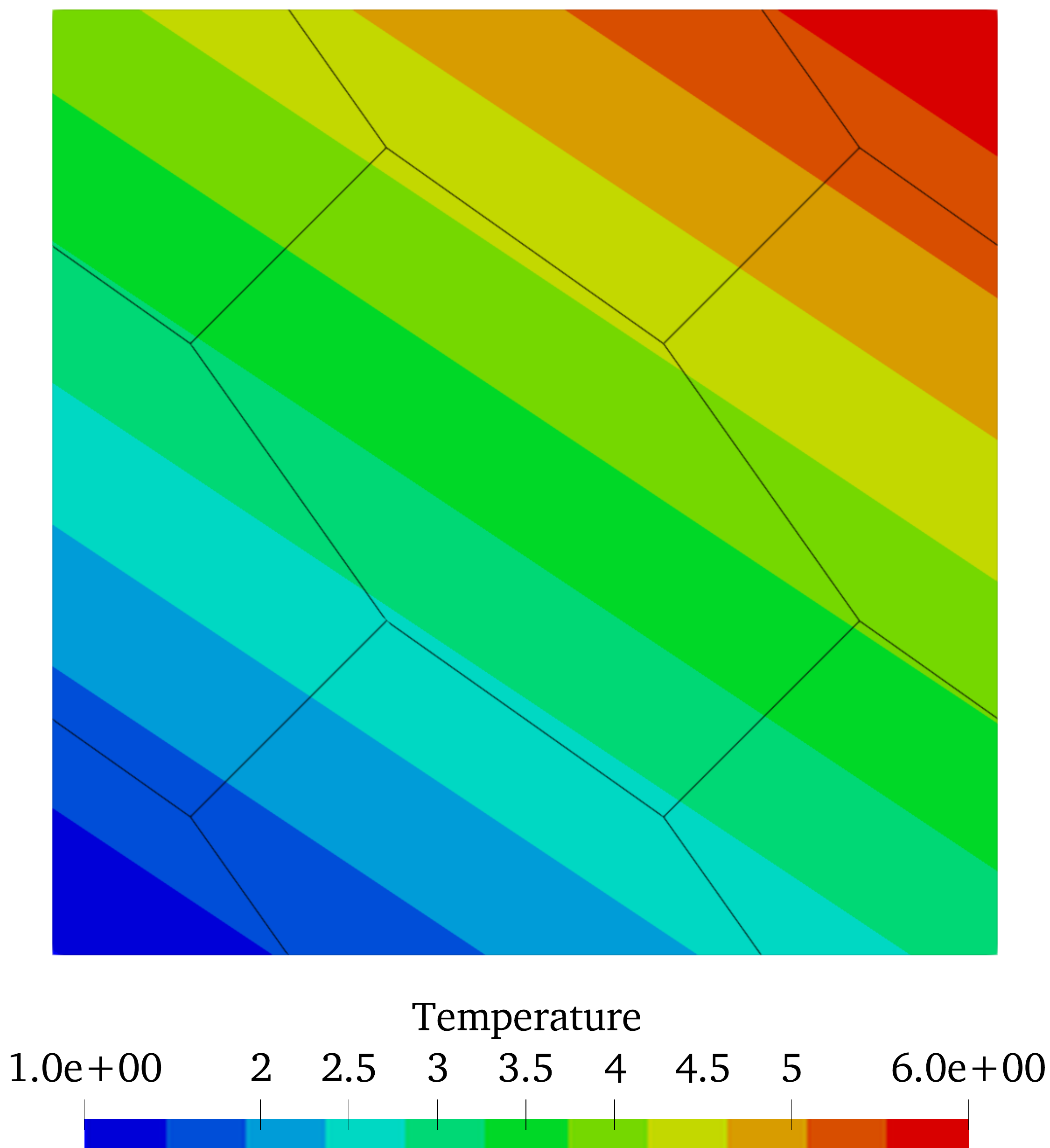}
\caption{}
\label{fig:patch_polyT}
\end{subfigure}
\hfill
\begin{subfigure}[b]{0.3\textwidth}
\centering
\includegraphics[width=1\textwidth]{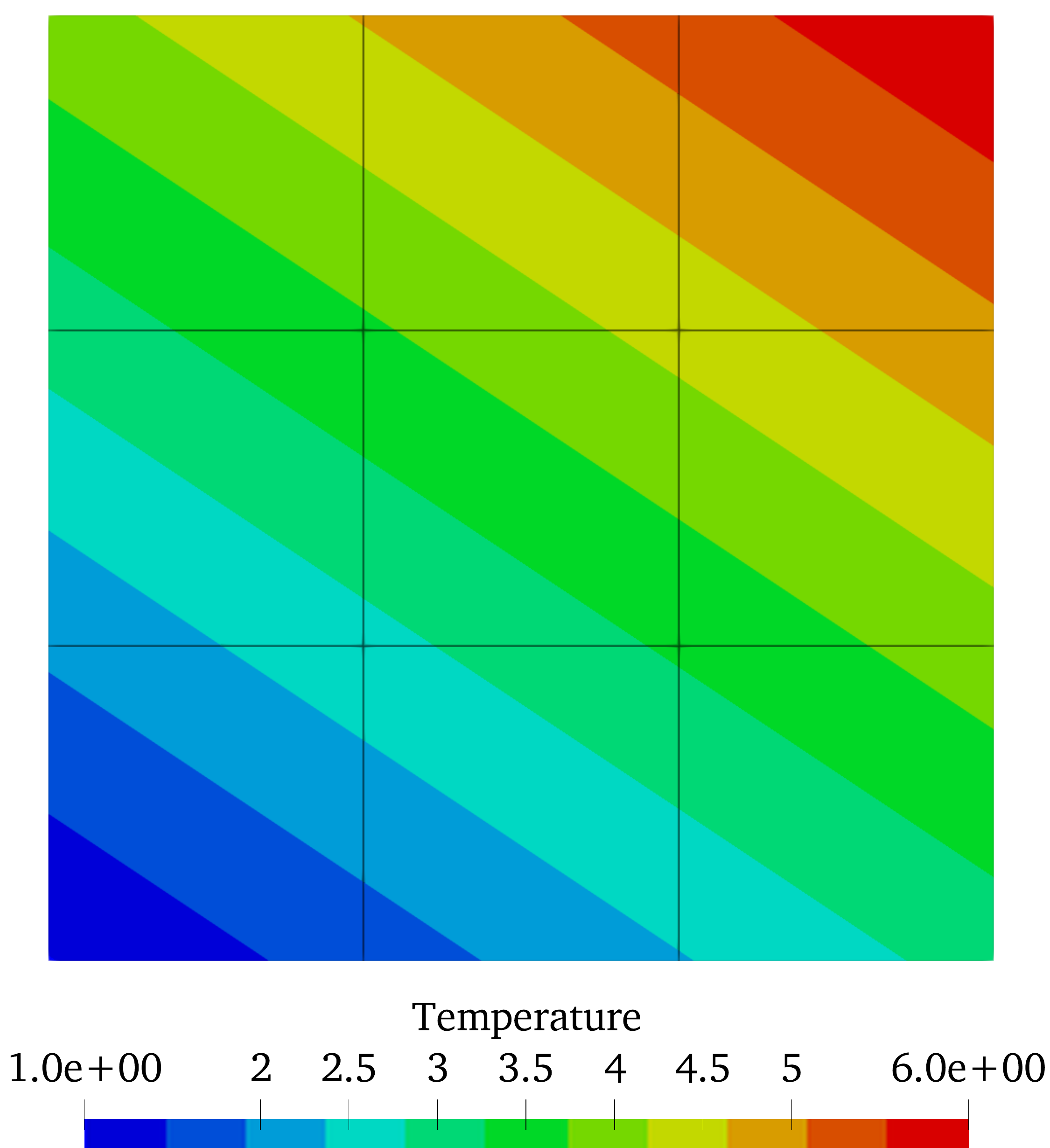}
\caption{}
\label{fig:patch_quadT}
\end{subfigure}
\hfill
\begin{subfigure}[b]{0.3\textwidth}
\centering
\includegraphics[width=1\textwidth]{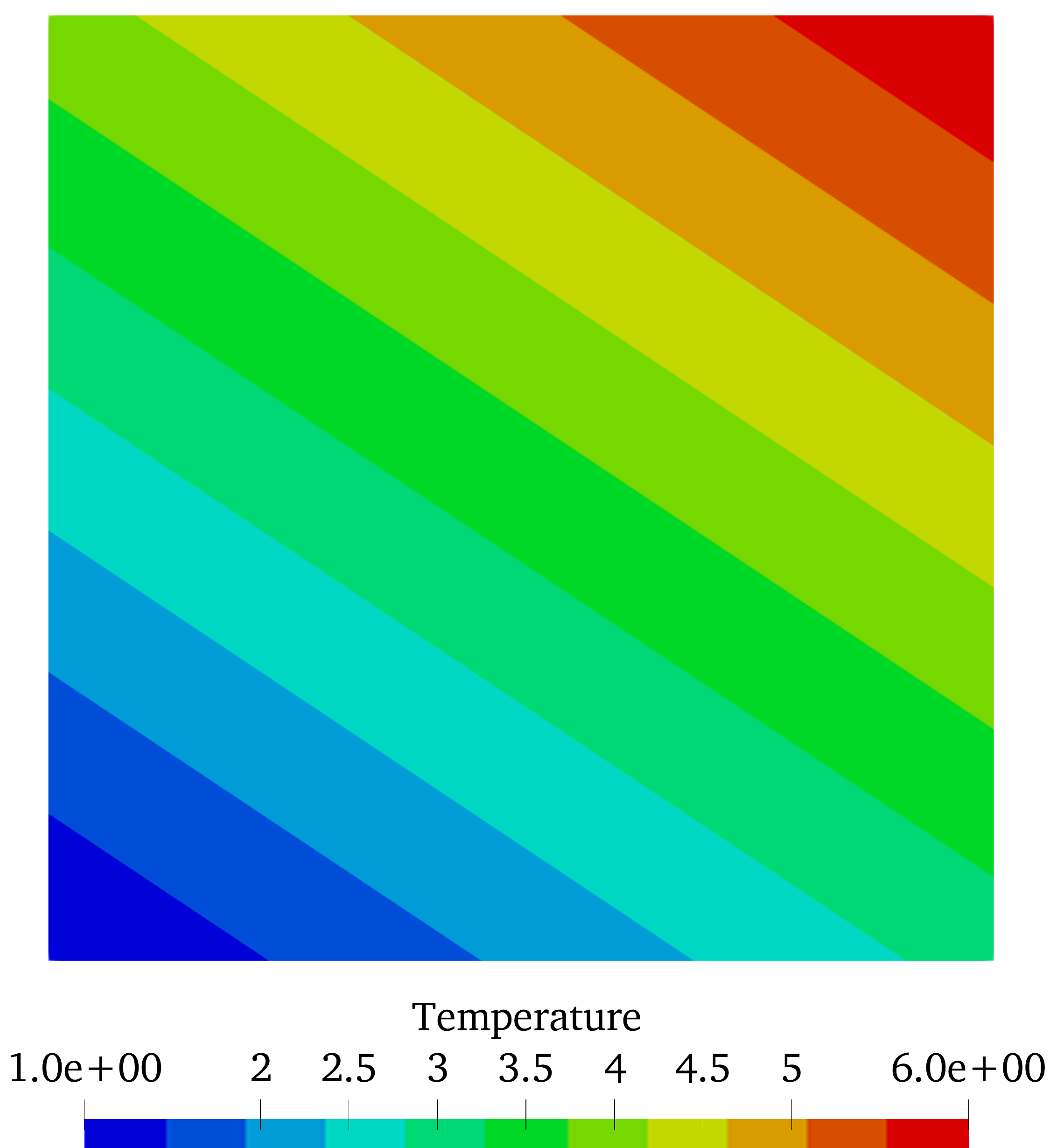}
\caption{}
\label{fig:patch_ref}
\end{subfigure}
\caption{Temperature contours for the linear patch test: (a) polygonal CS-FEM; (b) quadrilateral CS-FEM; and (c) analytical solution.}
\label{fig:linear_patch_temperature}
\end{figure}

Fig.~\ref{fig:linear_patch_temperature} presents the temperature contours obtained by the polygonal CS-FEM, the quadrilateral CS-FEM, and the analytical solution. The temperature distributions obtained by the two numerical methods are in excellent agreement with the analytical solution, indicating that both methods can accurately reproduce the prescribed linear temperature field.

\FloatBarrier

Tab.~\ref{tab:patch_temperature_error} lists the relative temperature errors of the two methods with respect to the analytical solution. The relative temperature errors of both methods are less than or close to $10^{-6}$, indicating that the numerical temperature fields agree very well with the analytical solution. Therefore, the proposed polygonal CS-FEM successfully passes the linear patch test and satisfies the first-order consistency requirement. The comparable accuracy between the polygonal and quadrilateral CS-FEM results also confirms the reliability of the polygonal formulation for linear heat-conduction problems.

\begin{table}[H]
\centering
\caption{Relative temperature errors for the linear patch test.}
\label{tab:patch_temperature_error}
\begin{tabular}{lc}
\toprule
Calculation method & Relative error \\
\midrule
Quadrilateral CS-FEM & $1.15 \times 10^{-6}$ \\
Polygonal CS-FEM     & $9.49 \times 10^{-7}$ \\
\bottomrule
\end{tabular}
\end{table}

\subsection{The steady-state heat conduction analysis}
\subsubsection{The steady-state heat conduction of square plate}
\begin{figure}[htbp]
 \centering
 \begin{subfigure}[b]{0.45\textwidth}
  \centering
  \includegraphics[width=1\textwidth]{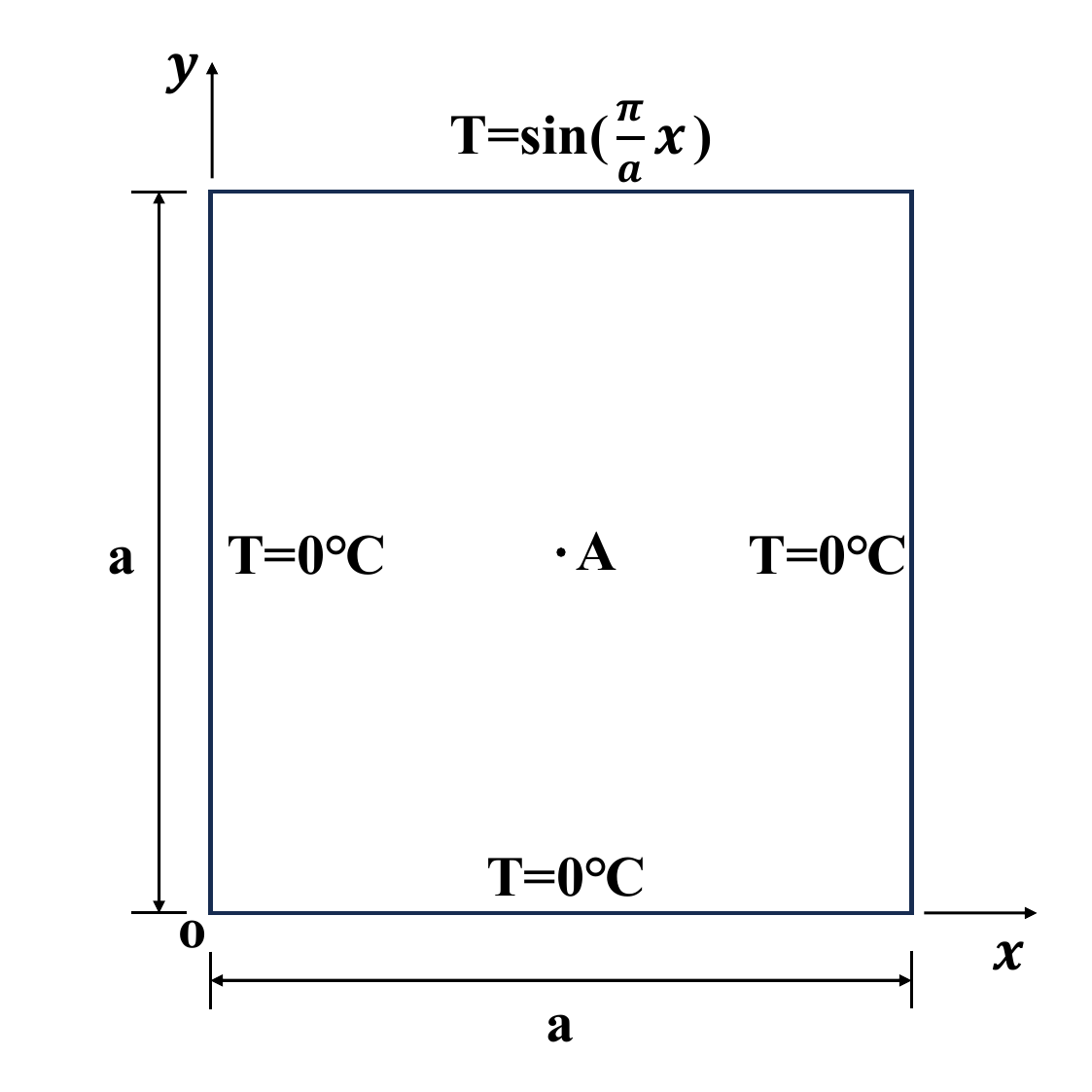}
  \caption{}
  \label{fig:geometry1}
 \end{subfigure}
 \hfill
 \begin{subfigure}[b]{0.45\textwidth}
  \centering
  \includegraphics[width=1\textwidth]{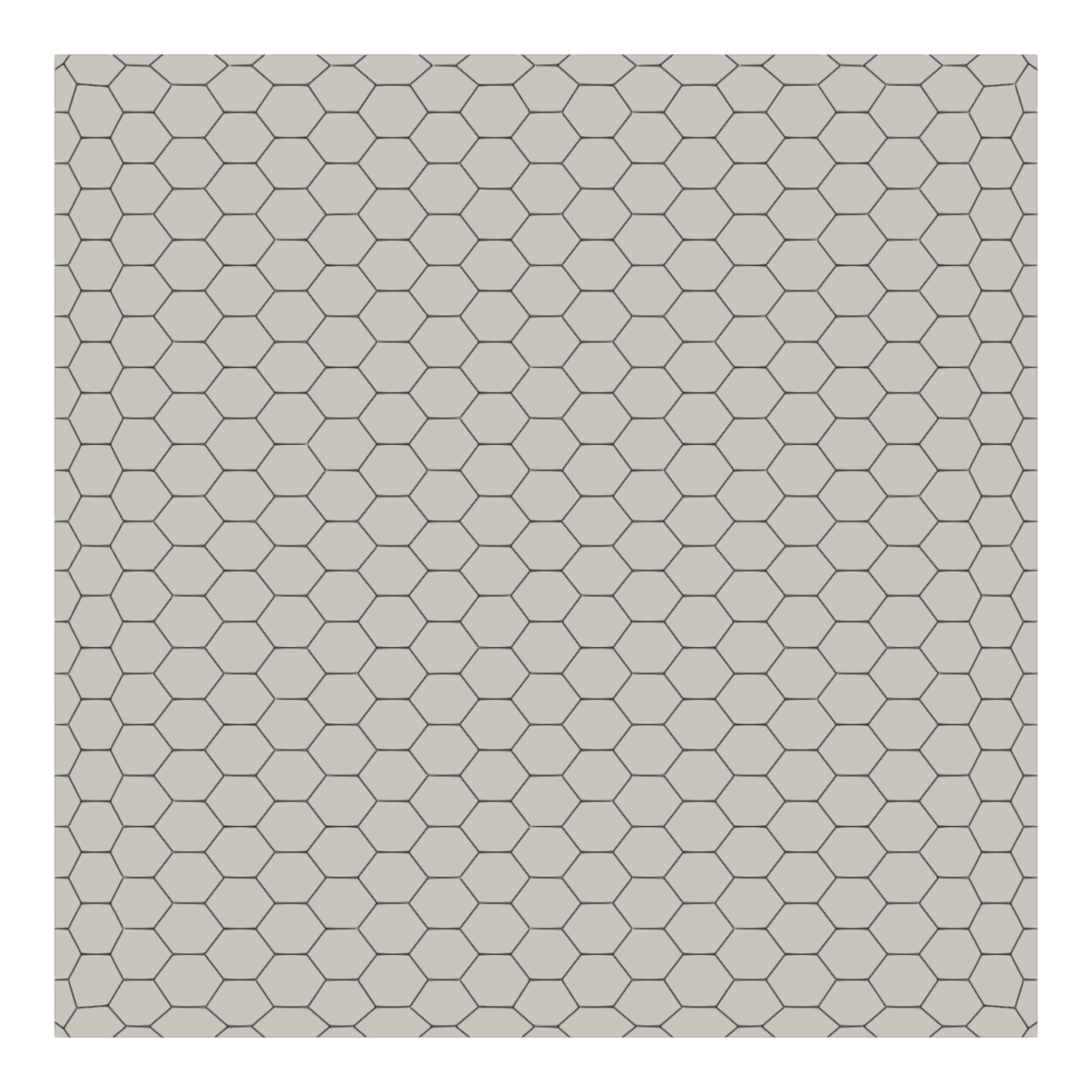}
  \caption{}
  \label{fig:1poly_M}
 \end{subfigure}
 \vspace{0cm}
  \begin{subfigure}[b]{0.45\textwidth}
  \centering
  \includegraphics[width=1\textwidth]{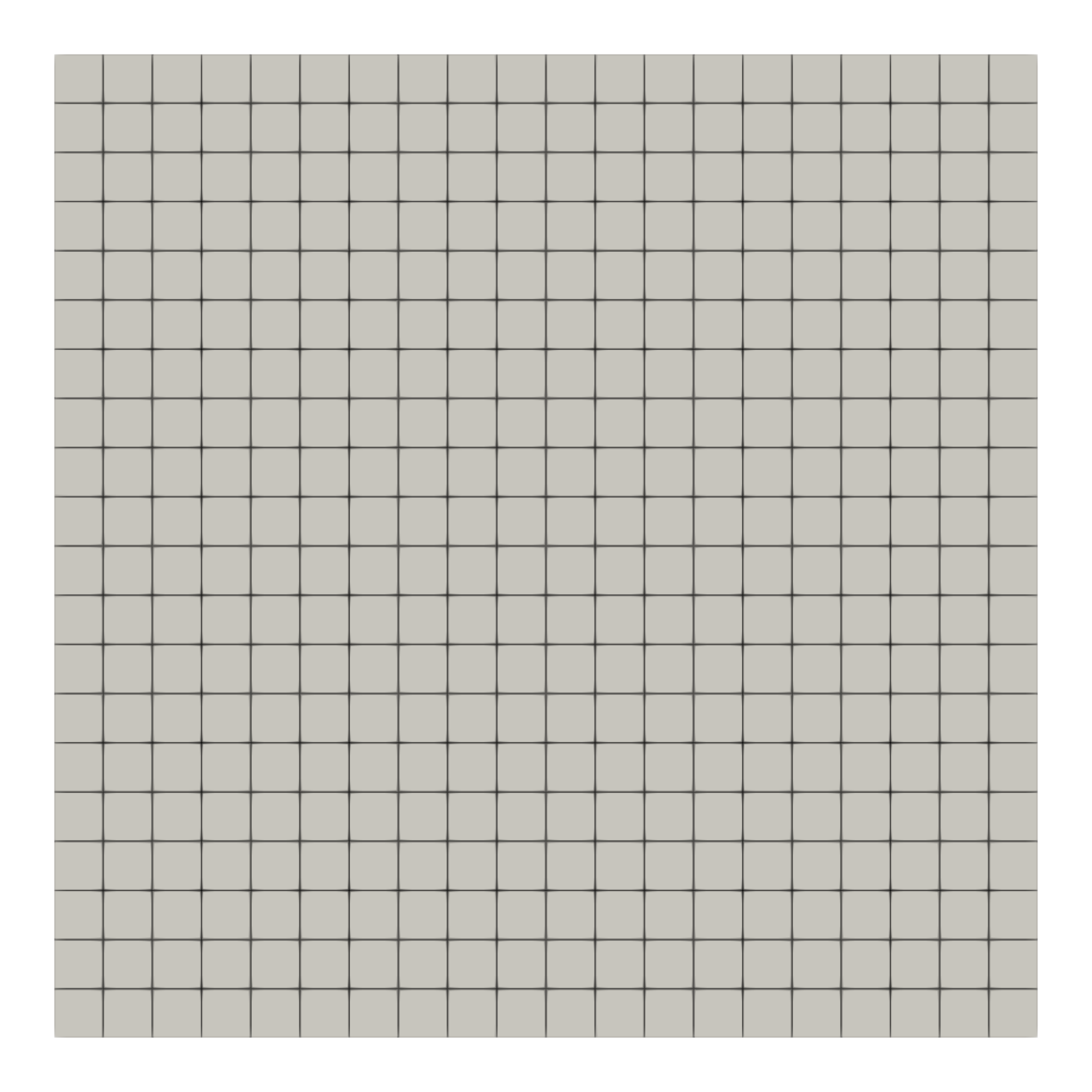}
  \caption{}
  \label{fig:quadM}
 \end{subfigure}
 \hfill
 \begin{subfigure}[b]{0.45\textwidth}
  \centering
  \includegraphics[width=1\textwidth]{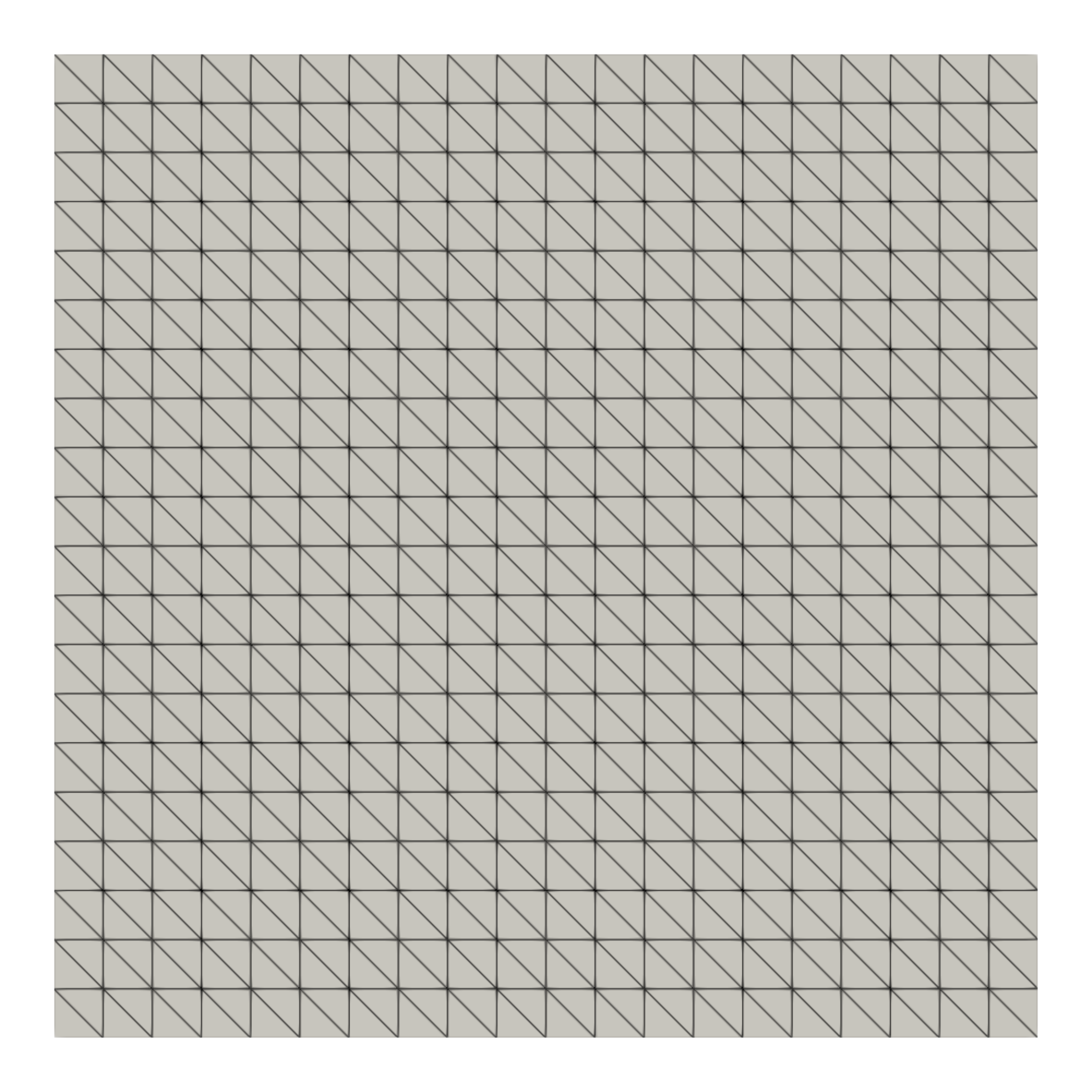}
  \caption{}
  \label{fig:triM}
 \end{subfigure}
\caption{Computational model of the square plate problem: (a) geometry and boundary conditions; (b) polygonal mesh; (c) quadrilateral mesh; and (d) triangular mesh.}
 \label{fig:The_model_of_the_square_plate}
\end{figure}

The first example considers a square plate with a side length of $a=2~\mathrm{m}$, subjected to prescribed temperature boundary conditions. The computational domain is discretized using triangular, quadrilateral, and polygonal elements, as illustrated in Fig.~\ref{fig:The_model_of_the_square_plate}. To assess the convergence behavior of the proposed method, a series of successively refined meshes with characteristic element sizes of $0.2~\mathrm{m}$, $0.1~\mathrm{m}$, $0.05~\mathrm{m}$, and $0.025~\mathrm{m}$ are employed. The temperature on the top boundary is prescribed as $T=\sin\left(\frac{\pi x}{a}\right),$ whereas zero temperature is imposed on the other three boundaries. The corresponding analytical solution \citep{nguyen_enhanced_2016} is given by

\begin{equation}
T=\frac{\sin (\pi x / a) \sinh (\pi y / a)}{\sinh (\pi)}
\end{equation}
\FloatBarrier
\FloatBarrier

\begin{figure}[htbp]
 \centering
 \begin{subfigure}[b]{0.3\textwidth}
  \centering
  \includegraphics[width=1\textwidth]{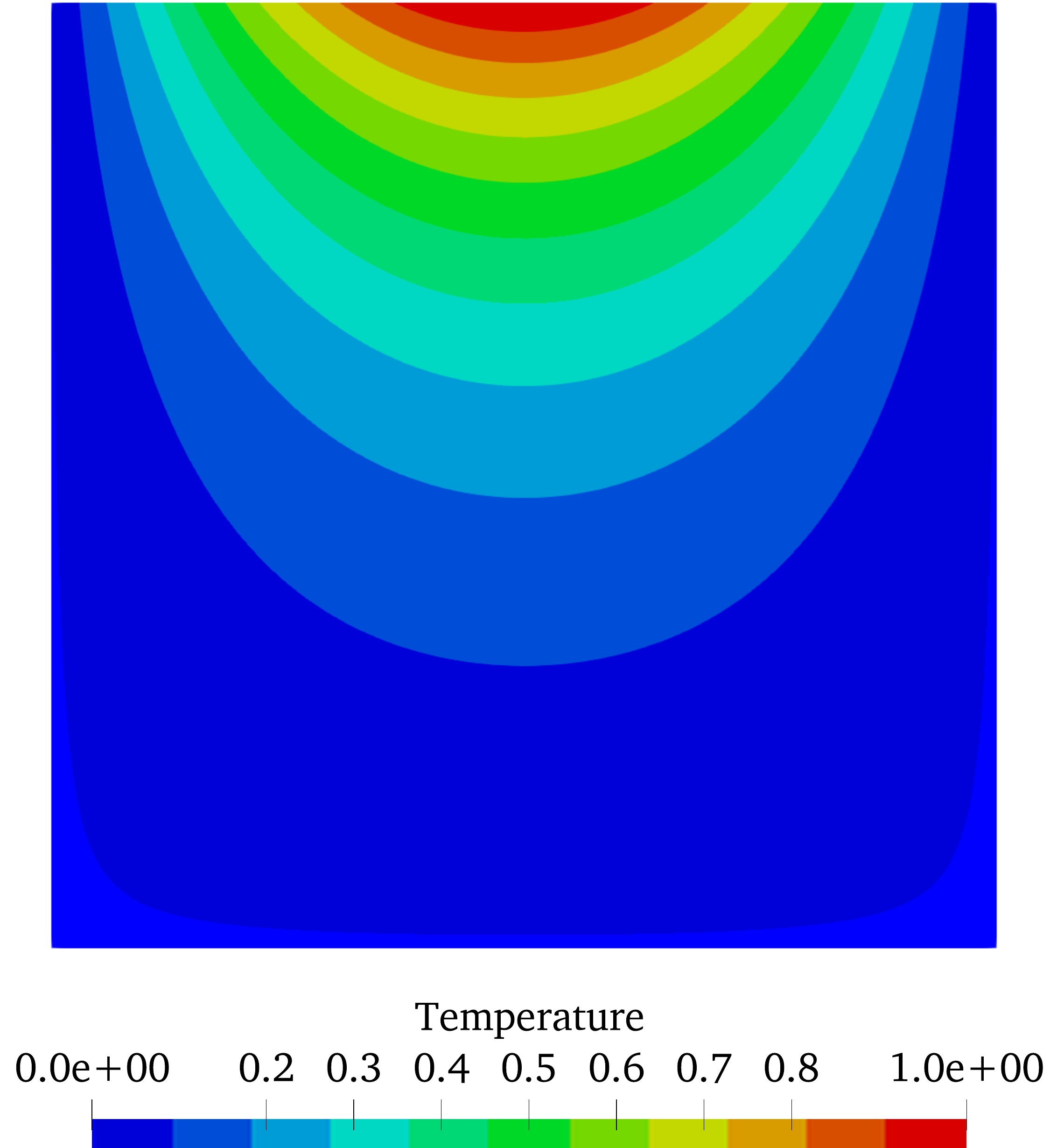}
  \caption{}
  \label{fig:1poly_T}
 \end{subfigure}
 \hfill
 \begin{subfigure}[b]{0.3\textwidth}
  \centering
  \includegraphics[width=1\textwidth]{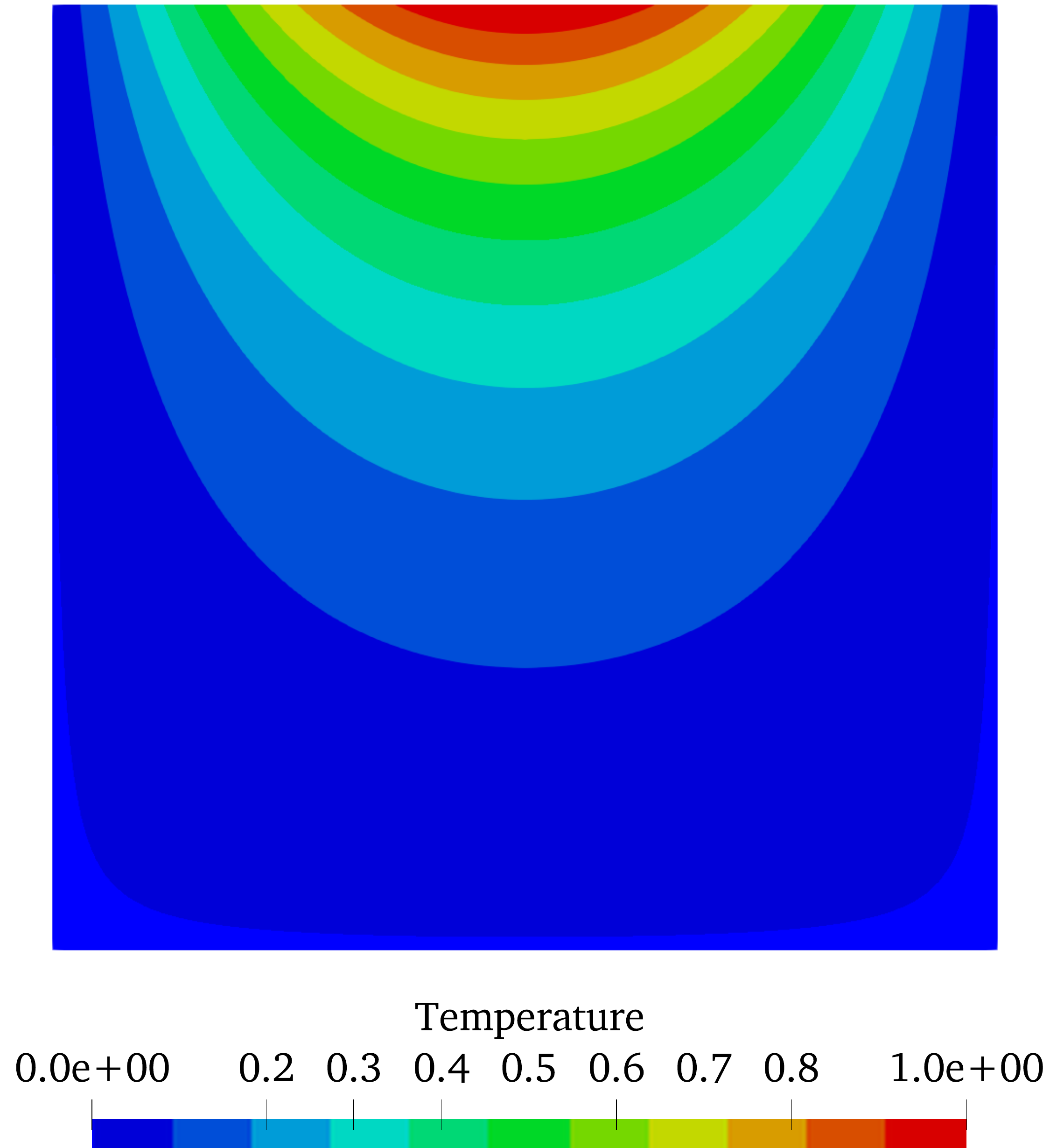}
  \caption{}
  \label{fig:1quad_T}
 \end{subfigure}
 \hfill
 \begin{subfigure}[b]{0.3\textwidth}
  \centering
  \includegraphics[width=1\textwidth]{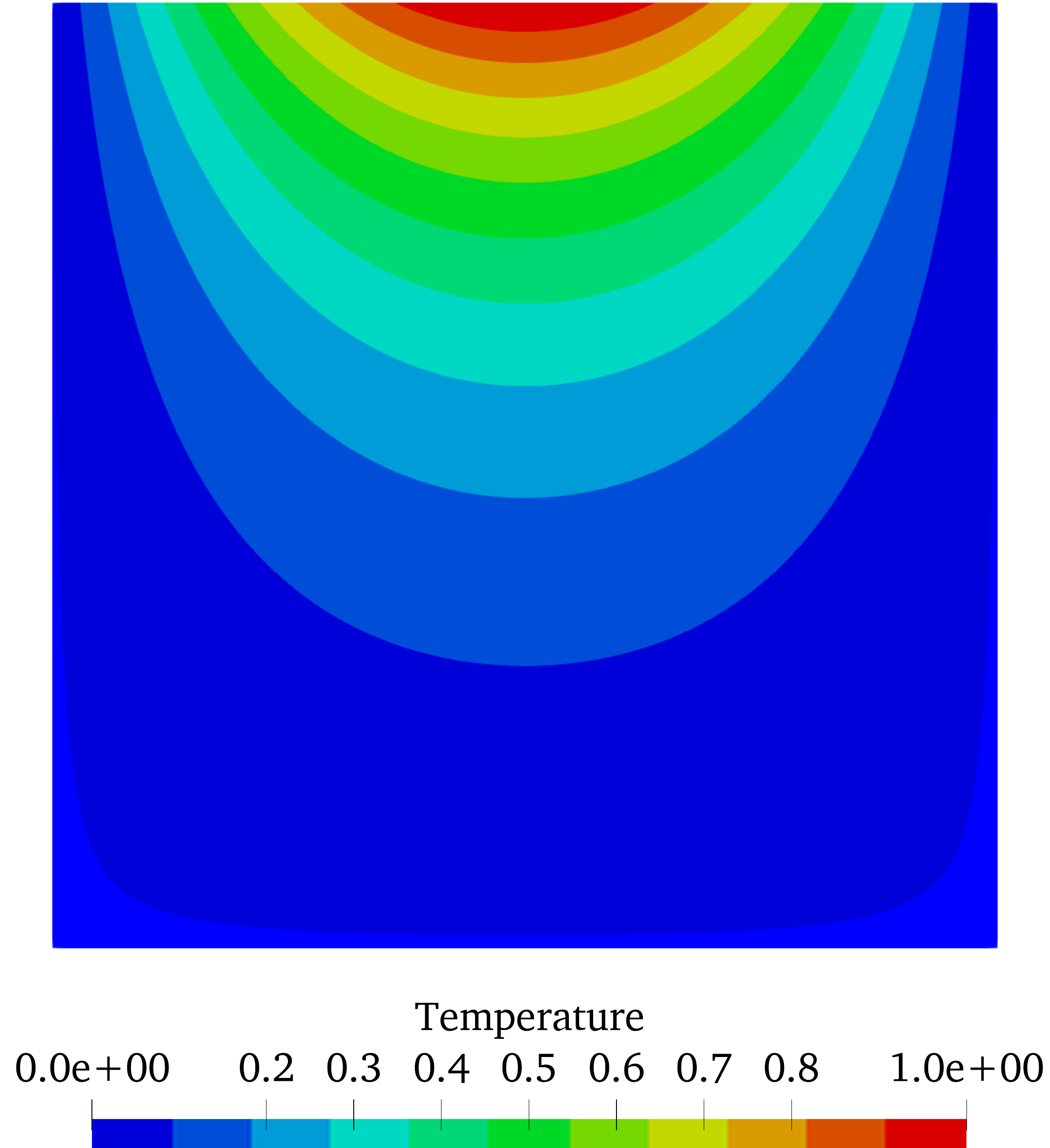}
  \caption{}
  \label{fig:1tri_T}
 \end{subfigure}
 \caption{Temperature contours for the square plate problem: (a) CS-FEM with polygonal elements; (b) FEM with quadrilateral elements; and (c) FEM with triangular elements.}
 \label{fig:temperature_contours_square_plate}
\end{figure}

Fig.~\ref{fig:temperature_contours_square_plate} shows the temperature contours obtained using the polygonal CS-FEM and conventional FEM with quadrilateral and triangular elements. It can be observed that the temperature fields predicted by the three discretization schemes are in excellent agreement, indicating that the proposed polygonal CS-FEM can accurately reproduce the analytical temperature distribution.

\begin{figure}[htbp]
  \centering
  \includegraphics[width=0.8\textwidth]{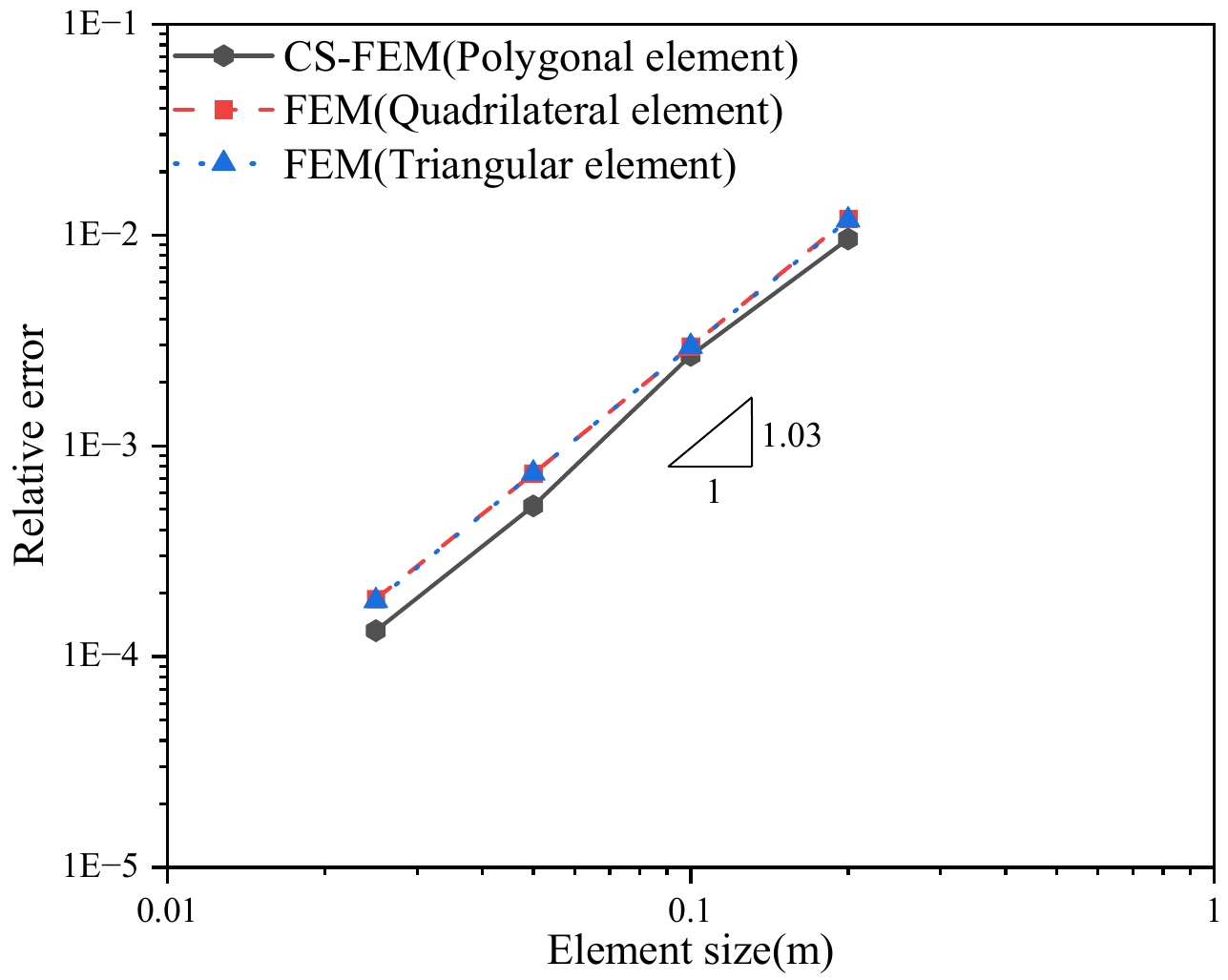}
  \caption{Convergence of the relative temperature error at point A.}
  \label{fig:1error} 
\end{figure}

The convergence behavior is further examined in Fig.~\ref{fig:1error}, where the relative temperature error at point A is plotted against the characteristic element size. Both FEM and CS-FEM exhibit stable convergence as the mesh is refined. Under the same element size, the polygonal CS-FEM produces slightly smaller errors than the conventional FEM models, demonstrating its improved accuracy for this benchmark problem.

\FloatBarrier

\FloatBarrier

To evaluate the influence of mesh distortion on the computational accuracy of the CS-FEM, nodal coordinate perturbations are introduced to the regular quadrilateral mesh in this study, with the specific expression \citep{litime2026} given as follows:
\begin{equation}
x_{ir} = x + a_{ir} \cdot r_c \cdot \Delta x, \quad y_{ir} = y + a_{ir} \cdot r_c \cdot \Delta y
\label{eq:perturbation}
\end{equation}
where $a_{ir}$ is a prescribed irregularity factor ranging from 0.0 to 0.5, $r_c$ is a random number uniformly distributed within $[-1.0, 1.0]$, and $\Delta x$ and $\Delta y$ denote the characteristic element spacings of the initial regular mesh in the $x$ and $y$ directions, respectively.

\begin{figure}[htbp]
\centering
\begin{subfigure}[b]{0.3\textwidth}
\centering
\includegraphics[width=1\textwidth]{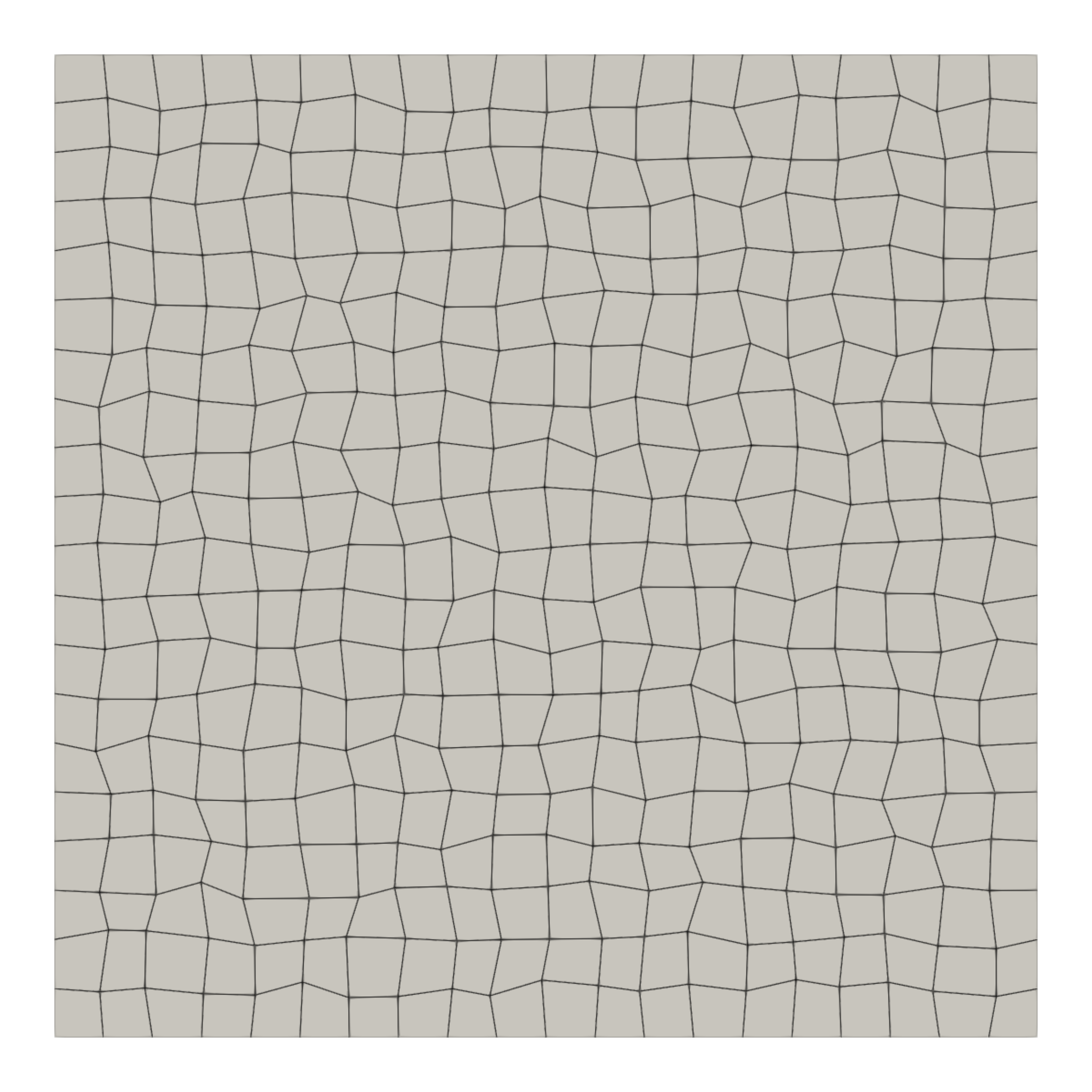}
\caption{}
\label{fig:a_02M}
\end{subfigure}
\hfill
\begin{subfigure}[b]{0.3\textwidth}
\centering
\includegraphics[width=1\textwidth]{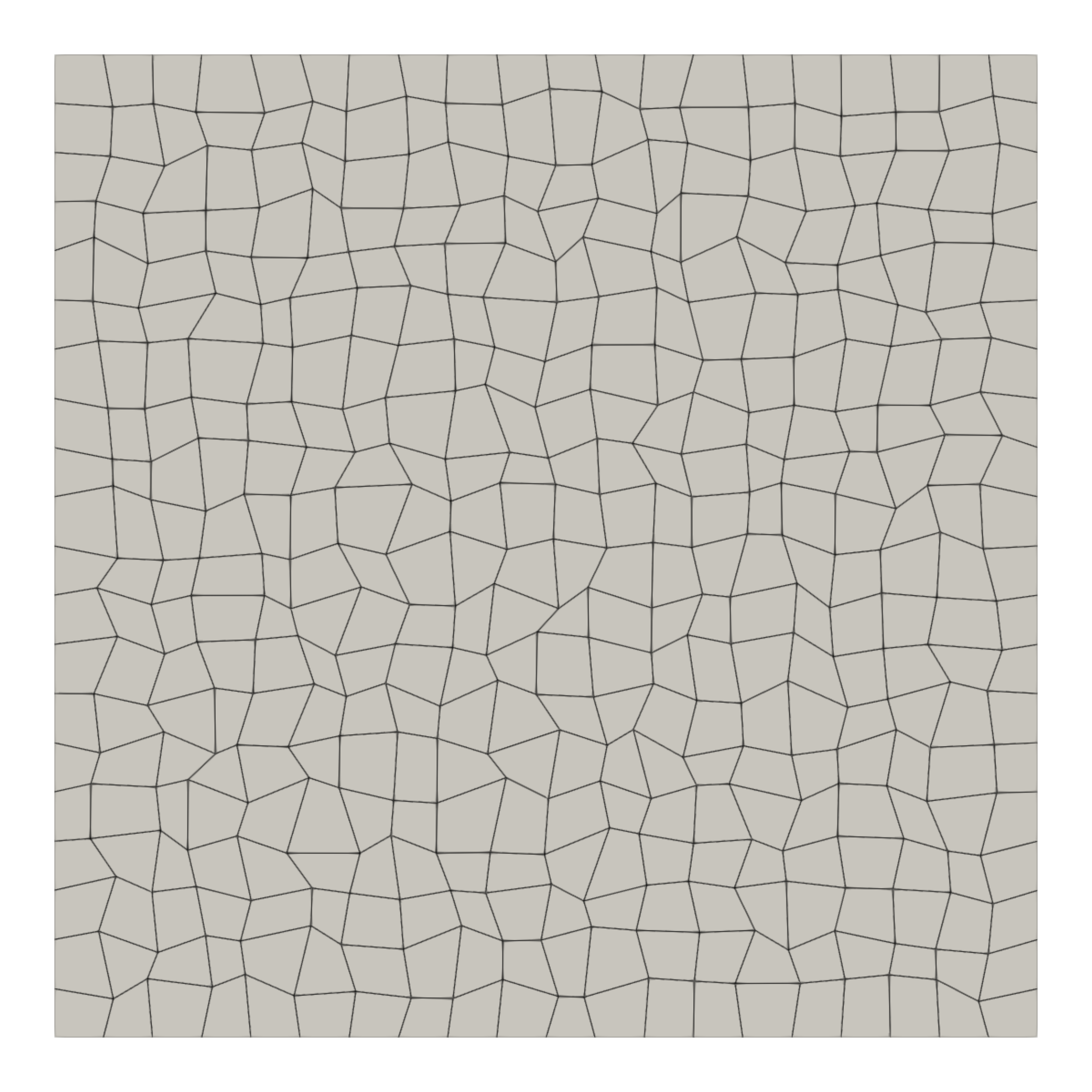}
\caption{}
\label{fig:a_03M}
\end{subfigure}
\hfill
\begin{subfigure}[b]{0.3\textwidth}
\centering
\includegraphics[width=1\textwidth]{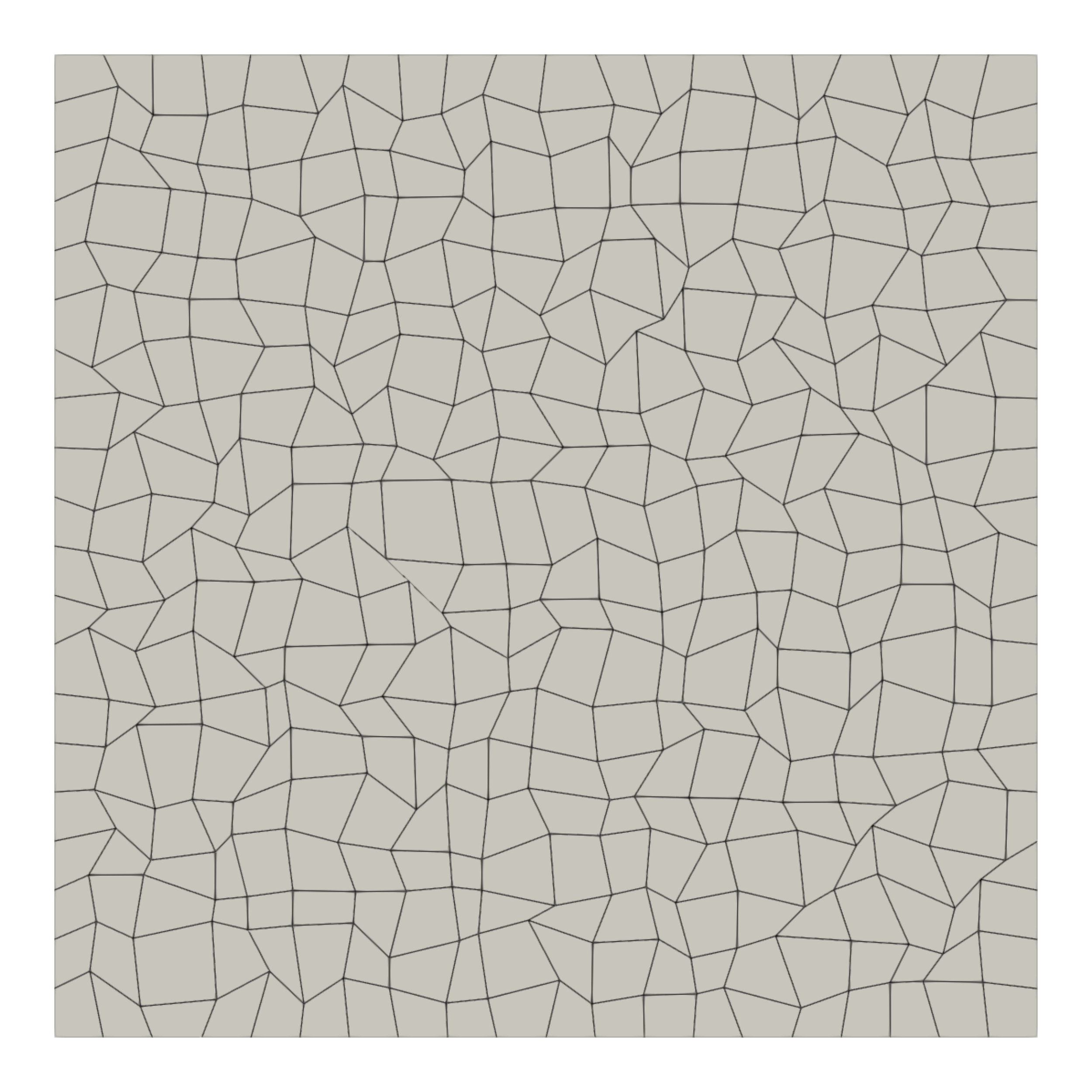}
\caption{}
\label{fig:a_04M}
\end{subfigure}
\caption{Mesh discretizations with varying irregularity factors: (a) $a_{ir}=0.2$; (b) $a_{ir}=0.3$; and (c) $a_{ir}=0.4$.}
\label{fig:Mesh discretizations with varying irregularity factors}
\end{figure}

Fig.~\ref{fig:Mesh discretizations with varying irregularity factors} presents the mesh configurations under irregularity factors of 0.2, 0.3, and 0.4, respectively.

\FloatBarrier

\begin{figure}[htbp]
\centering
\begin{subfigure}[b]{0.3\textwidth}
\centering
\includegraphics[width=1\textwidth]{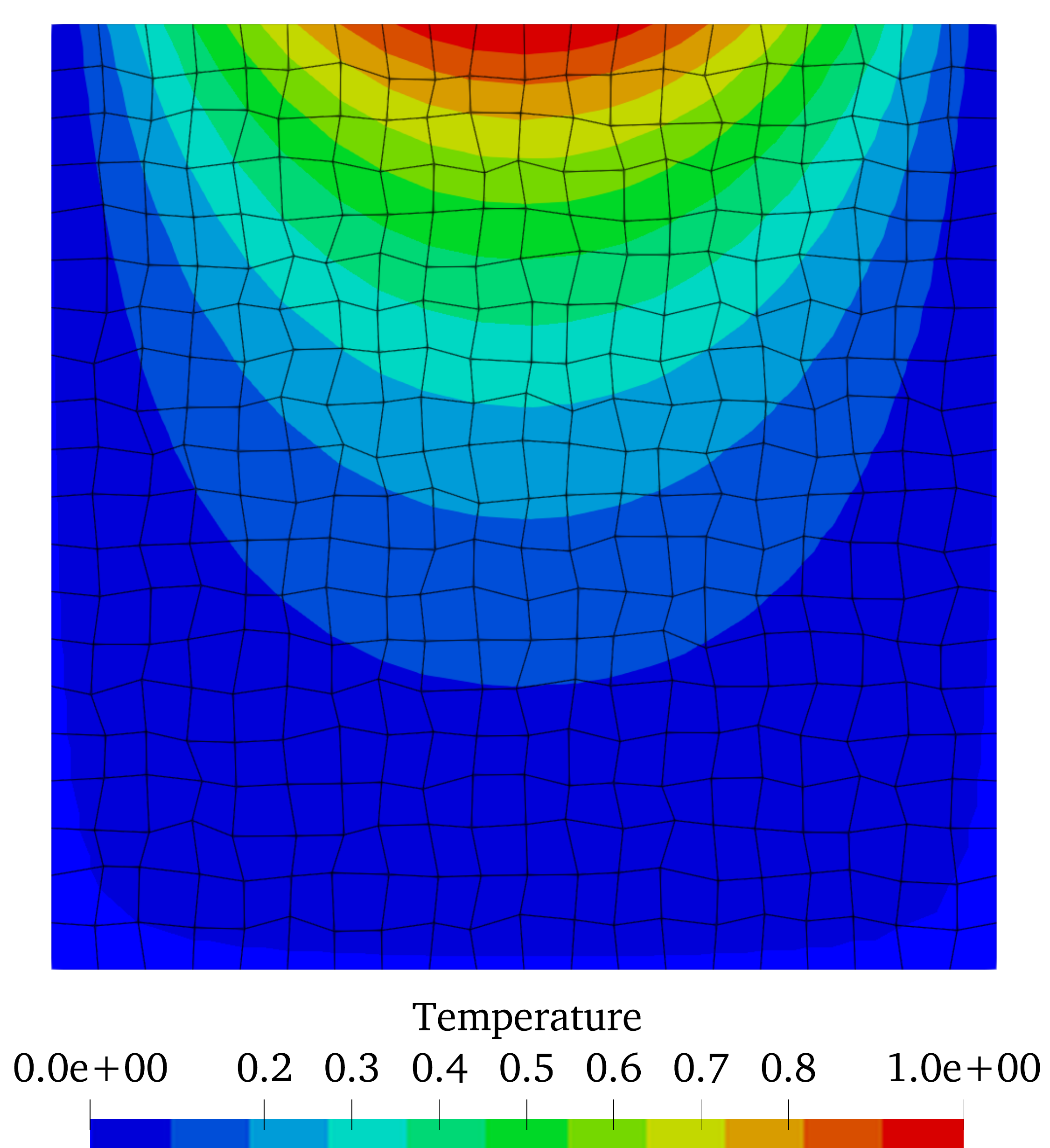}
\caption{}
\label{fig:a_02T}
\end{subfigure}
\hfill
\begin{subfigure}[b]{0.3\textwidth}
\centering
\includegraphics[width=1\textwidth]{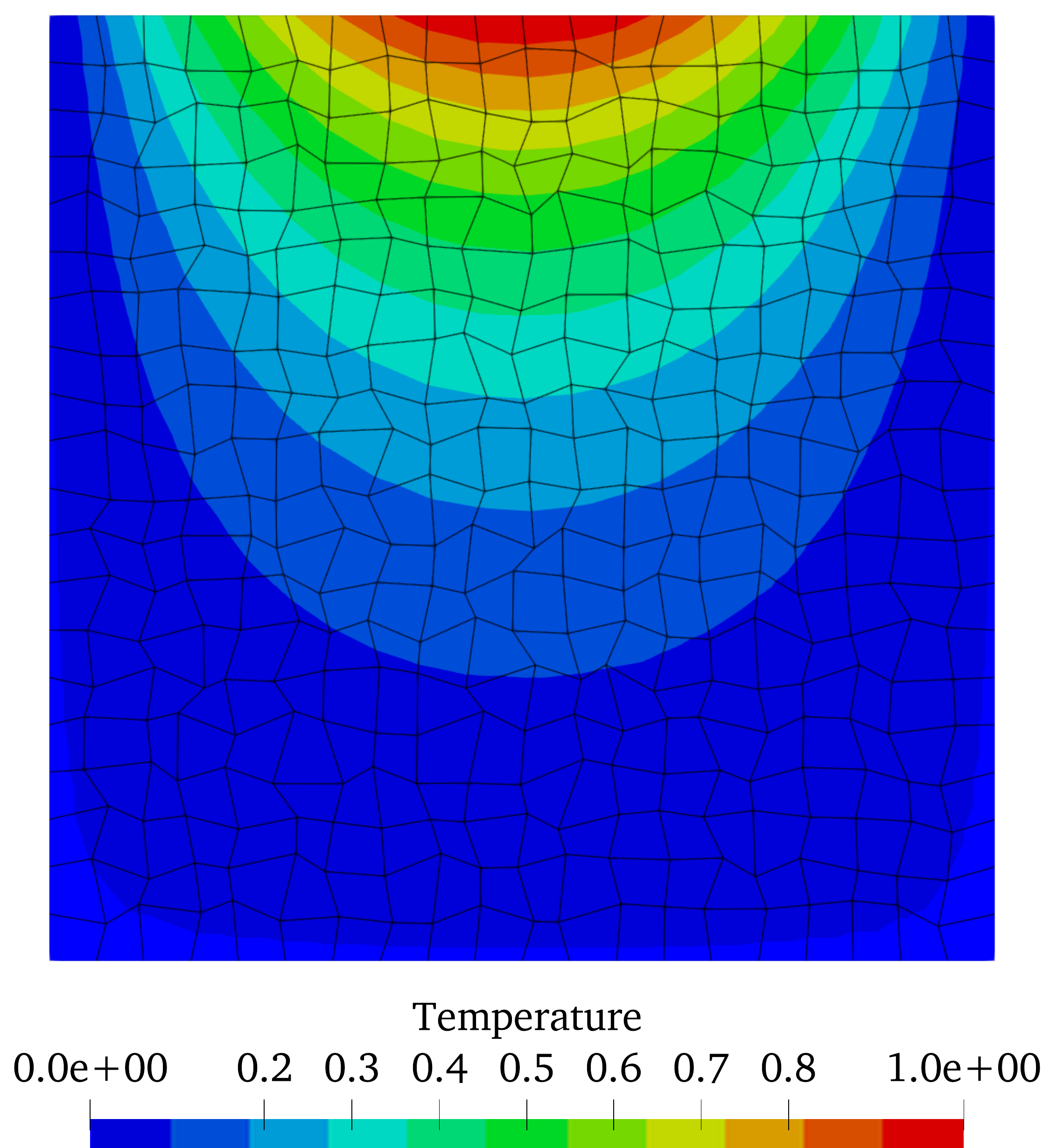}
\caption{}
\label{fig:a_03T}
\end{subfigure}
\hfill
\begin{subfigure}[b]{0.3\textwidth}
\centering
\includegraphics[width=1\textwidth]{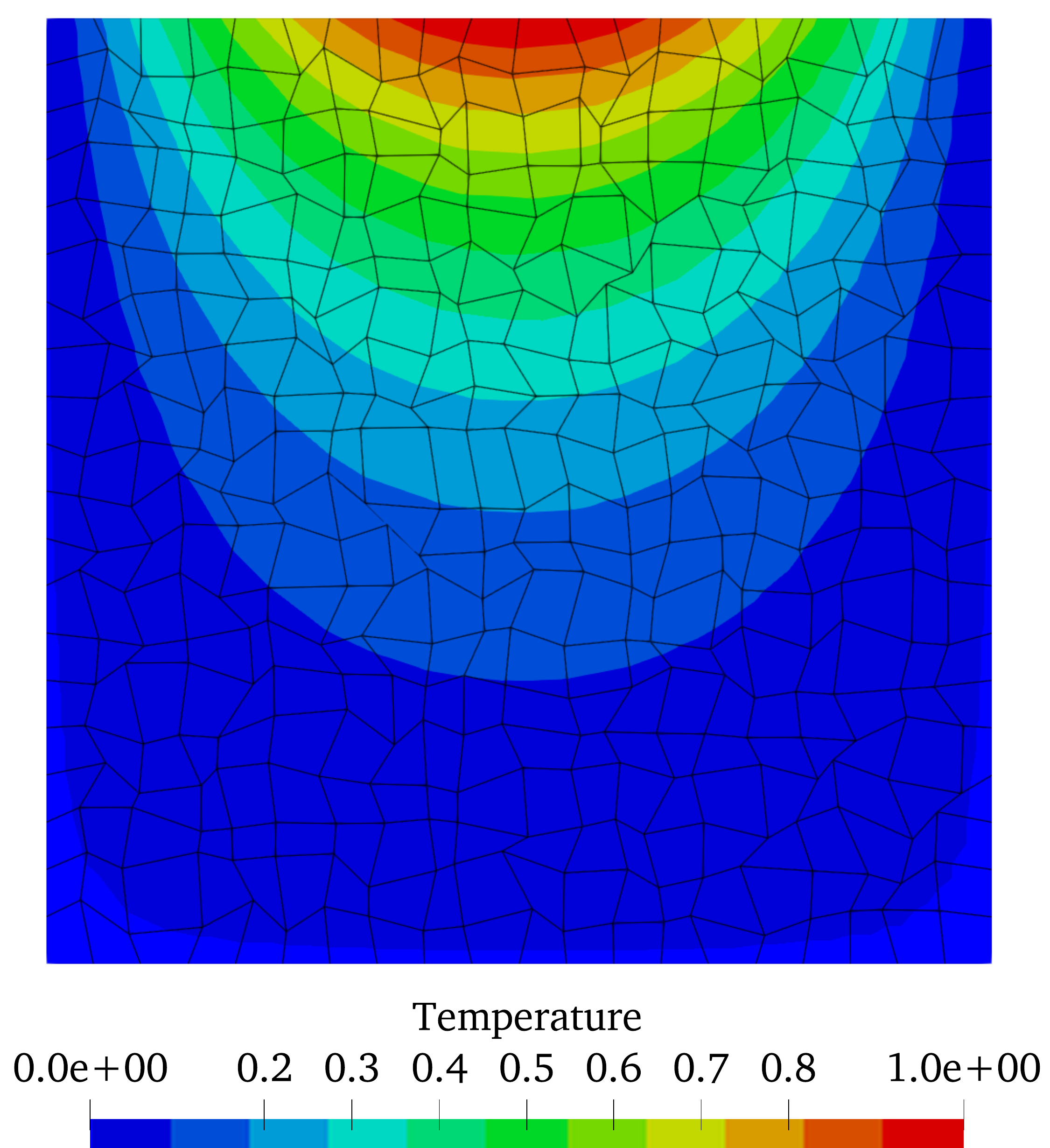}
\caption{}
\label{fig:a_04T}
\end{subfigure}
\caption{Temperature contours with varying irregularity factors:  (a) $a_{ir}=0.2$; (b) $a_{ir}=0.3$; and (c) $a_{ir}=0.4$.}
\label{fig:Temperature contours with varying irregularity factors}
\end{figure}

Fig.~\ref{fig:Temperature contours with varying irregularity factors} present the steady-state temperature contours computed by the CS-FEM under irregularity factors of 
 $a_{ir}=0.2,0.3$ and $0.4$,respectively. A qualitative comparison reveals that the overall temperature distributions,isotherm patterns,and locations of hot and cold regions remain highly consistent across all three cases. Together with the relative error data shown in Tab.~\ref{tab:error}, As the irregularity factor increases from 0.2 to 0.4, the relative error rises merely from 0.38\% to 0.44\%. This indicates that the CS-FEM maintains high computational accuracy even when the mesh is severely distorted, demonstrating its excellent robustness against mesh deformation.

\FloatBarrier

\begin{table}[htbp]
\centering
\caption{Relative errors under different irregularity factors.}
\label{tab:error}
\begin{tabular}{p{3.5cm} c c c}
\toprule
$\alpha_{ir}$ & 0.2 & 0.3 & 0.4 \\
\midrule
Relative error (\%) & 0.38 & 0.41 & 0.44 \\
\bottomrule
\end{tabular}
\end{table}

\subsubsection{Steady-state heat conduction in a square plate with multiple circular holes}

\begin{figure}[htbp]
 \centering
 \begin{subfigure}[b]{0.34\textwidth}
  \centering
  \includegraphics[width=1\textwidth]{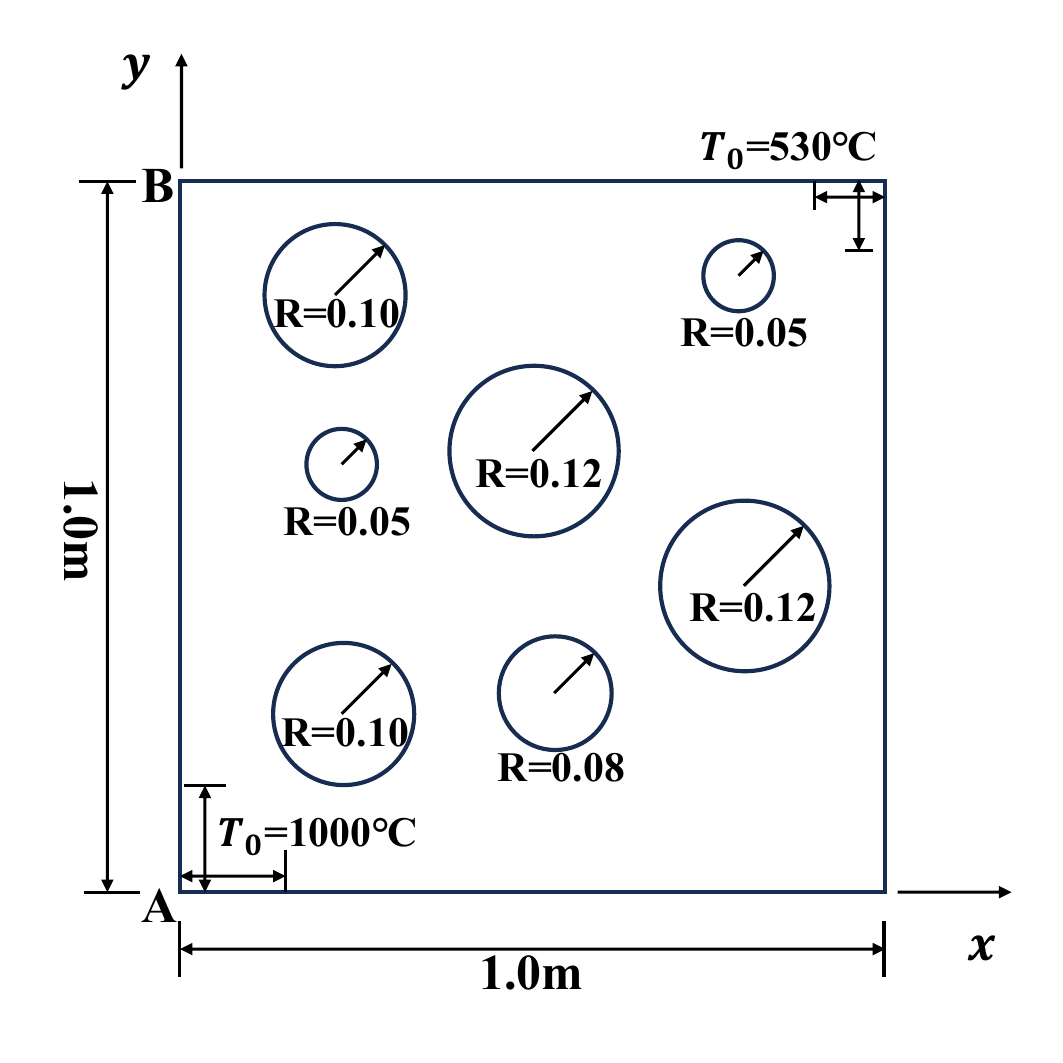}
  \caption{}
  \label{fig:geometry2}
 \end{subfigure}
 \hfill
 \begin{subfigure}[b]{0.3\textwidth}
  \centering
  \includegraphics[width=1\textwidth]{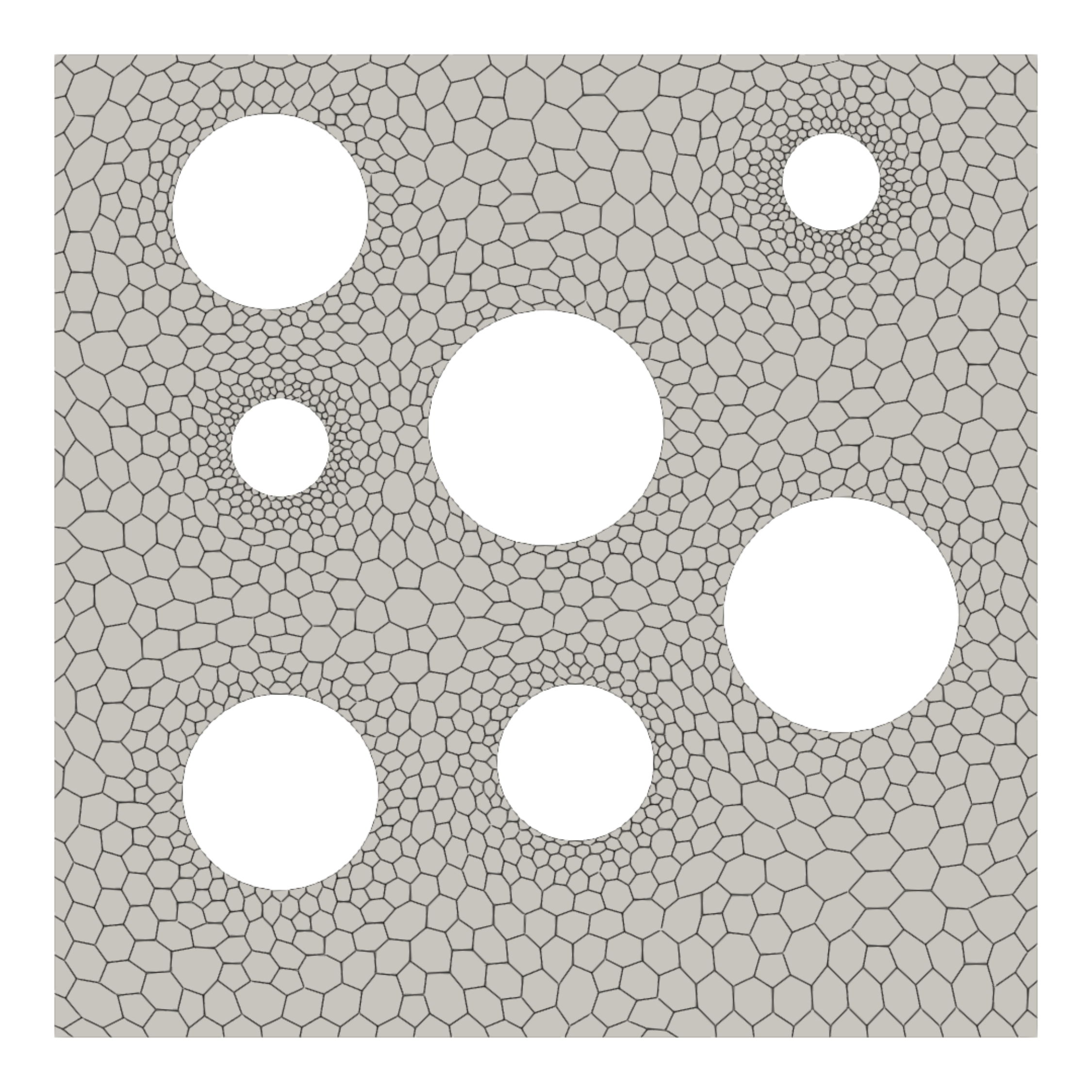}
  \caption{}
  \label{fig:2poly_M}
 \end{subfigure}
 \hfill
 \begin{subfigure}[b]{0.3\textwidth}
  \centering
  \includegraphics[width=1\textwidth]{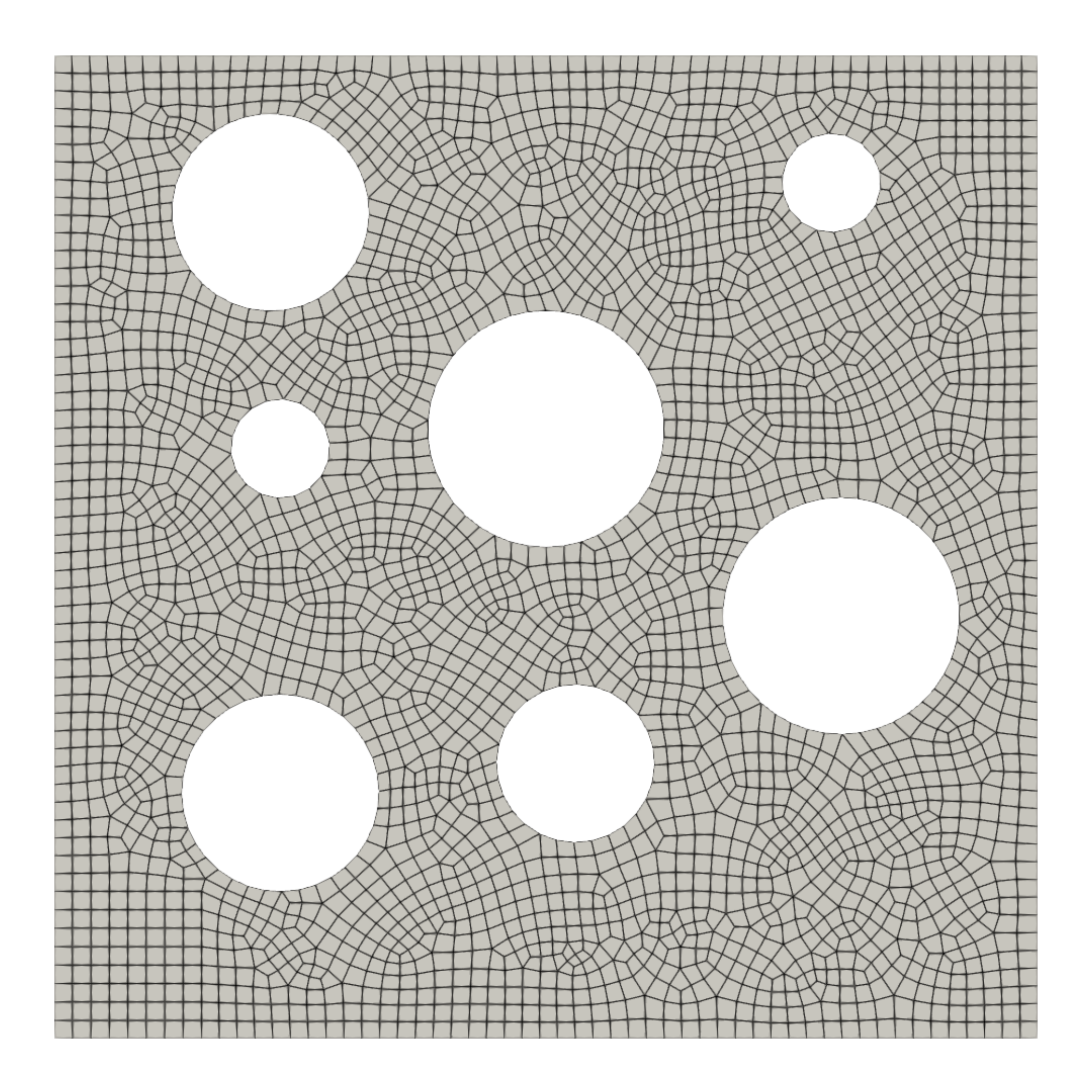}
  \caption{}
  \label{fig:2quad_M}
 \end{subfigure}
 \caption{Computational model of the square plate with multiple circular holes: (a) geometry and boundary conditions; (b) polygonal mesh; and (c) quadrilateral mesh.}
 \label{fig:The_model_of_the_square_plate_with_multiple_holes}
\end{figure}

This example considers steady-state heat conduction in a square plate containing seven circular holes with different radii, as illustrated in Fig.~\ref{fig:The_model_of_the_square_plate_with_multiple_holes}. The side length of the square plate is $L=1~\mathrm{m}$. A prescribed temperature of $T_0=1000~^{\circ}\mathrm{C}$ is applied to a boundary segment of length $0.15~\mathrm{m}$ near the lower-left corner, whereas a prescribed temperature of $T_1=530~^{\circ}\mathrm{C}$ is imposed on a boundary segment of length $0.1~\mathrm{m}$ near the upper-right corner. The remaining boundaries are treated as thermally insulated unless otherwise specified.

The computational model is discretized using polygonal elements with 3226 nodes. For comparison, a conventional FEM model using quadrilateral elements with 3245 nodes is also established, and its numerical solution is taken as the reference solution.

\FloatBarrier

\begin{figure}[htbp]
 \centering
 \begin{subfigure}[b]{0.45\textwidth}
  \centering
  \includegraphics[width=1\textwidth]{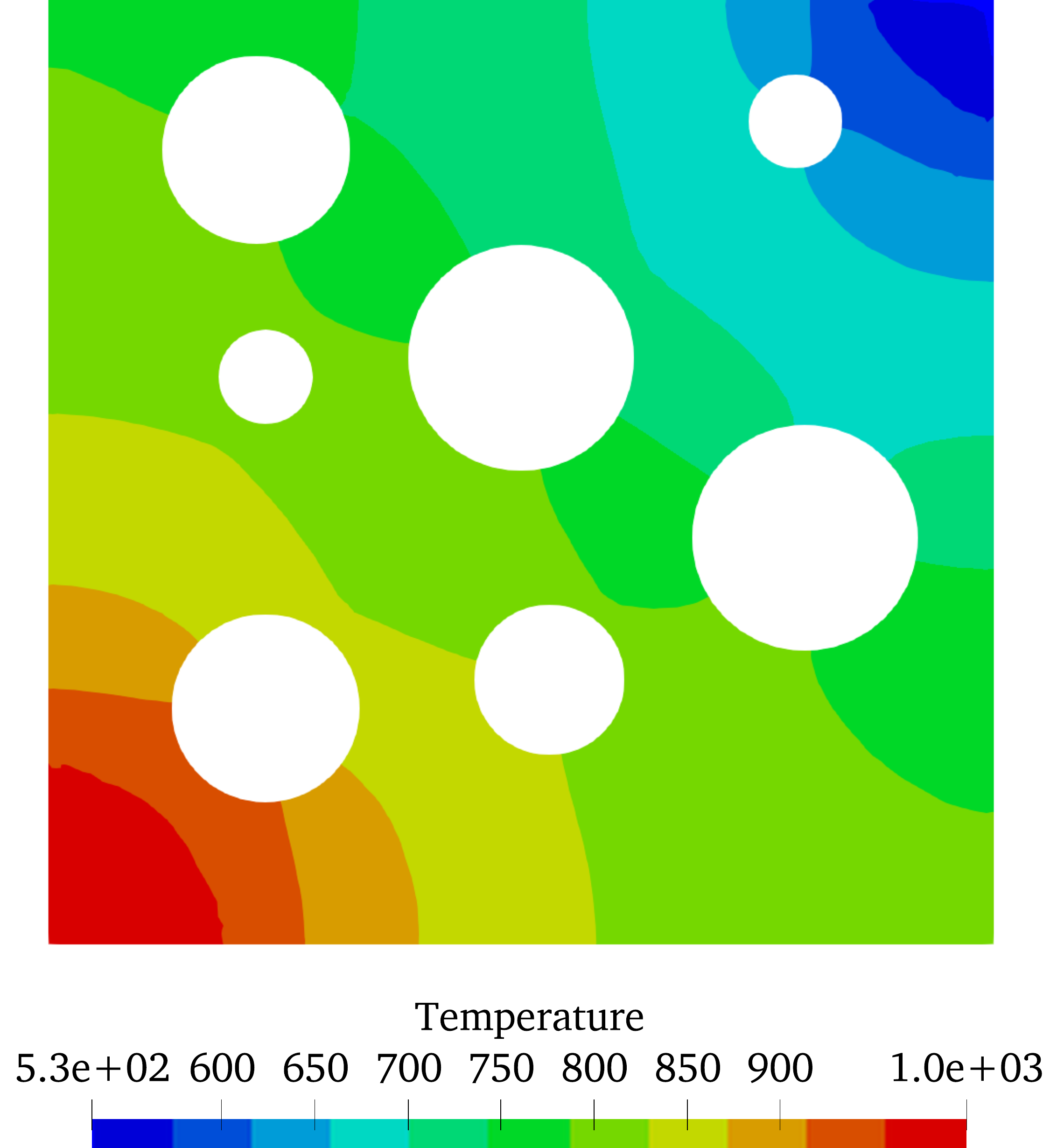}
  \caption{}
  \label{fig:2poly_T}
 \end{subfigure}
 \hfill
 \begin{subfigure}[b]{0.45\textwidth}
  \centering
  \includegraphics[width=1\textwidth]{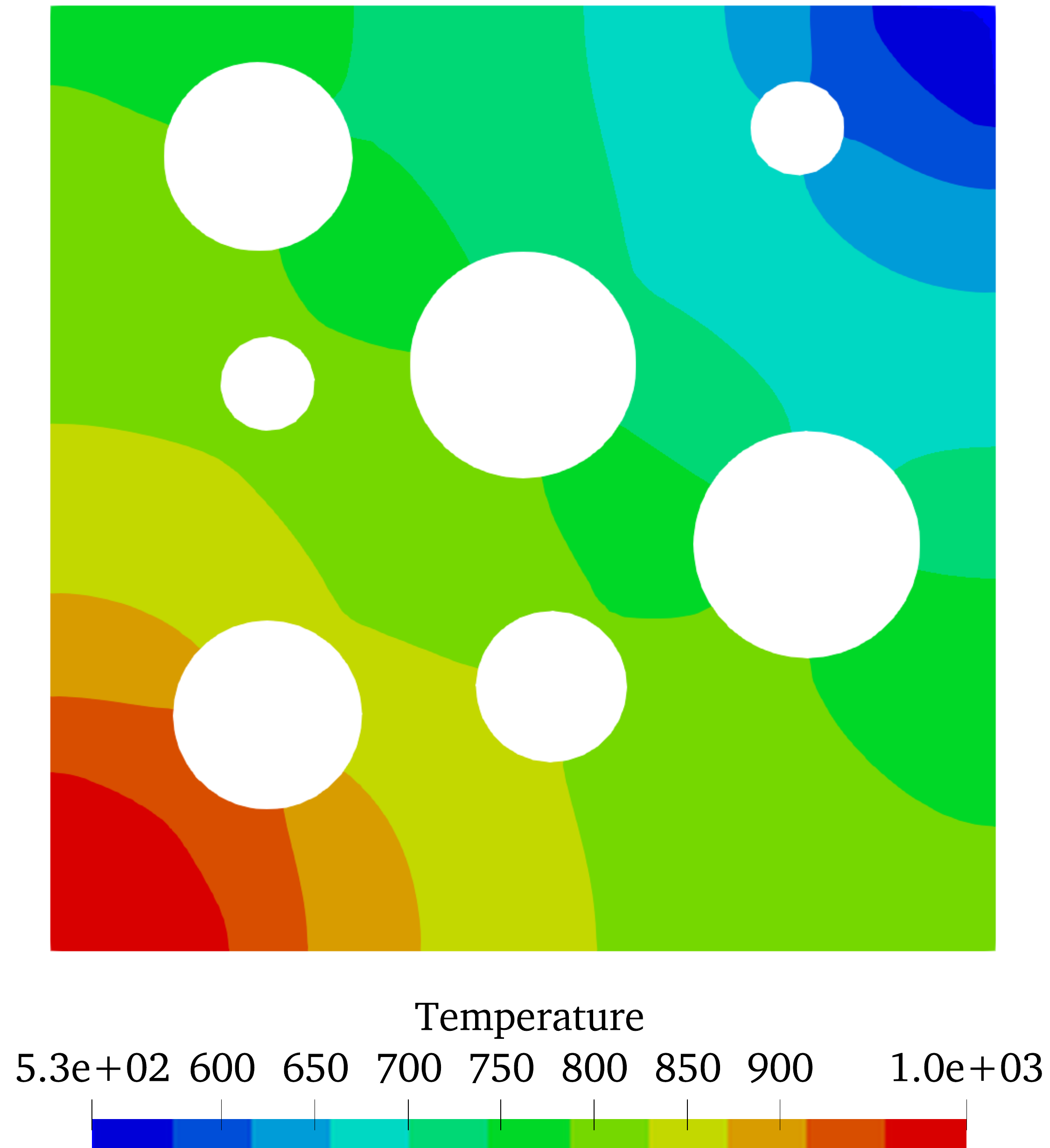}
  \caption{}
  \label{fig:2quad_T}
 \end{subfigure}
 \caption{Temperature contours for the square plate with multiple circular holes: (a) polygonal CS-FEM; and (b) FEM reference solution.}
 \label{fig:Contour_plots_of_temperature}
\end{figure}

To further assess the accuracy of the proposed method, the temperature distribution and the temperature variation along line AB are compared with the reference solution. Fig.~\ref{fig:Contour_plots_of_temperature} presents the temperature contours obtained using the polygonal CS-FEM and the conventional FEM. The two temperature fields show very similar distributions, indicating that the proposed method can accurately capture the global heat-transfer behavior in a multiply connected domain.

\begin{figure}[htbp]
  \centering
  \includegraphics[width=0.8\textwidth]{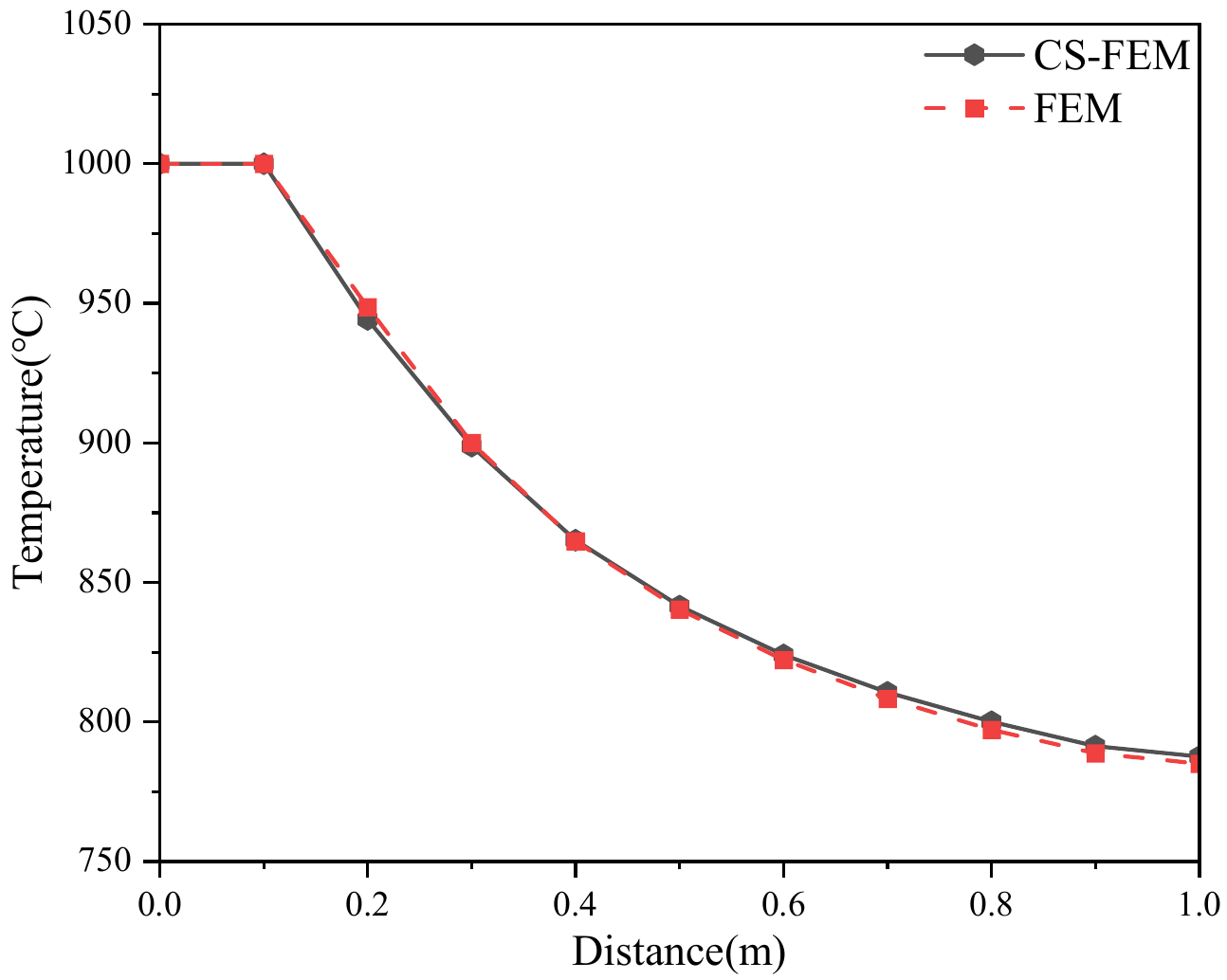}
  \caption{Comparison of temperature profiles along line AB.}
  \label{fig:ABT} 
\end{figure}

The temperature profiles along line AB are further compared in Fig.~\ref{fig:ABT}. It can be observed that the results obtained by the polygonal CS-FEM agree well with the reference solution throughout the selected line. These comparisons demonstrate the accuracy and applicability of the proposed polygonal CS-FEM for steady-state heat-conduction problems involving complex geometries.

\FloatBarrier

\FloatBarrier

\subsection{Transient heat-conduction analysis}

\subsubsection{A square plate without heat source}

\begin{figure}[htbp]
 \centering
 \begin{subfigure}[b]{0.45\textwidth}
  \centering
  \includegraphics[width=0.95\textwidth]{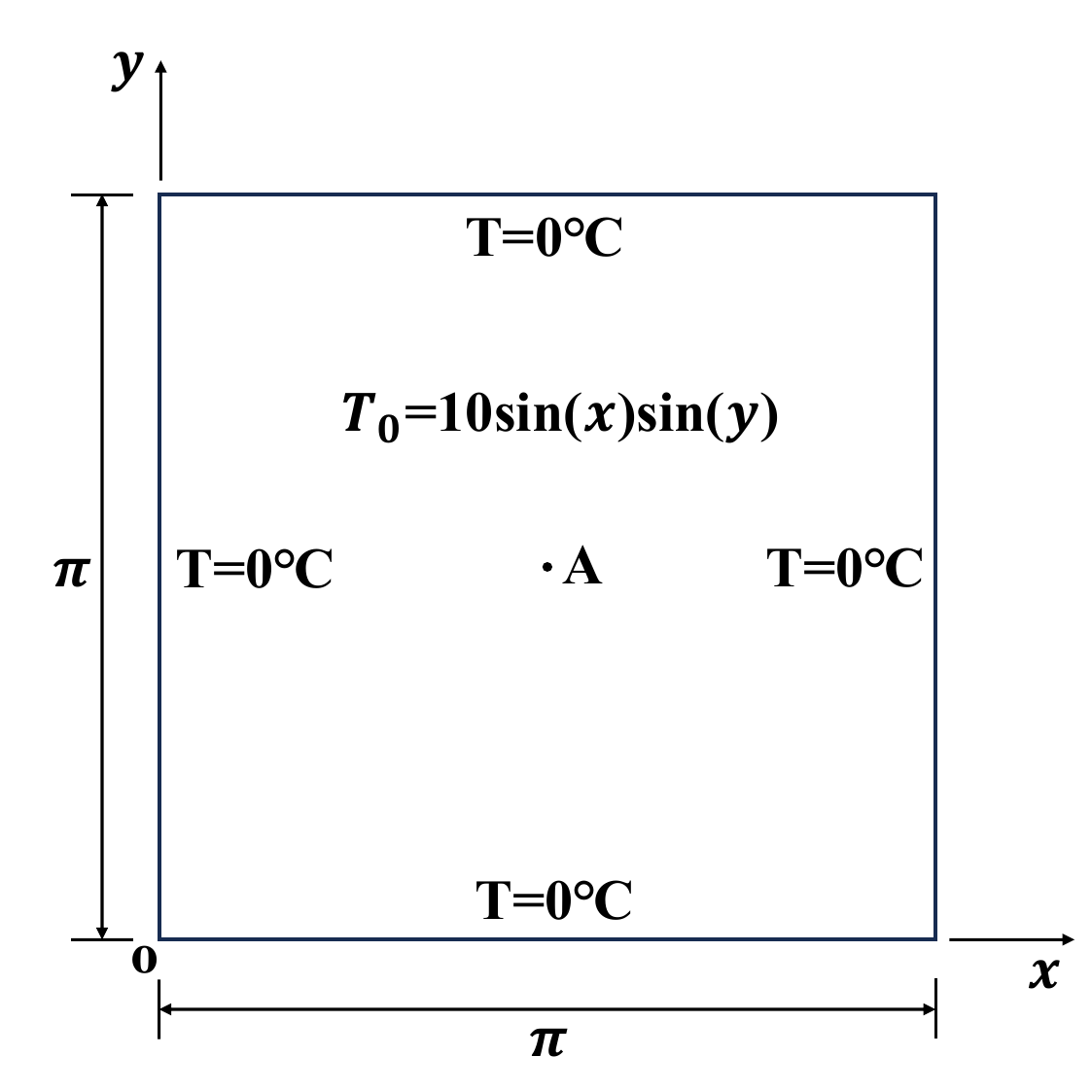}
  \caption{}
  \label{fig:geometry3}
 \end{subfigure}
 \hfill
 \begin{subfigure}[b]{0.45\textwidth}
  \centering
  \includegraphics[width=0.95\textwidth]{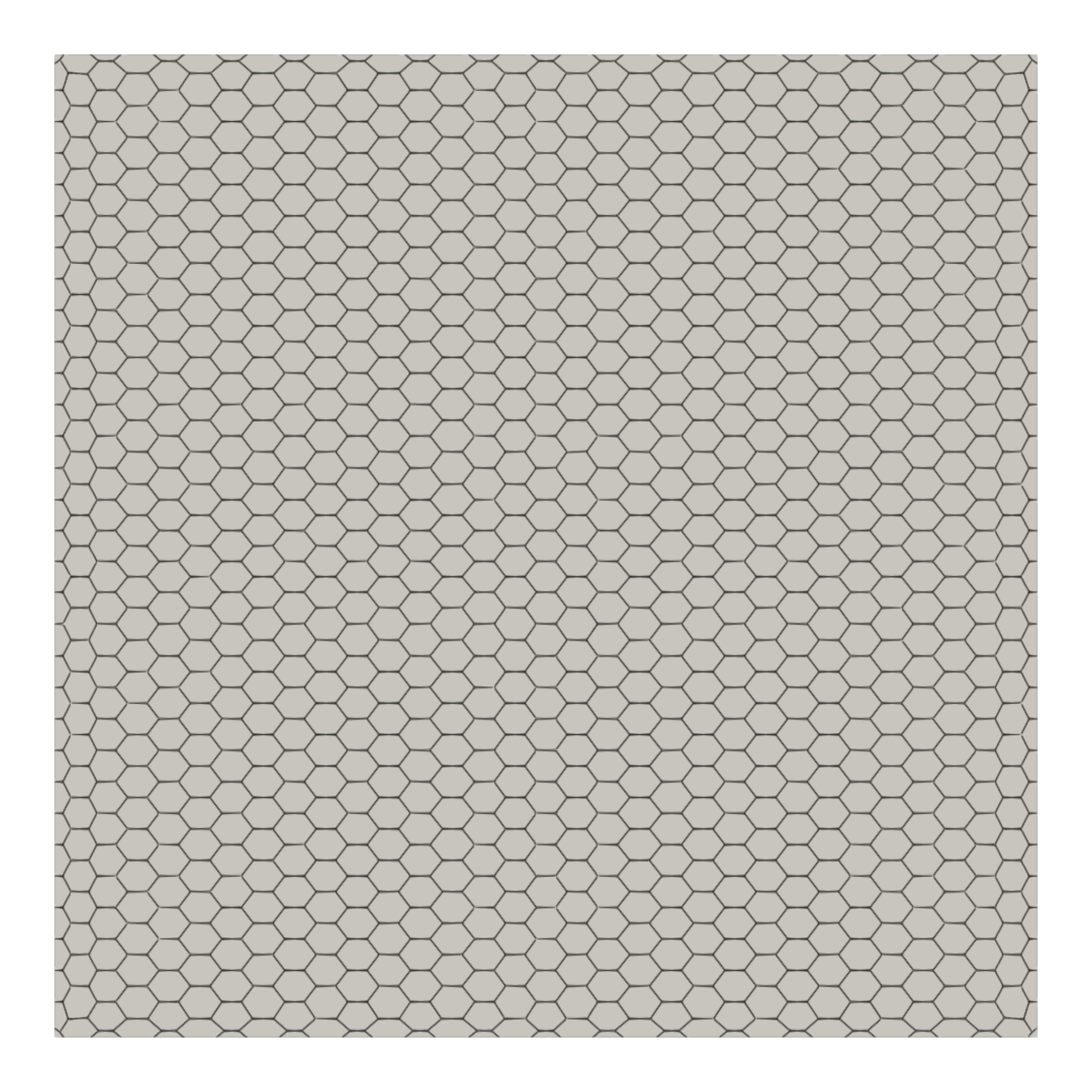}
  \caption{}
  \label{fig:3poly_M}
 \end{subfigure}
 \vspace{0cm}
 \begin{subfigure}[b]{0.45\textwidth}
  \centering
  \includegraphics[width=0.95\textwidth]{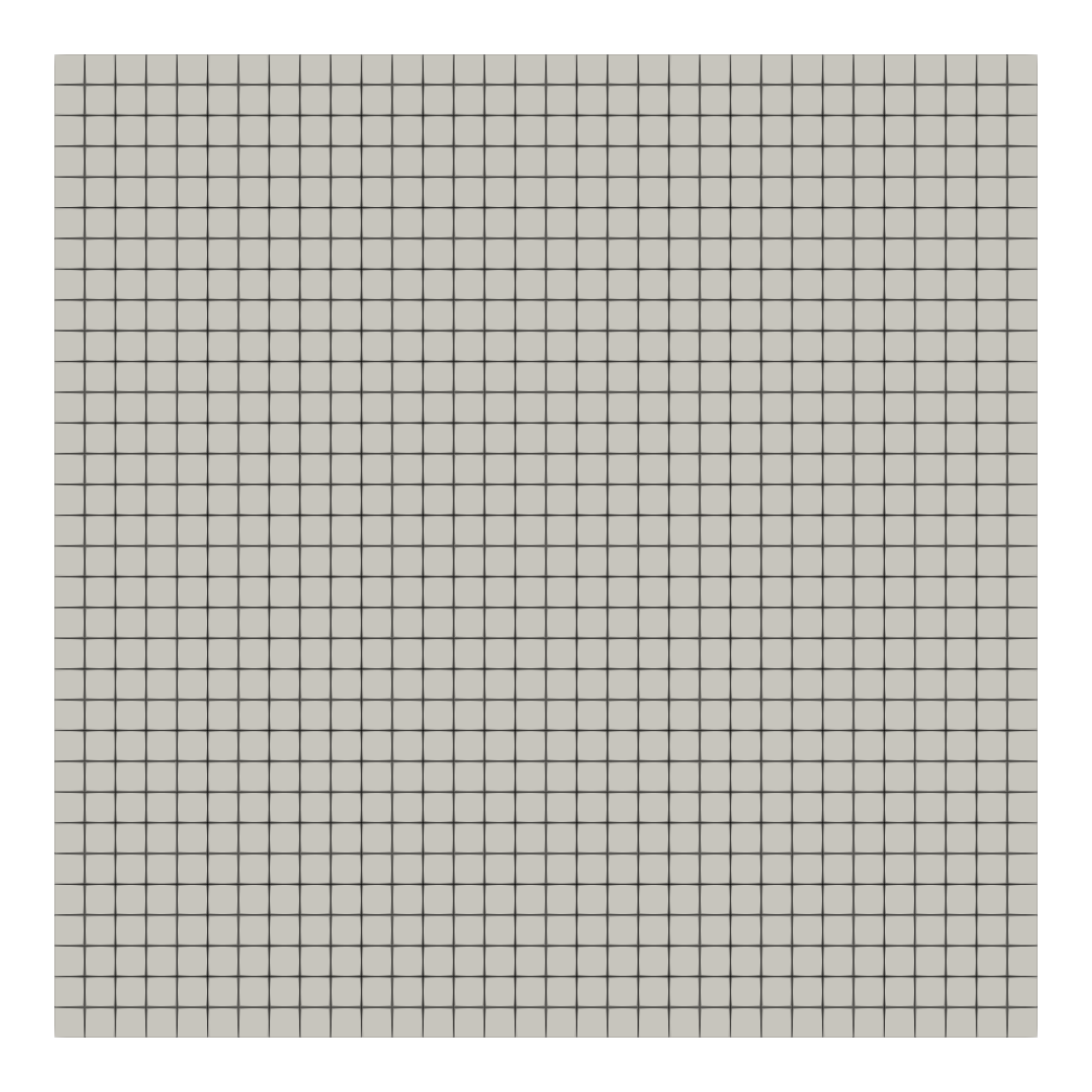}
  \caption{}
  \label{fig:3quad_M}
 \end{subfigure}
 \hfill
 \begin{subfigure}[b]{0.45\textwidth}
  \centering
  \includegraphics[width=0.95\textwidth]{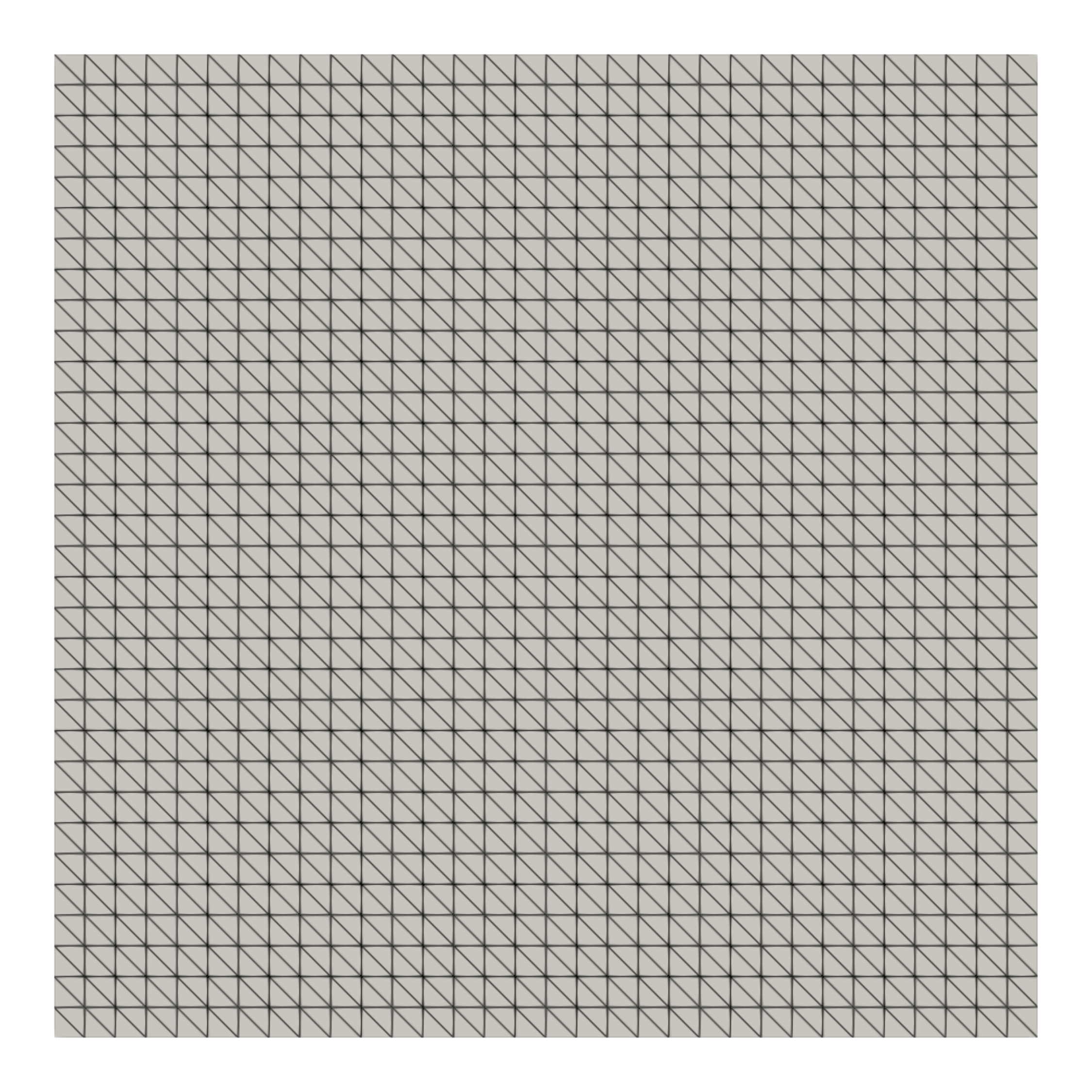}
  \caption{}
  \label{fig:3tri_M}
 \end{subfigure}
 \caption{Computational model of the square plate without heat source: (a) geometry and boundary conditions; (b) polygonal mesh; (c) quadrilateral mesh; and (d) triangular mesh.}
 \label{fig:The_model_of_the_square_plate_without_heat_source}
\end{figure}

This example is used to verify the accuracy of the proposed polygonal CS-FEM for transient heat-conduction analysis. A square domain with dimensions of $\pi \times \pi$ is considered, as shown in Fig.~\ref{fig:The_model_of_the_square_plate_without_heat_source}. Zero-temperature Dirichlet boundary conditions are prescribed on all sides of the square domain. The initial temperature field is given by
\begin{equation}
T(x,y,0)=10\sin(x)\sin(y).
\end{equation}
For this problem, the analytical solution \cite{lin_transient_2017} is expressed as
\begin{equation}
T(x,y,t)=10e^{-2t}\sin(x)\sin(y).
\end{equation}

\FloatBarrier

\begin{figure}[htbp]
 \centering
 \begin{subfigure}[b]{0.32\textwidth}
  \centering
  \includegraphics[width=1\textwidth]{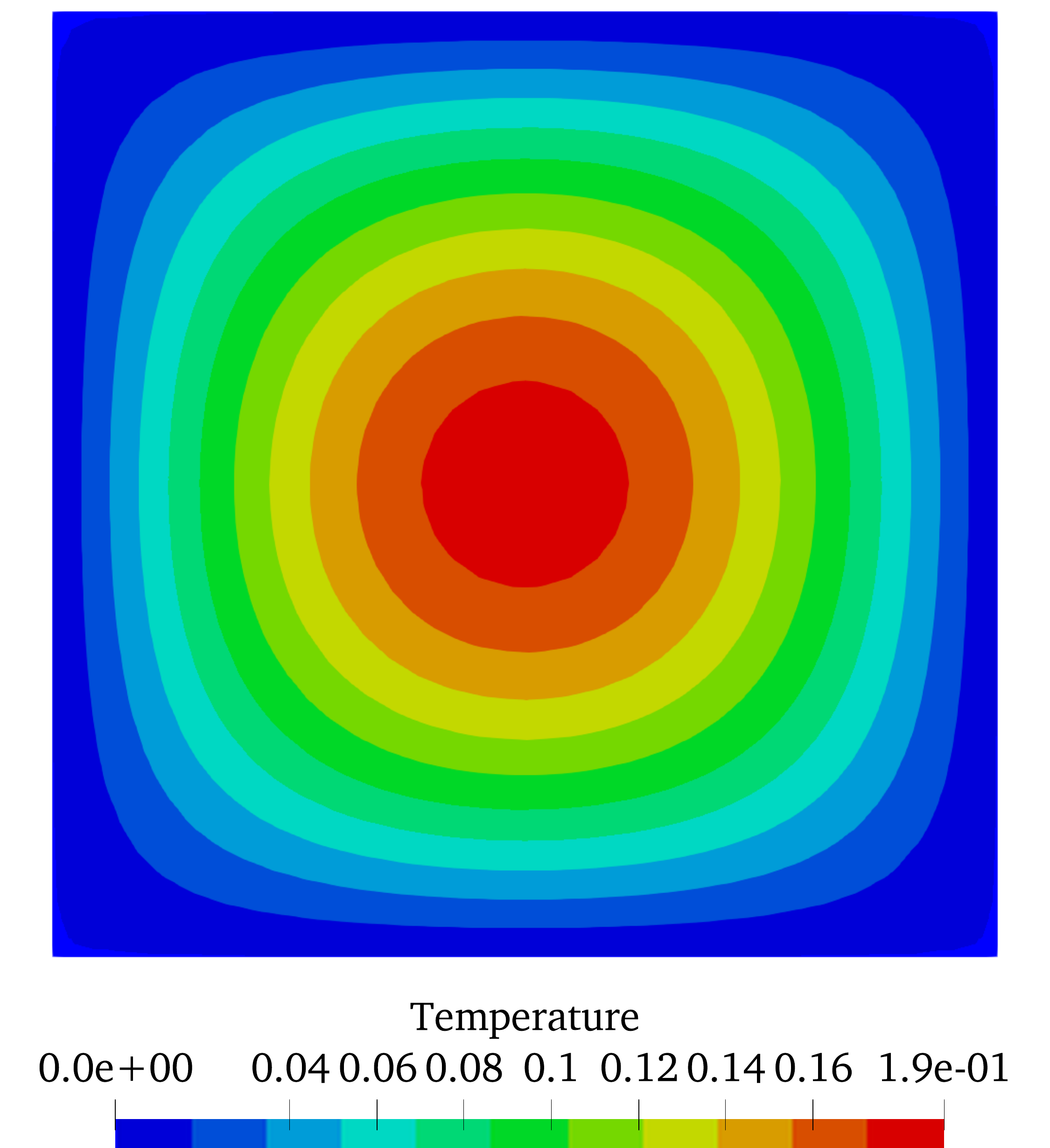}
  \caption{}
  \label{fig:3poly_T}
 \end{subfigure}
 \hfill
 \begin{subfigure}[b]{0.32\textwidth}
  \centering
  \includegraphics[width=1\textwidth]{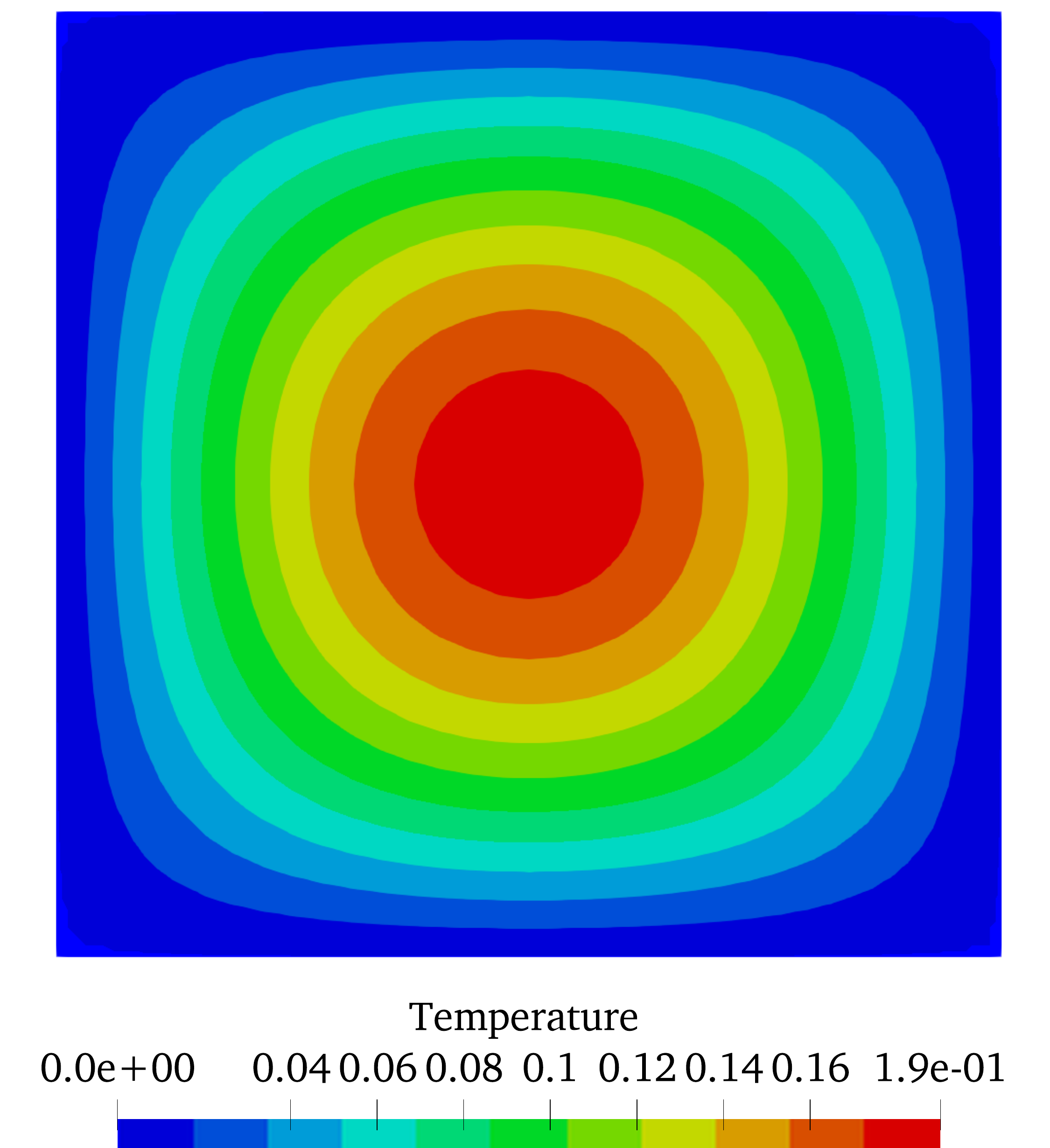}
  \caption{}
  \label{fig:3quad_T}
 \end{subfigure}
 \hfill
 \begin{subfigure}[b]{0.32\textwidth}
  \centering
  \includegraphics[width=1\textwidth]{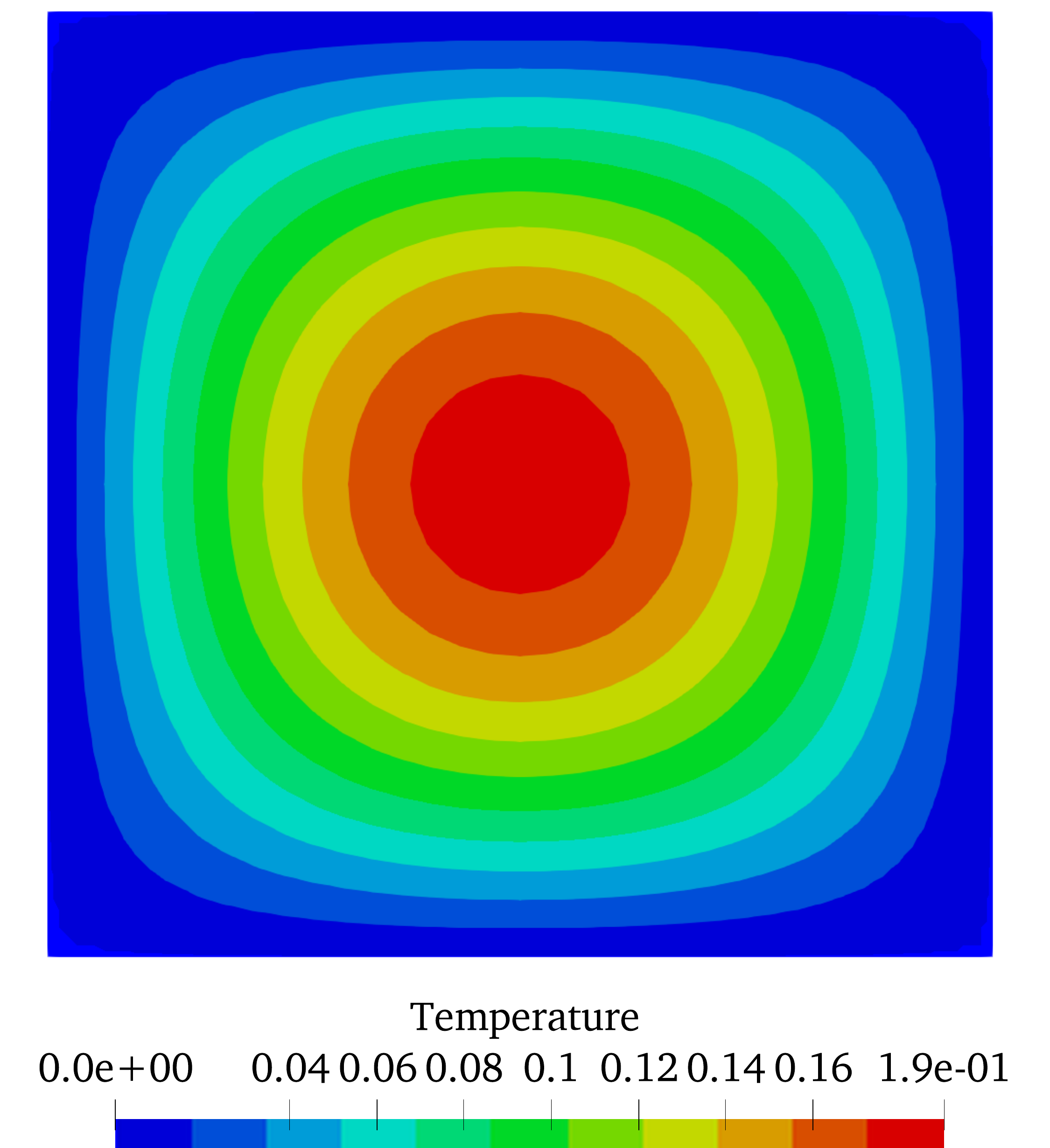}
  \caption{}
  \label{fig:3tri_T}
 \end{subfigure}
 \caption{Temperature distributions at $t=2~\mathrm{s}$: (a) CS-FEM with polygonal elements; (b) FEM with quadrilateral elements; and (c) FEM with triangular elements.}
 \label{fig:The_temperature_distribution}
\end{figure}

In the numerical analysis, a time step of $\Delta t=0.001~\mathrm{s}$ is used, and the total simulation time is set to $t=2~\mathrm{s}$. Fig.~\ref{fig:The_temperature_distribution} shows the temperature contours at $t=2~\mathrm{s}$ obtained using the polygonal CS-FEM and conventional FEM with quadrilateral and triangular elements. Very good agreement can be observed among the three numerical results.

\begin{figure}[htbp]
  \centering
  \includegraphics[width=0.75\textwidth]{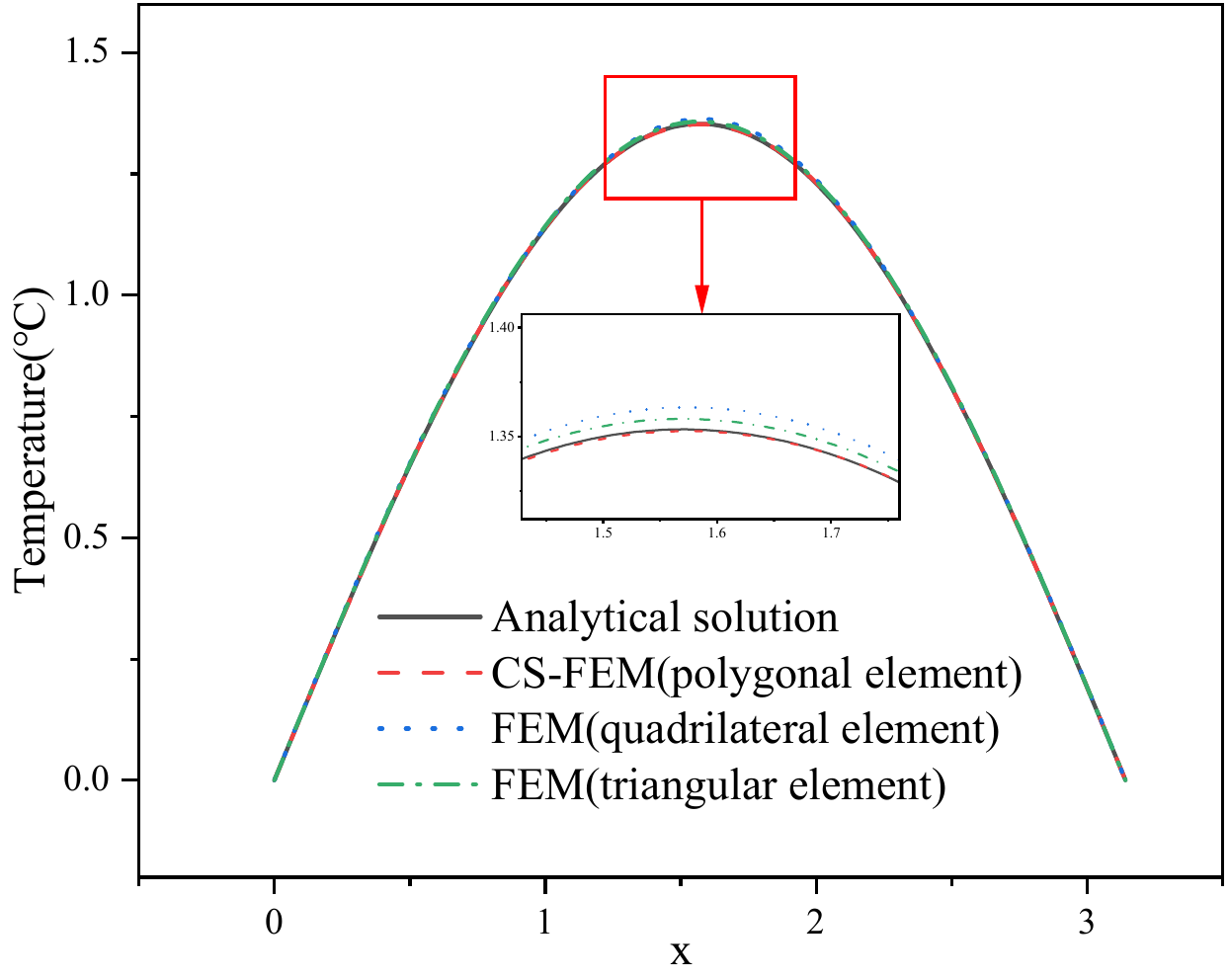}
  \caption{Comparison of temperature distributions along $y=\pi/2$ at $t=1~\mathrm{s}$.}
  \label{fig:xT} 
\end{figure}

The temperature variation along the line $y=\pi/2$ at $t=1~\mathrm{s}$ is further compared with the analytical solution, as shown in Fig.~\ref{fig:xT}. The temperature profile obtained by the polygonal CS-FEM agrees closely with the analytical solution, confirming the capability of the proposed method to accurately model transient heat conduction.

\FloatBarrier

\FloatBarrier

\begin{figure}[htbp]
 \centering
 \begin{subfigure}[b]{0.4\textwidth}
  \centering
  \includegraphics[width=1\textwidth]{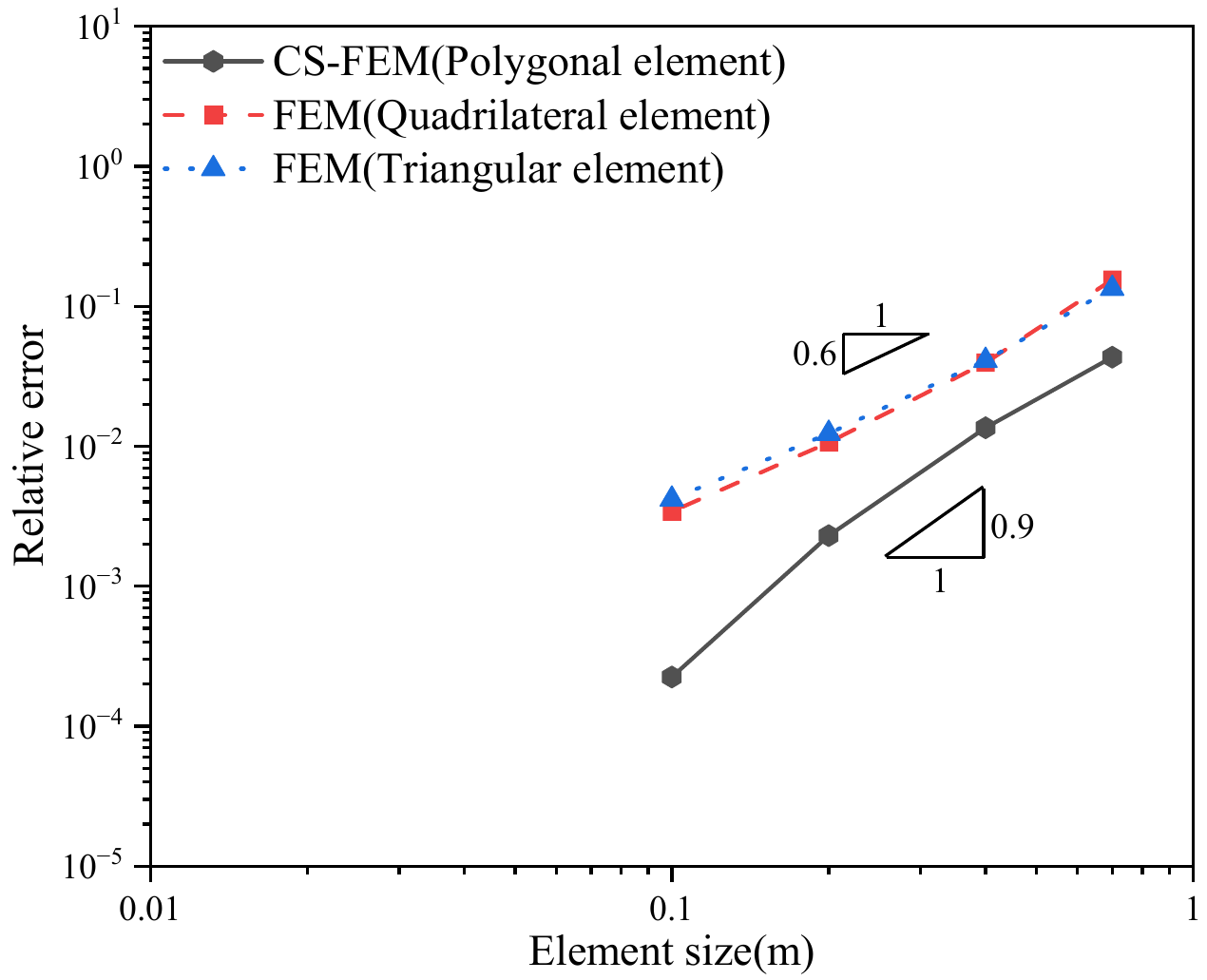}
  \caption{}
  \label{fig:0.5serror}
 \end{subfigure}
 \hfill
 \begin{subfigure}[b]{0.4\textwidth}
  \centering
  \includegraphics[width=1\textwidth]{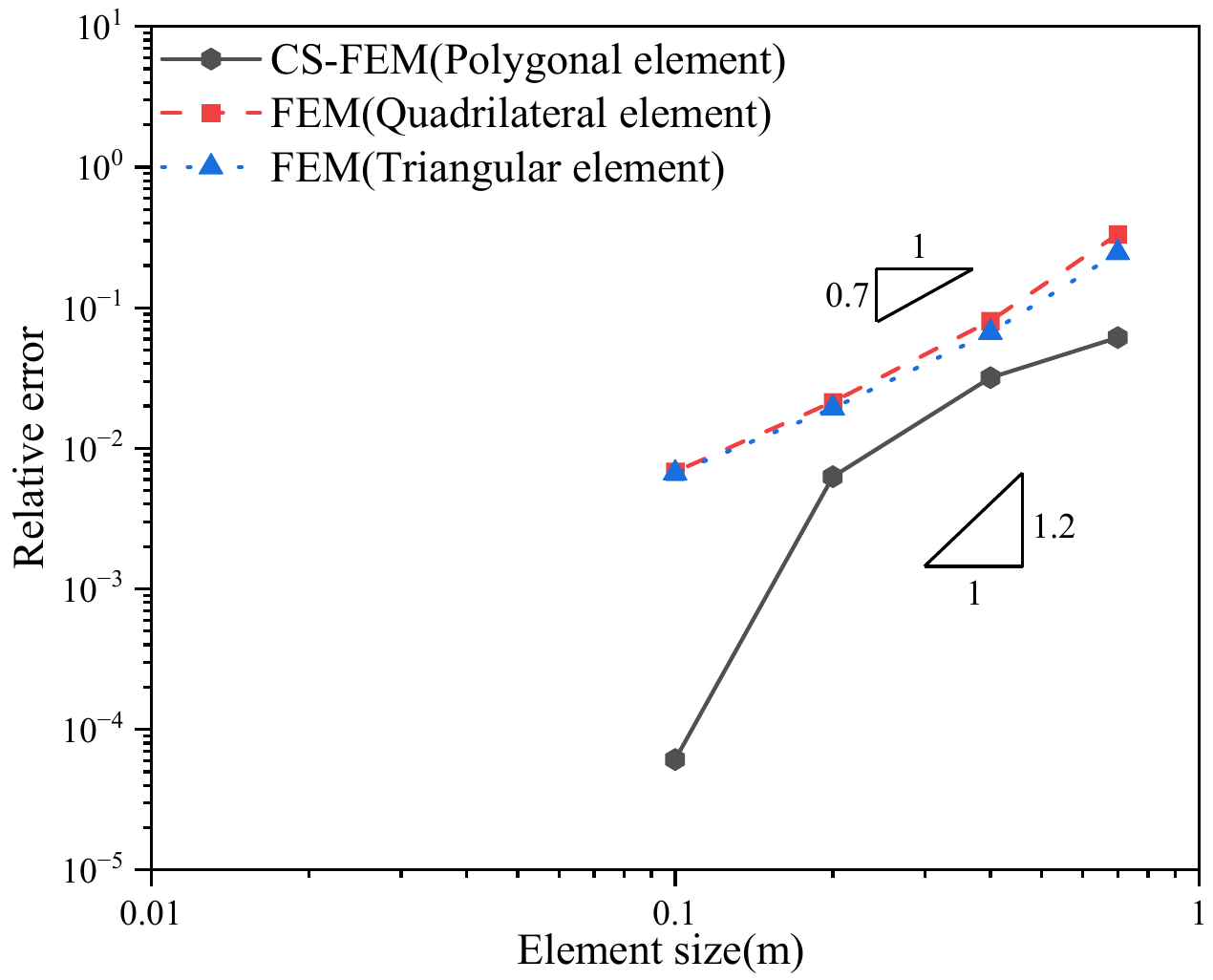}
  \caption{}
  \label{fig:1.0serror}
 \end{subfigure}
 \vspace{0cm}
 \begin{subfigure}[b]{0.4\textwidth}
  \centering
  \includegraphics[width=1\textwidth]{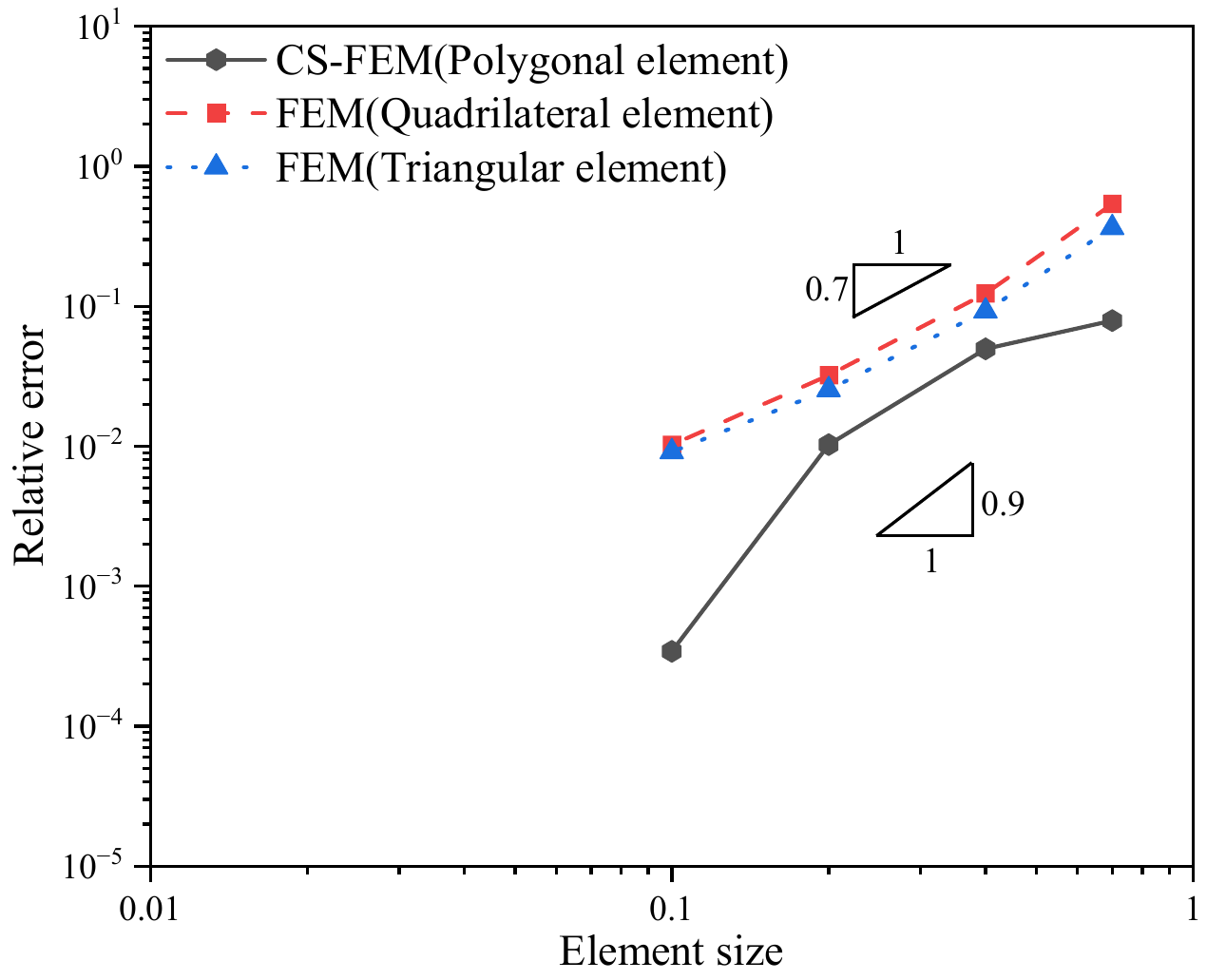}
  \caption{}
  \label{fig:1.5serror}
 \end{subfigure}
 \hfill
 \begin{subfigure}[b]{0.4\textwidth}
  \centering
  \includegraphics[width=1\textwidth]{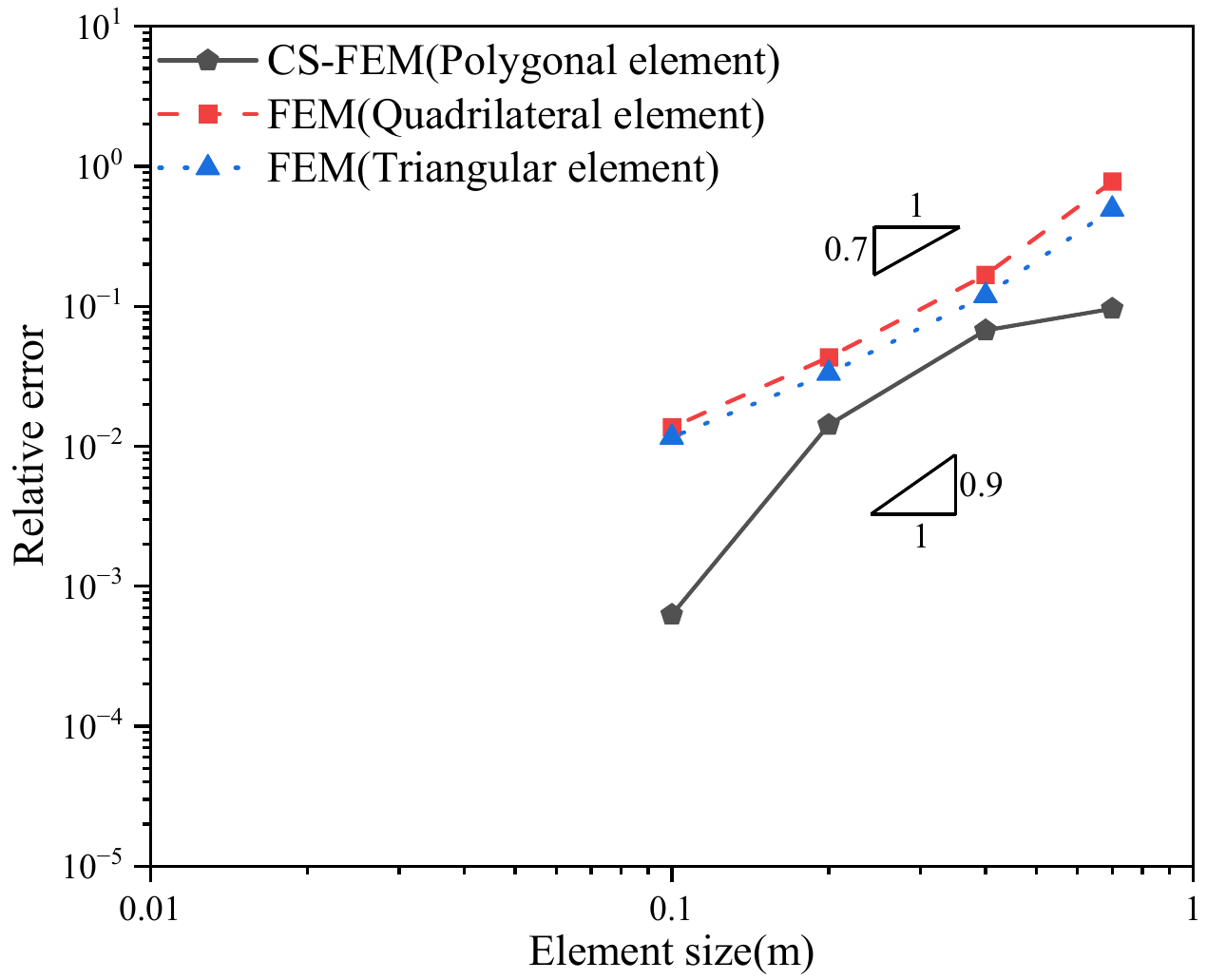}
  \caption{}
  \label{fig:2.0serror}
 \end{subfigure}
 \caption{Convergence of the relative temperature error at different time instants: (a) $t=0.5~\mathrm{s}$; (b) $t=1.0~\mathrm{s}$; (c) $t=1.5~\mathrm{s}$; and (d) $t=2.0~\mathrm{s}$.}
 \label{fig:Convergence_of_the_relative_error_in_the_temperature_at_the_different_times}
\end{figure}

Fig.~\ref{fig:Convergence_of_the_relative_error_in_the_temperature_at_the_different_times} presents the convergence curves of the relative temperature error at different time instants. It can be observed that the relative errors decrease consistently with mesh refinement for all element types. Compared with conventional triangular and quadrilateral FEM, the polygonal CS-FEM achieves smaller errors under the same mesh size and exhibits a favorable convergence trend. These results further demonstrate the accuracy and robustness of the proposed method for transient heat-conduction problems.

\FloatBarrier

\subsubsection{A square plate with a trefoil-shaped cavity}

\begin{figure}[htbp]
 \centering
 \begin{subfigure}[b]{0.34\textwidth}
  \centering
  \includegraphics[width=1\textwidth]{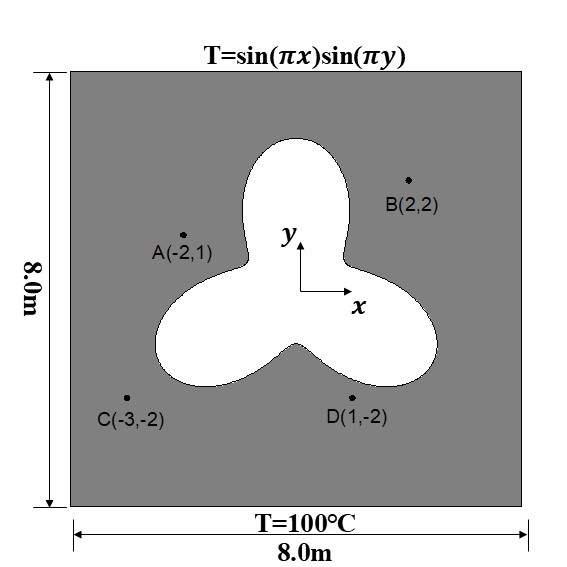}
  \caption{}
  \label{fig:geometry4}
 \end{subfigure}
 \hfill
 \begin{subfigure}[b]{0.3\textwidth}
  \centering
  \includegraphics[width=1\textwidth]{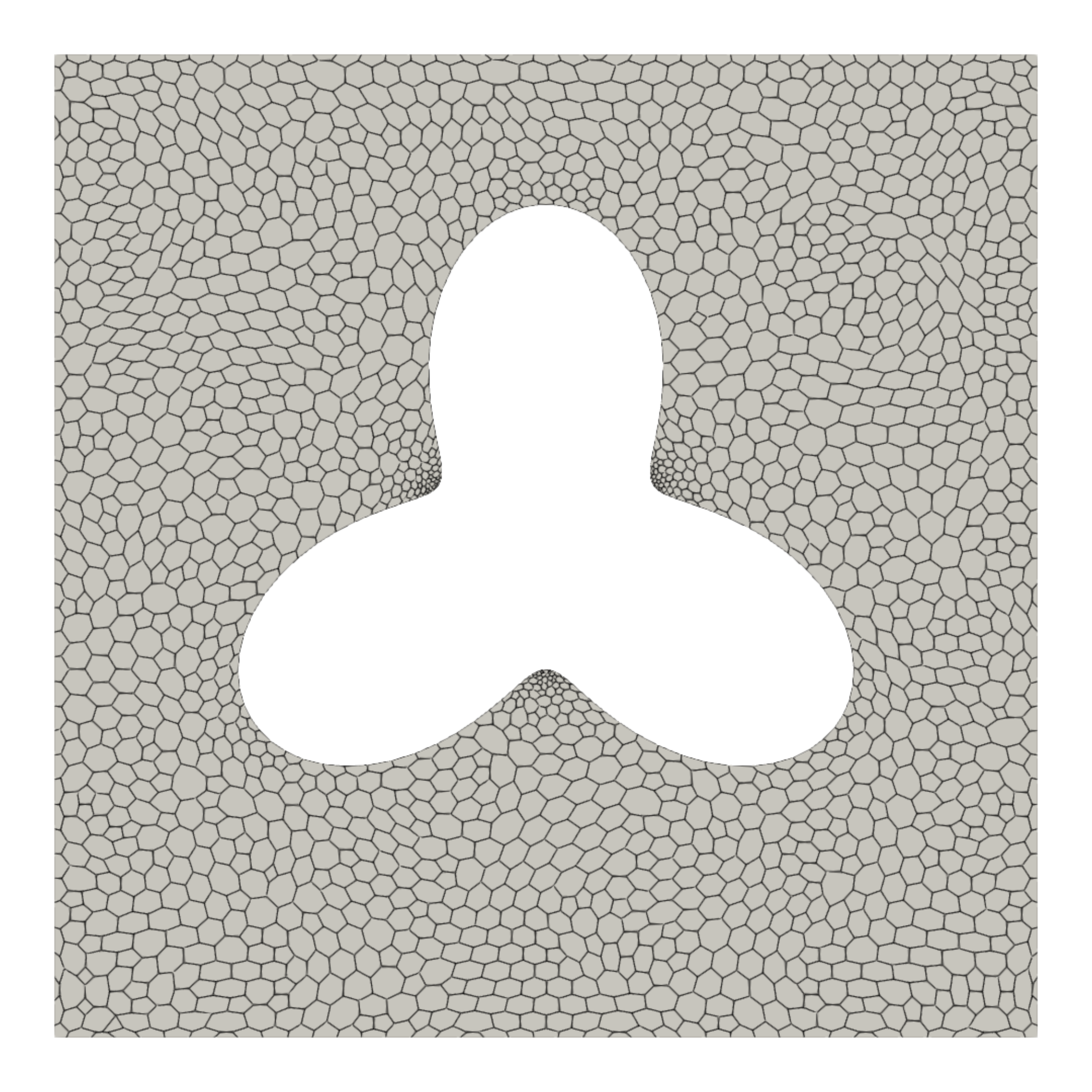}
  \caption{}
  \label{fig:4poly_M}
 \end{subfigure}
 \hfill
 \begin{subfigure}[b]{0.3\textwidth}
  \centering
  \includegraphics[width=1\textwidth]{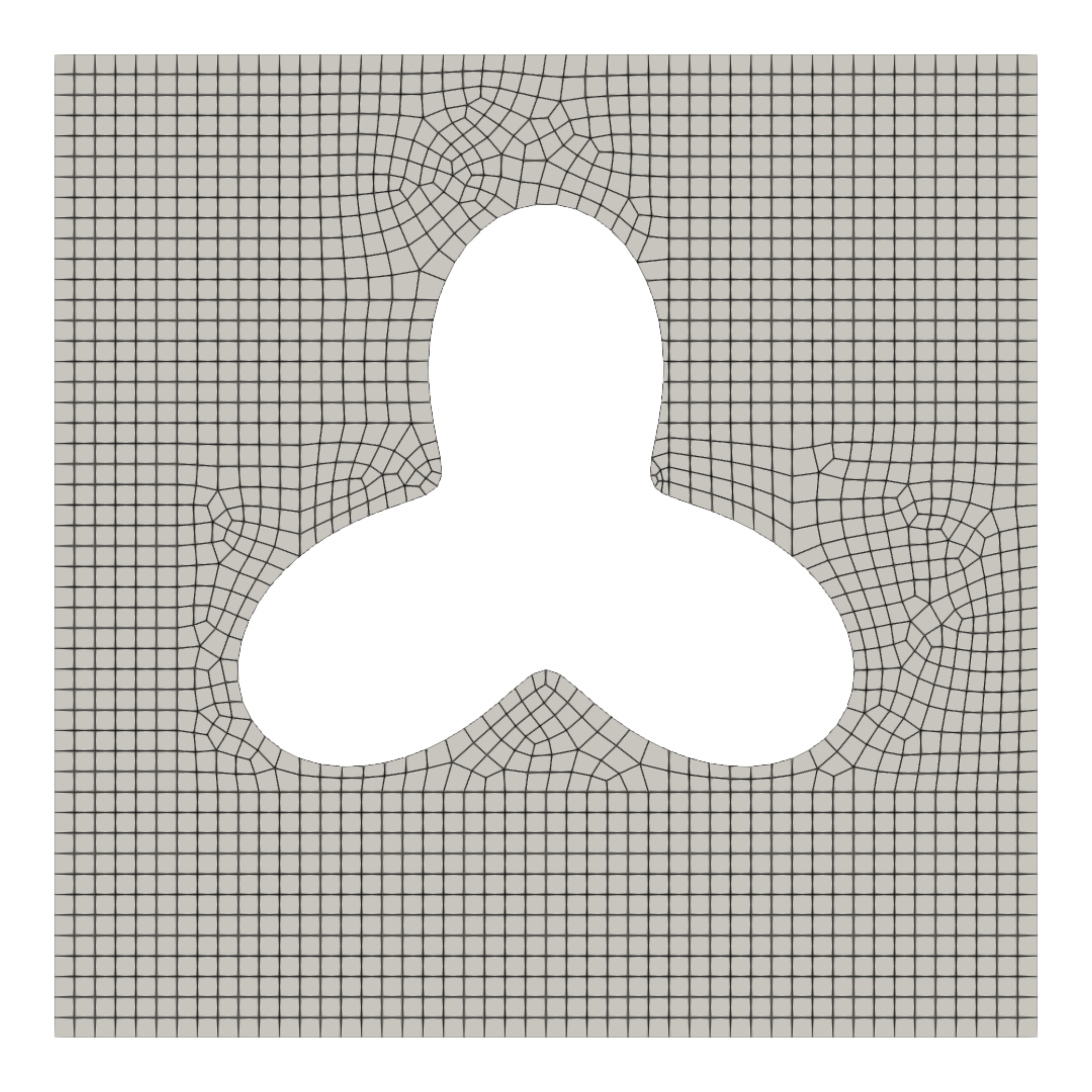}
  \caption{}
  \label{fig:4quad_M}
 \end{subfigure}
 \caption{Computational model of the square plate with a trefoil-shaped cavity: (a) geometry and boundary conditions; (b) polygonal mesh; and (c) quadrilateral mesh.}
 \label{fig:The_model_of_the_square_plate_with_a_trefoil-shape_cavity}
\end{figure}

This example considers transient heat conduction in a square plate containing a trefoil-shaped cavity, as shown in Fig.~\ref{fig:The_model_of_the_square_plate_with_a_trefoil-shape_cavity}. The side length of the square plate is $L=8.0~\mathrm{m}$. To quantitatively evaluate the performance of the proposed method, four representative points, denoted as A, B, C, and D, are selected in the computational domain. The temperature on the top boundary is prescribed as
\begin{equation}
T=\sin(\pi x)\sin(\pi y),
\end{equation}
whereas a constant temperature of $100~^{\circ}\mathrm{C}$ is applied on the bottom boundary. The time step is set to $\Delta t=1~\mathrm{s}$, and the total simulation time is $t=30~\mathrm{s}$.

The computational domain is discretized using 1951 polygonal elements. Since no analytical solution is available for this complex geometry, a conventional ABAQUS FEM solution obtained using a fine quadrilateral mesh with 1956 elements is adopted as the reference solution.

\FloatBarrier

\begin{figure}[htbp]
 \centering
 \begin{subfigure}[b]{0.32\textwidth}
  \centering
  \includegraphics[width=1\textwidth]{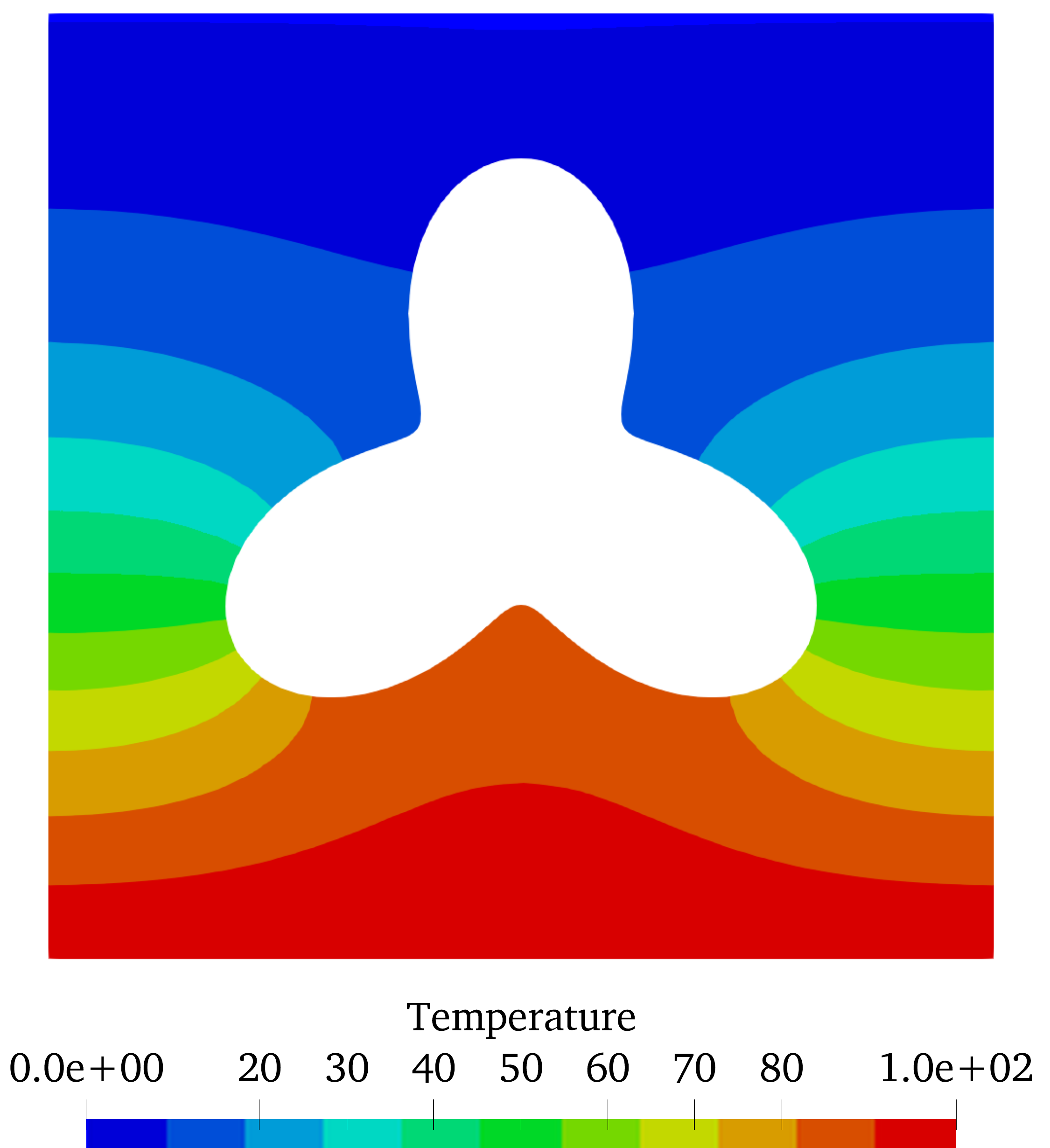}
  \caption{}
  \label{fig:4poly_T(10s)}
 \end{subfigure}
 \hfill
 \begin{subfigure}[b]{0.32\textwidth}
  \centering
  \includegraphics[width=1\textwidth]{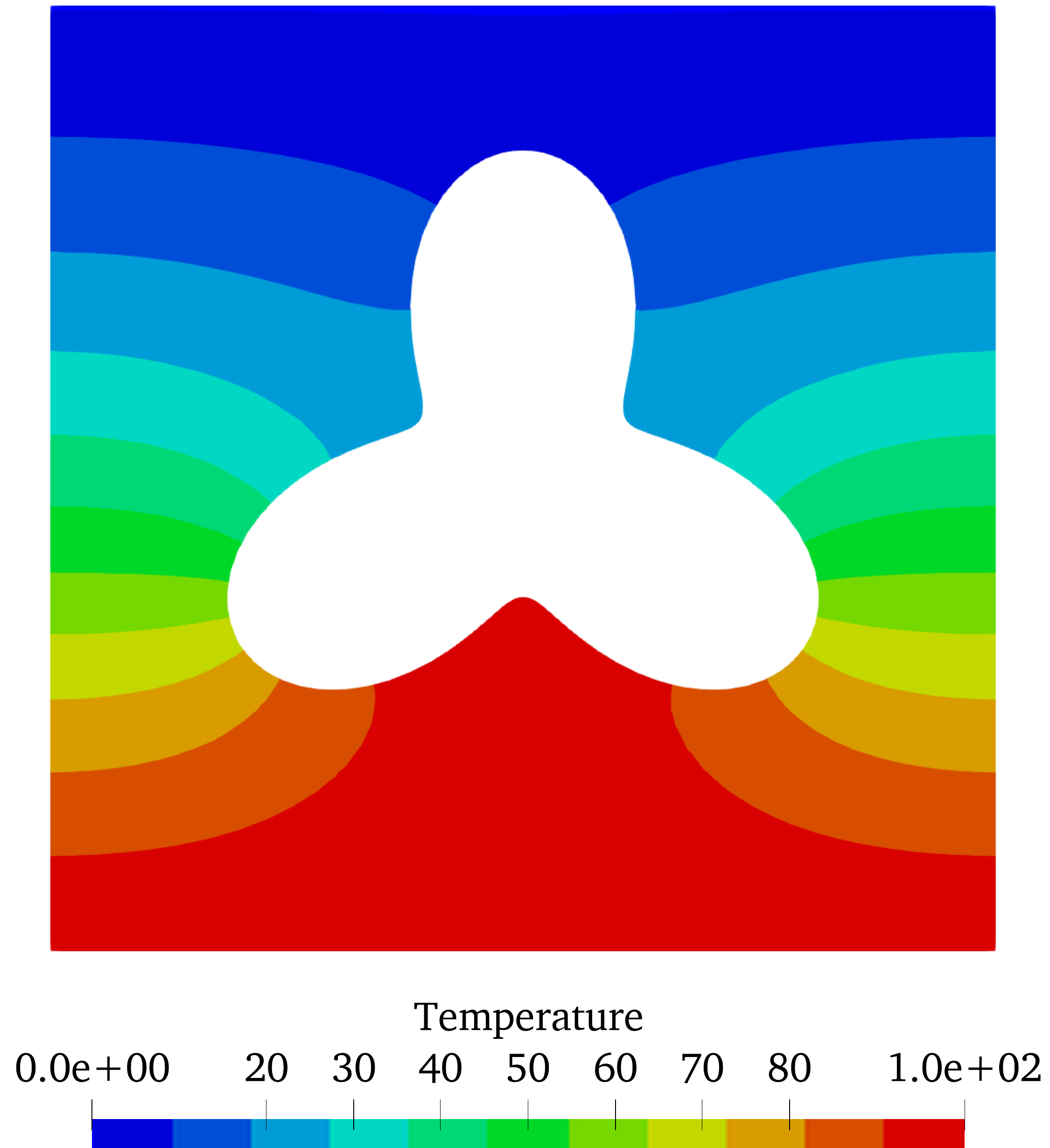}
  \caption{}
  \label{fig:4poly_T(20s)}
 \end{subfigure}
 \hfill
 \begin{subfigure}[b]{0.32\textwidth}
  \centering
  \includegraphics[width=1\textwidth]{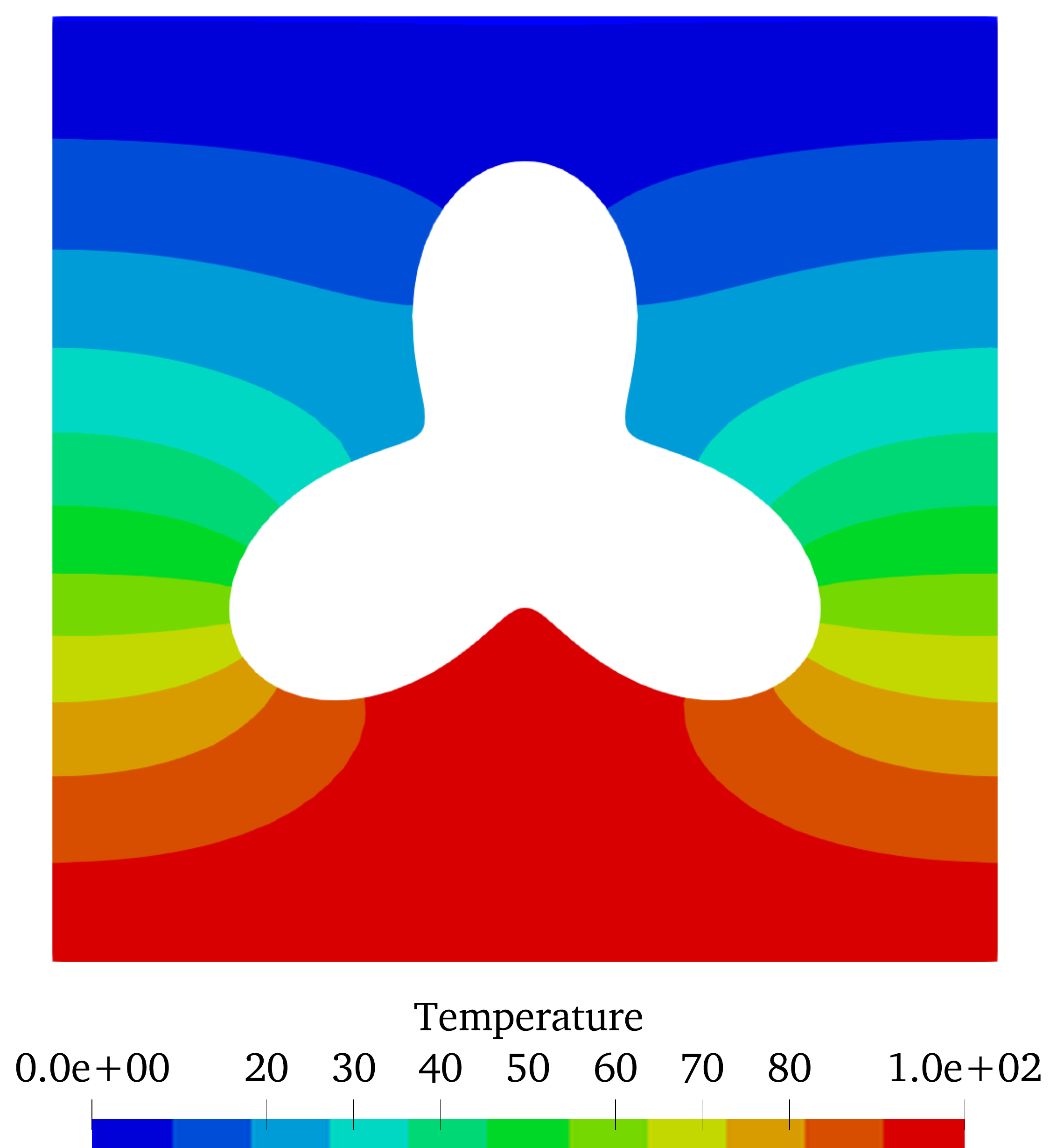}
  \caption{}
  \label{fig:4poly_T(30s)}
 \end{subfigure}
 \vspace{0cm}
 \begin{subfigure}[b]{0.32\textwidth}
  \centering
  \includegraphics[width=1\textwidth]{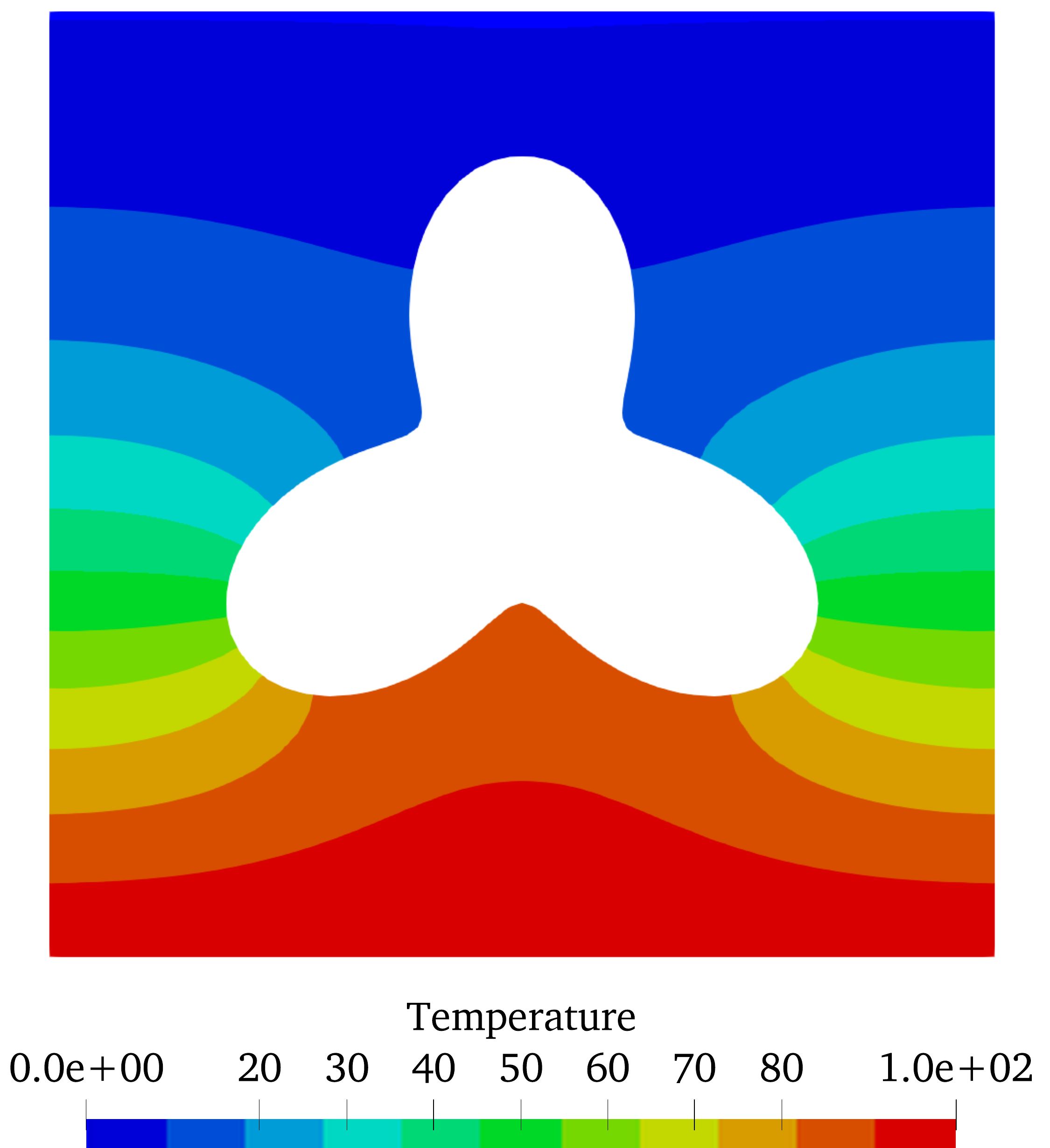}
  \caption{}
  \label{fig:4quad_T(10s)}
 \end{subfigure}
 \hfill
 \begin{subfigure}[b]{0.32\textwidth}
  \centering
  \includegraphics[width=1\textwidth]{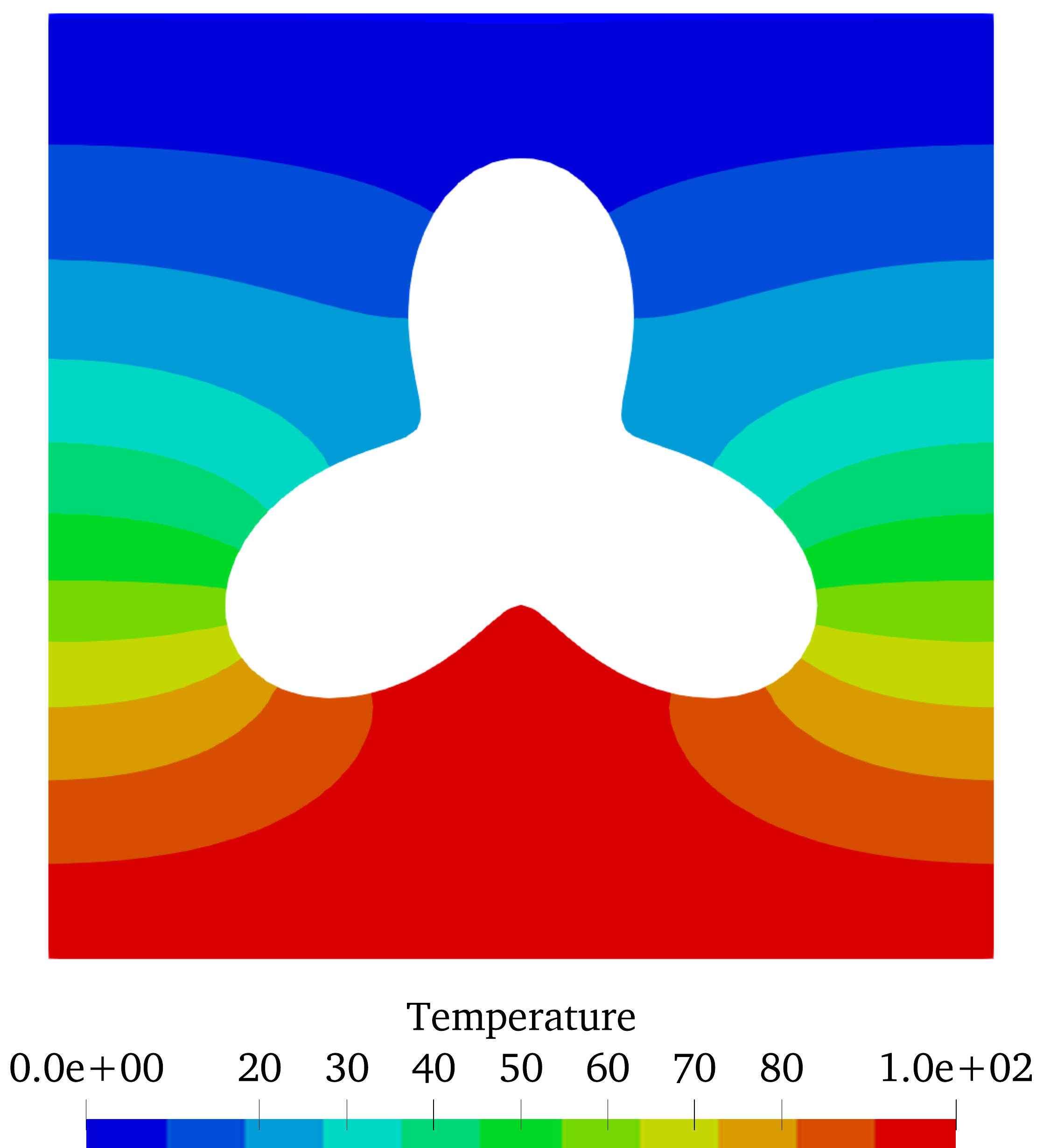}
  \caption{}
  \label{fig:4quad_T(20s)}
 \end{subfigure}
 \hfill
 \begin{subfigure}[b]{0.32\textwidth}
  \centering
  \includegraphics[width=1\textwidth]{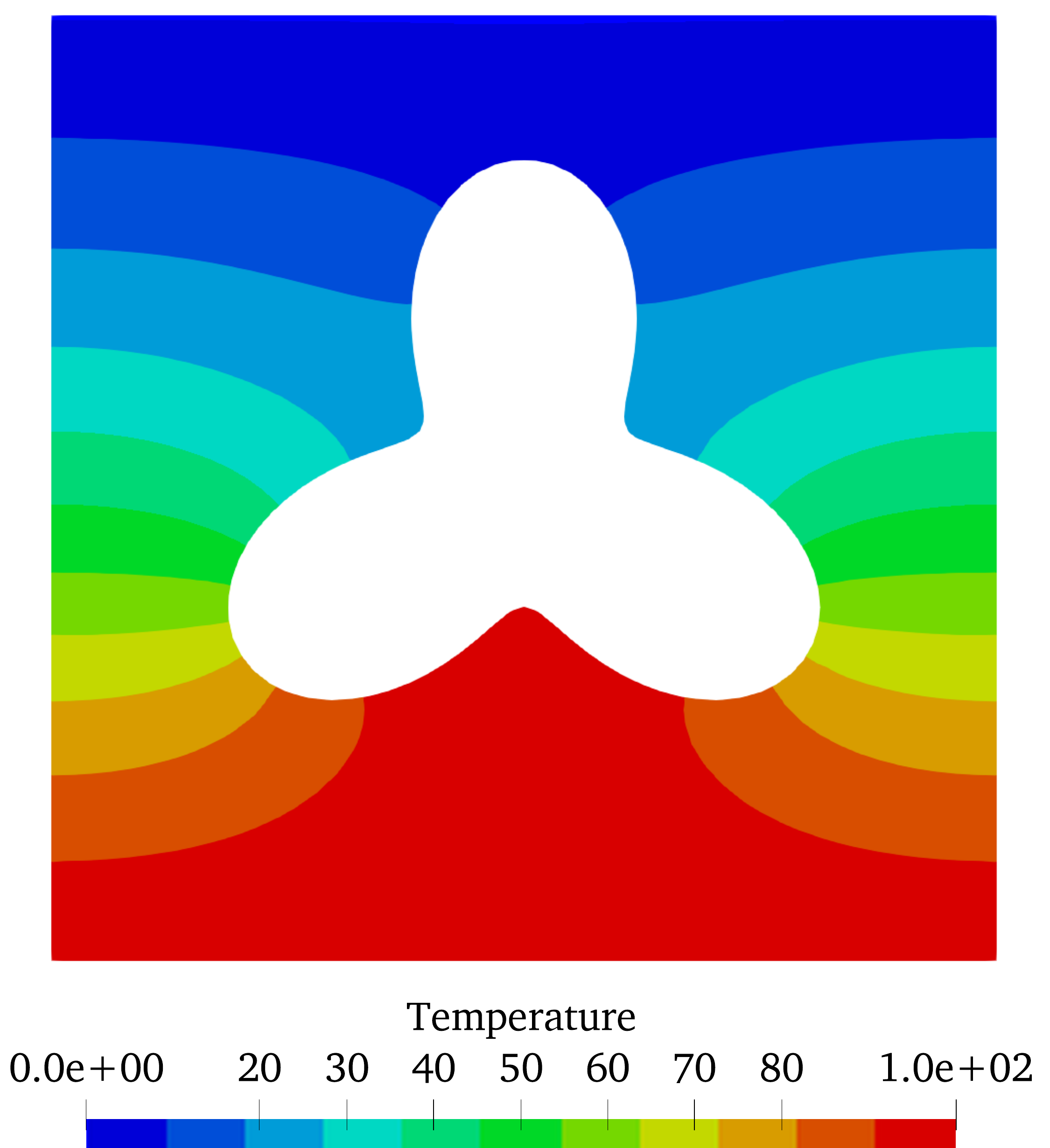}
  \caption{}
  \label{fig:4quad_T(30s)}
 \end{subfigure}
 \caption{Temperature distributions at different time instants obtained using the polygonal CS-FEM and FEM: (a) CS-FEM at $t=10~\mathrm{s}$; (b) CS-FEM at $t=20~\mathrm{s}$; (c) CS-FEM at $t=30~\mathrm{s}$; (d) FEM at $t=10~\mathrm{s}$; (e) FEM at $t=20~\mathrm{s}$; and (f) FEM at $t=30~\mathrm{s}$.}
 \label{fig:The_temperature_distribution_at_different_time_using_CS-FEM_and_FEM}
\end{figure}

Fig.~\ref{fig:The_temperature_distribution_at_different_time_using_CS-FEM_and_FEM} compares the temperature distributions obtained by the polygonal CS-FEM and FEM at $t=10~\mathrm{s}$, $20~\mathrm{s}$, and $30~\mathrm{s}$. The temperature contours predicted by the two methods are essentially identical at all selected time instants, indicating that the proposed method can accurately reproduce the transient temperature evolution in a complex domain.

\begin{figure}[htbp]
  \centering
  \includegraphics[width=0.8\textwidth]{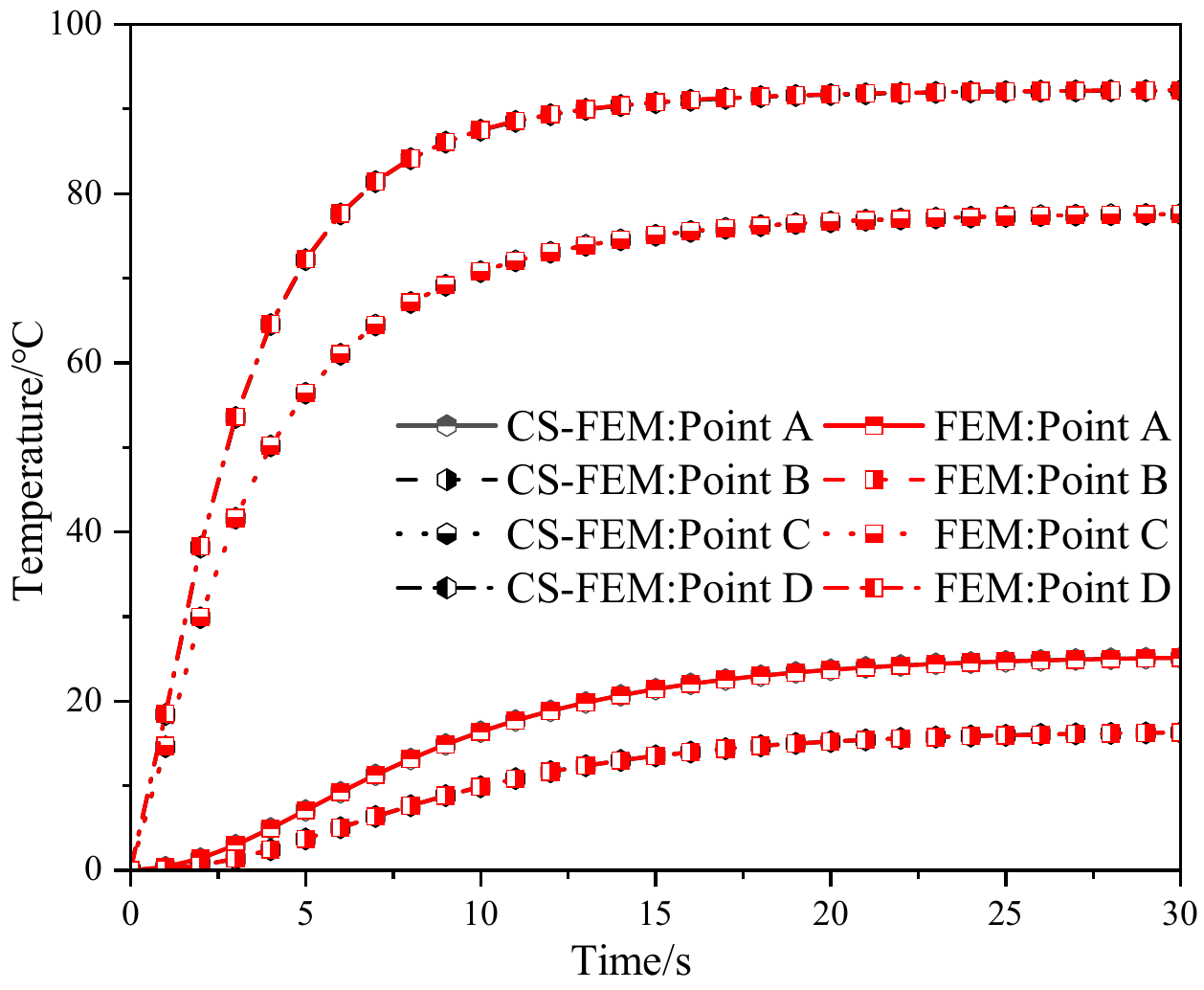}
  \caption{Temperature histories at the four selected points.}
  \label{fig:ABCDT}
\end{figure}

Fig.~\ref{fig:ABCDT} further presents the temperature histories at the four selected points. The temperature responses calculated by the polygonal CS-FEM closely match the reference solution throughout the simulation. These results demonstrate the effectiveness of the proposed method for transient heat-conduction problems involving complex geometries.

\FloatBarrier

\FloatBarrier

\subsubsection{A square domain subjected to heat flux}

\begin{figure}[htbp]
 \centering
 \begin{subfigure}[b]{0.34\textwidth}
  \centering
  \includegraphics[width=1\textwidth]{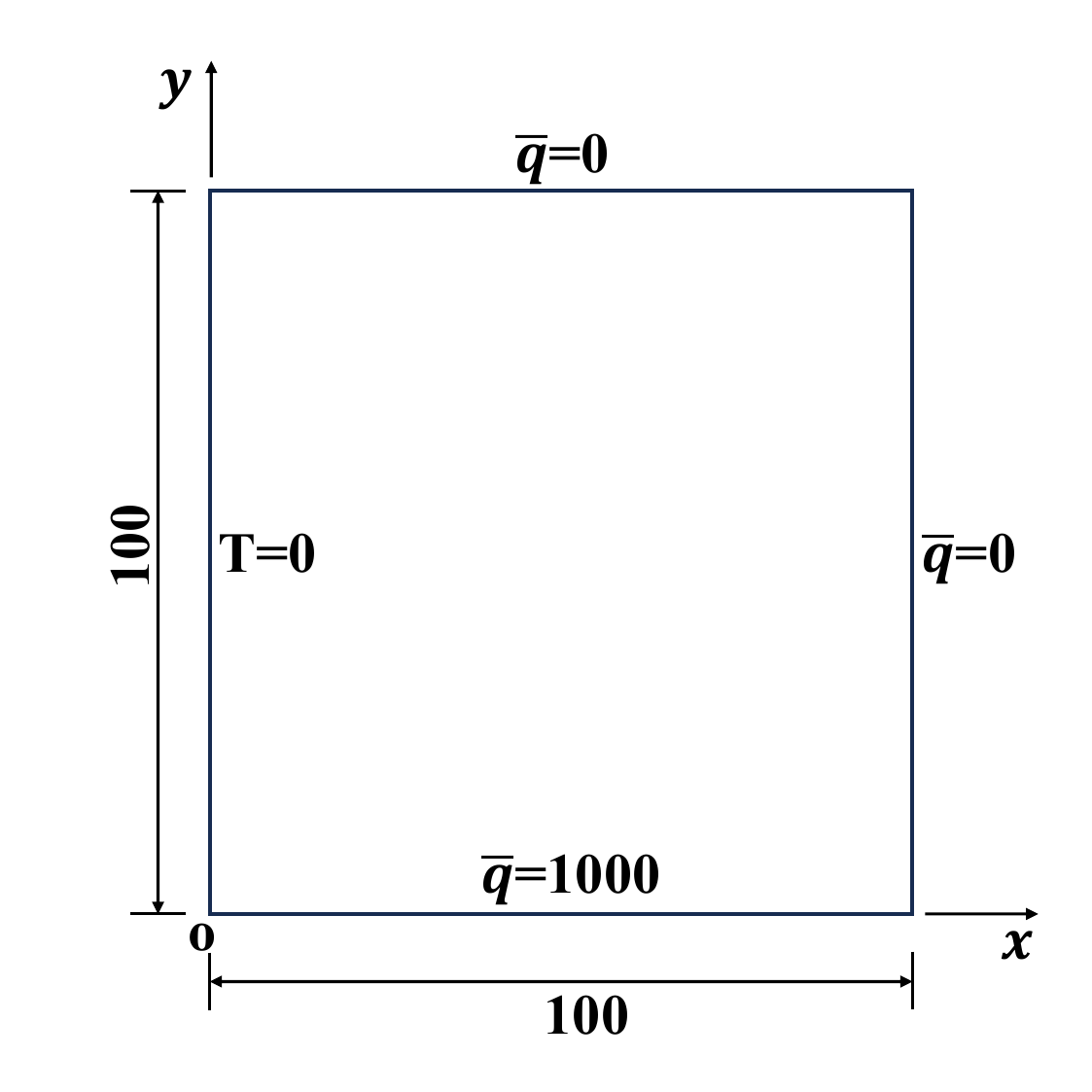}
  \caption{}
  \label{fig:geo}
 \end{subfigure}
 \hfill
 \begin{subfigure}[b]{0.3\textwidth}
  \centering
  \includegraphics[width=1\textwidth]{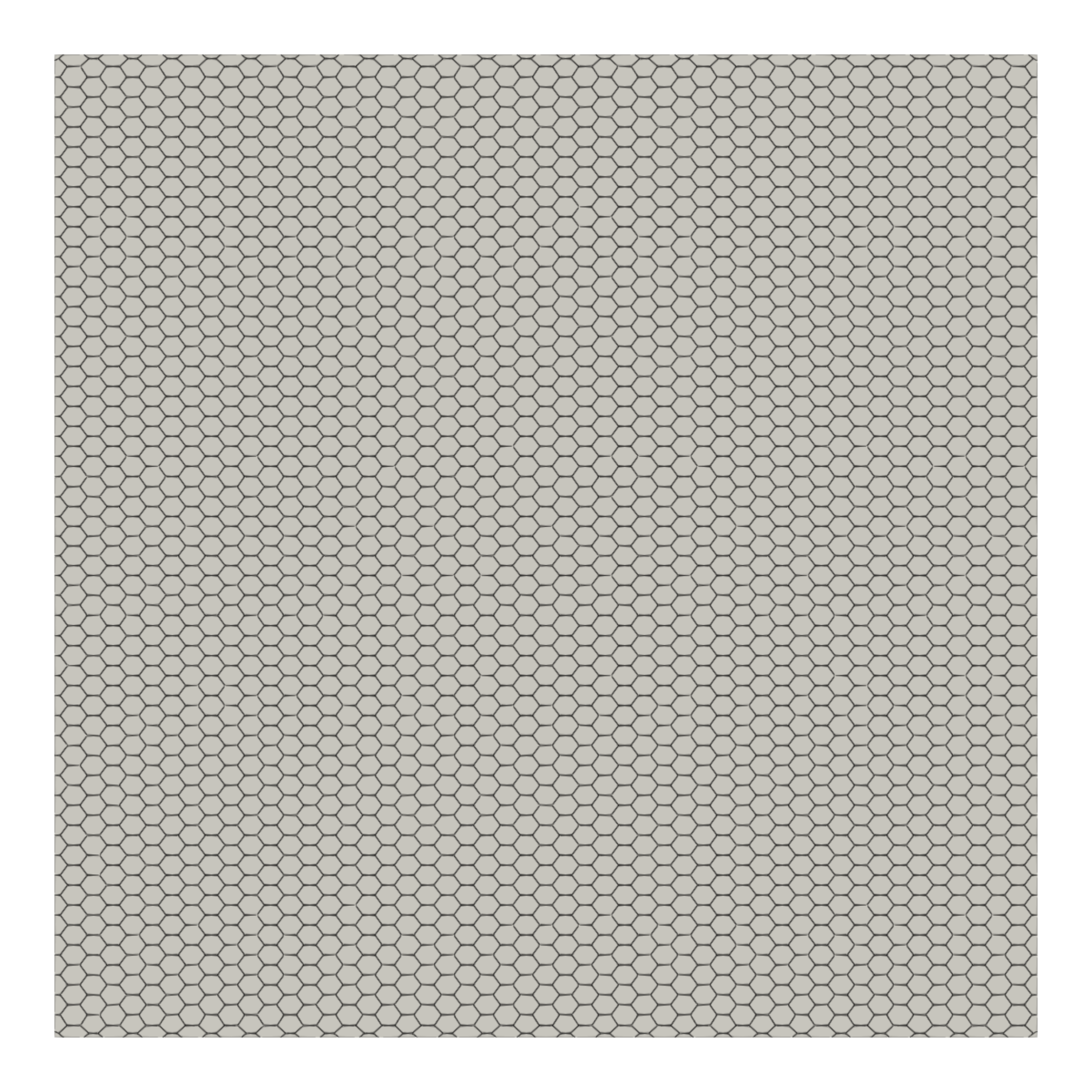}
  \caption{}
  \label{fig:pm}
 \end{subfigure}
 \hfill
 \begin{subfigure}[b]{0.3\textwidth}
  \centering
  \includegraphics[width=1\textwidth]{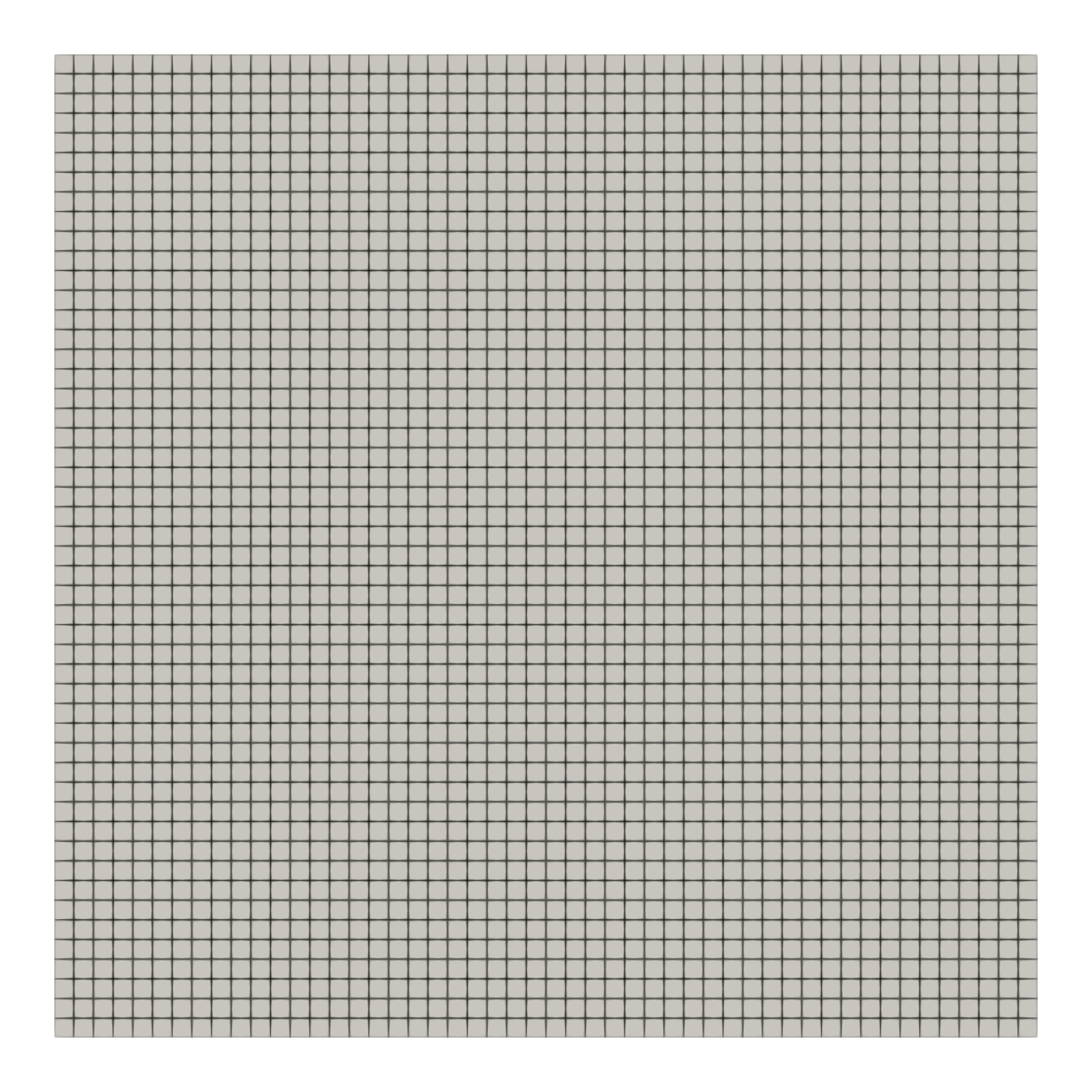}
  \caption{}
  \label{fig:qm}
 \end{subfigure}
 \caption{Computational model of the square domain subjected to heat flux: (a) geometry and boundary conditions; (b) polygonal mesh; and (c) quadrilateral mesh.}
 \label{fig:The_model_of_the_square_domain_with_heat_flux}
\end{figure}

A square domain with dimensions of $100 \times 100$ is considered as a benchmark problem to further examine the applicability of the proposed method under heat-flux boundary conditions. The geometry and boundary conditions are shown in Fig.~\ref{fig:The_model_of_the_square_domain_with_heat_flux}. The initial temperature is set to zero, and the thermal conductivity is specified as $k=1000$. A prescribed temperature of $T=0$ is imposed on the left boundary. A heat flux of $q=1000$ is applied along the bottom boundary, while the remaining boundaries are assumed to be thermally insulated. The time step is set to $\Delta t=0.02~\mathrm{s}$, and the total simulation time is $t=2~\mathrm{s}$.

The domain is discretized using 2576 polygonal elements for the polygonal CS-FEM analysis. For comparison, an ABAQUS FEM model using 2500 quadrilateral elements is established, and its solution is used as the reference solution.

\FloatBarrier

\begin{figure}[htbp]
 \centering
 \begin{subfigure}[b]{0.45\textwidth}
  \centering
  \includegraphics[width=1\textwidth]{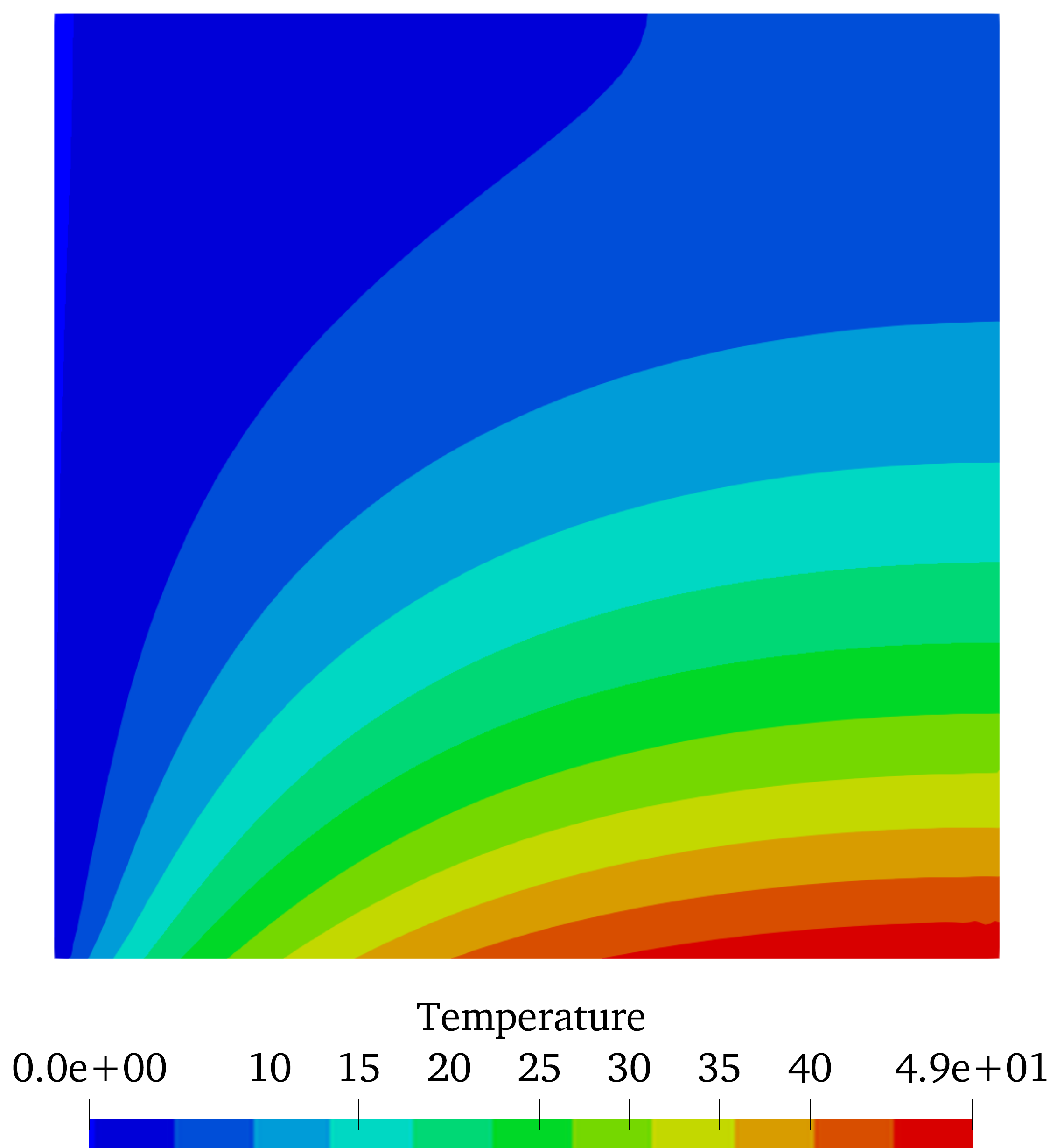}
  \caption{}
  \label{fig:pt}
 \end{subfigure}
 \hfill
 \begin{subfigure}[b]{0.45\textwidth}
  \centering
  \includegraphics[width=1\textwidth]{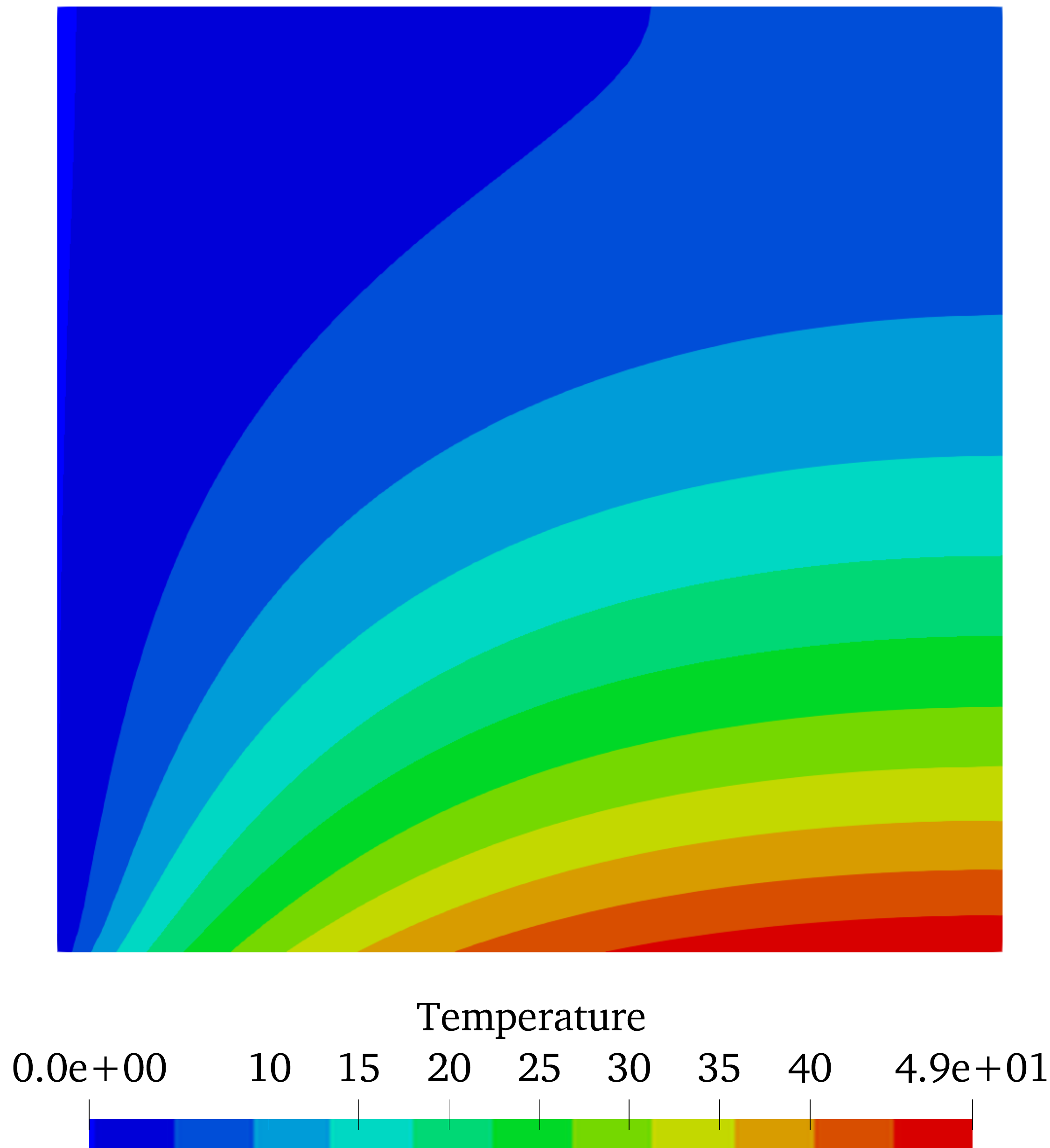}
  \caption{}
  \label{fig:qt}
 \end{subfigure}
 \caption{Temperature contour comparison at $t=2~\mathrm{s}$: (a) polygonal CS-FEM; and (b) FEM reference solution.}
 \label{fig:Temperature_contour_comparison}
\end{figure}

\begin{figure}[htbp]
 \centering
 \begin{subfigure}[b]{0.45\textwidth}
  \centering
  \includegraphics[width=1\textwidth]{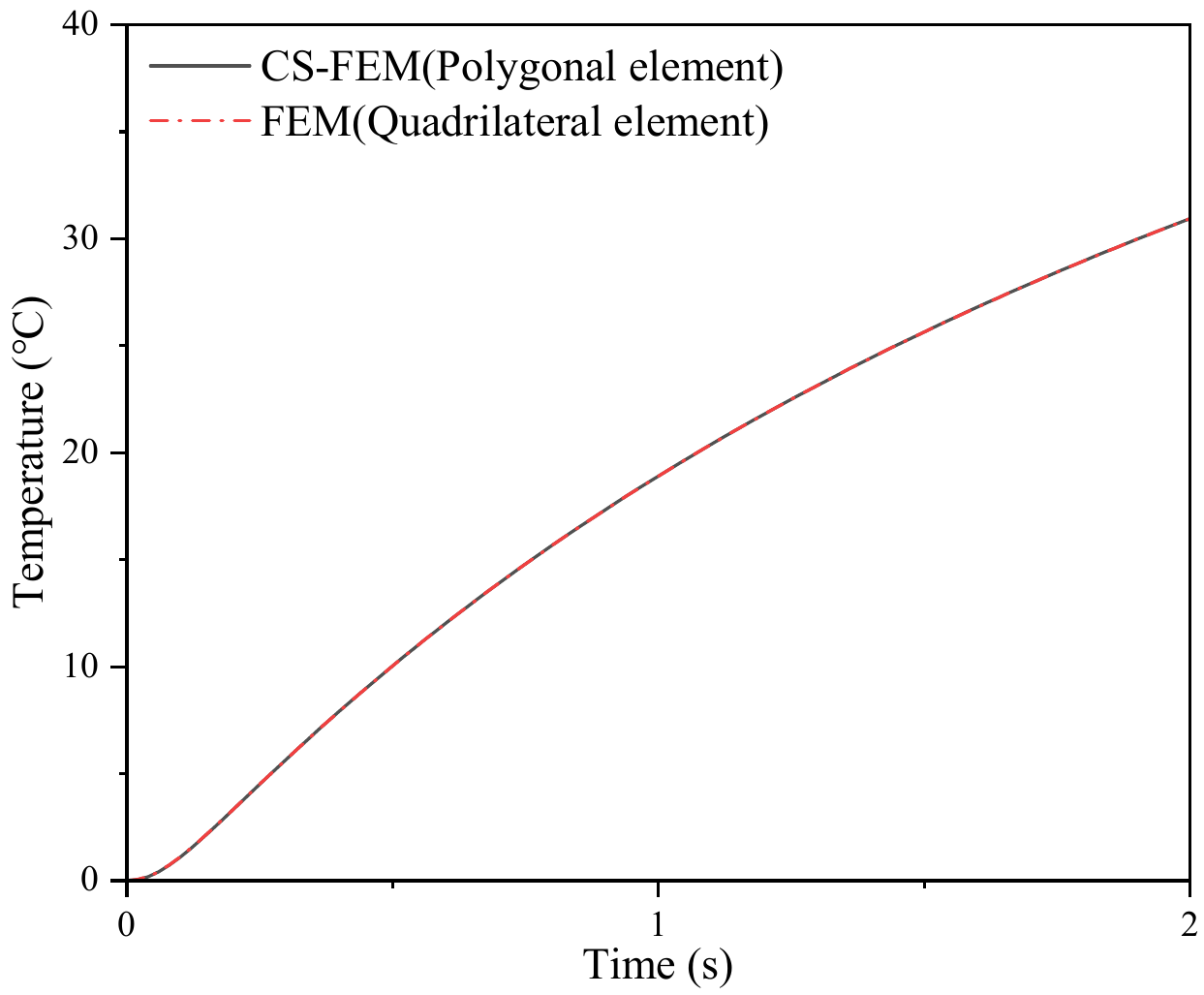}
  \caption{}
  \label{fig:AT}
 \end{subfigure}
 \hfill
 \begin{subfigure}[b]{0.45\textwidth}
  \centering
  \includegraphics[width=1\textwidth]{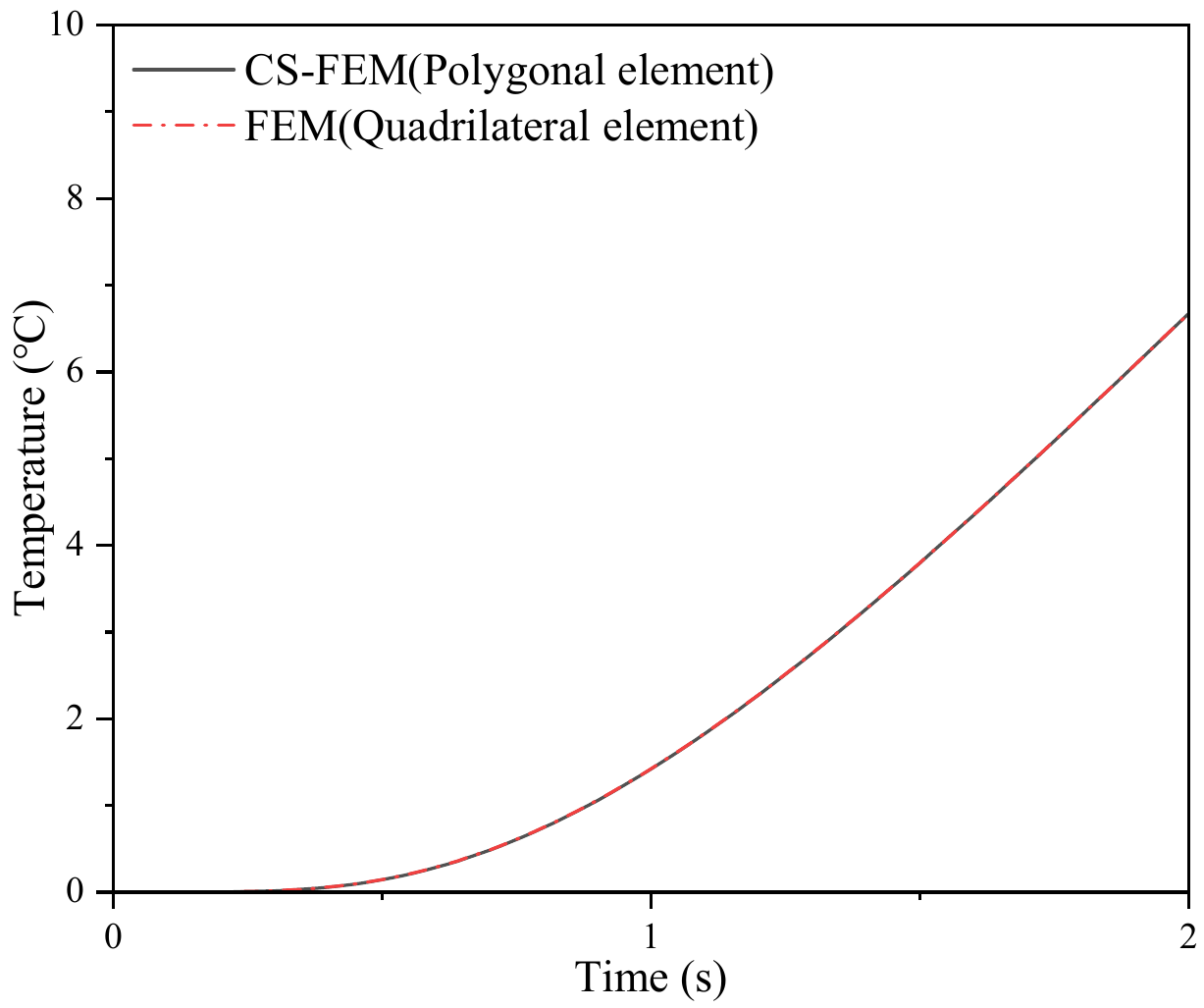}
  \caption{}
  \label{fig:BT}
 \end{subfigure}
 \caption{Temperature histories at points A and B: (a) point A; and (b) point B.}
 \label{fig:Variation_of_temperature_at_point_A_and_B}
\end{figure}

The temperature contours at $t=2~\mathrm{s}$ obtained using the polygonal CS-FEM and FEM are presented in Fig.~\ref{fig:Temperature_contour_comparison}. The comparison shows that the two temperature fields are nearly identical. Furthermore, the temperature histories at points A $(x=100, y=20)$ and B $(x=100, y=80)$ are plotted in Fig.~\ref{fig:Variation_of_temperature_at_point_A_and_B}. It can be observed that the results obtained by the proposed polygonal CS-FEM are in good agreement with the reference FEM solution, indicating that the proposed method remains stable and accurate for transient heat-conduction problems involving prescribed heat-flux boundaries.

\begin{figure}[htbp]
  \centering
  \includegraphics[width=0.8\textwidth]{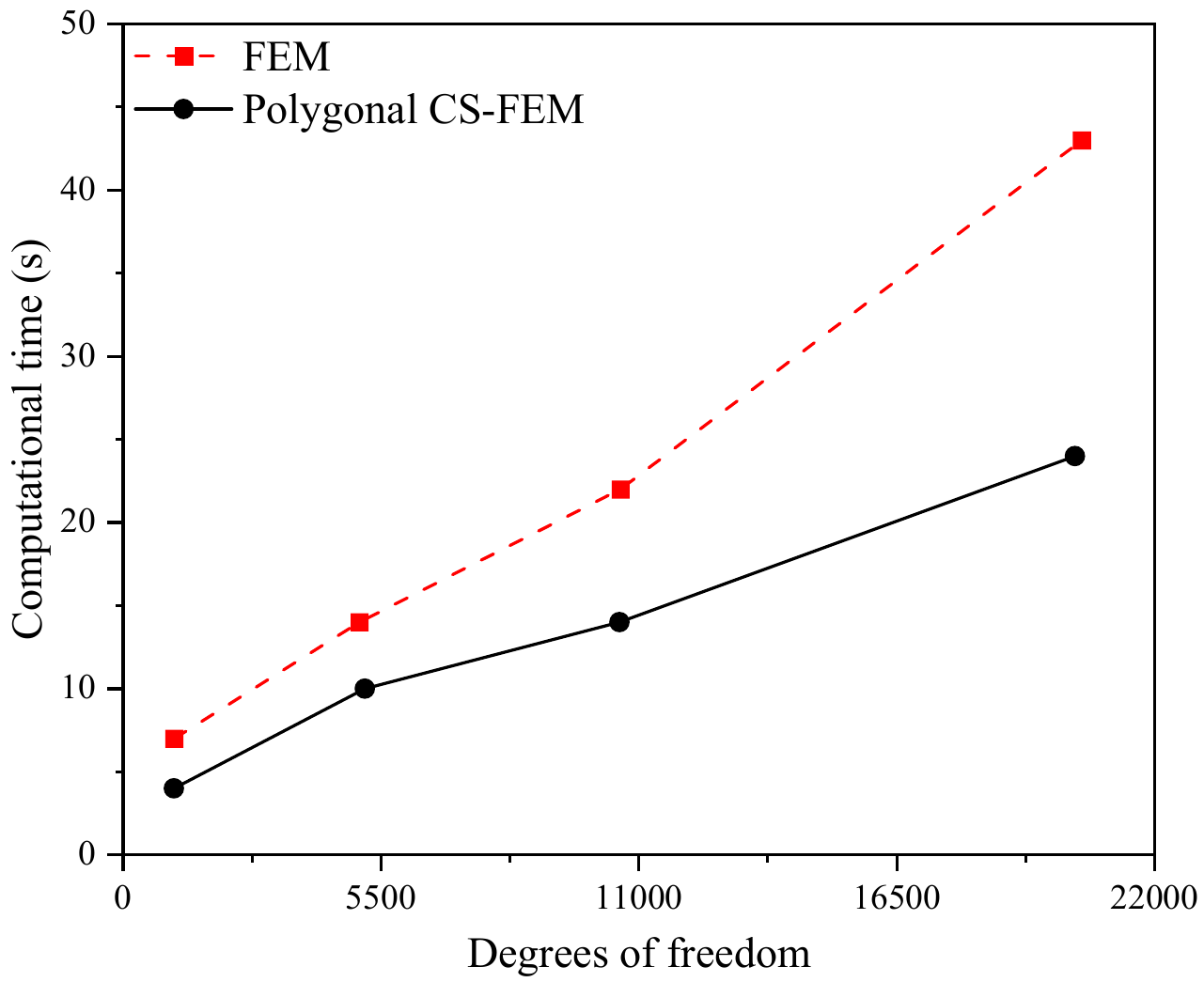}
  \caption{Comparison of computational time between the polygonal CS-FEM and FEM.}
  \label{fig:computational_time} 
\end{figure}

The comparison of computational time between the polygonal CS-FEM and FEM is shown in Fig.~\ref{fig:computational_time}. The computational time of both methods increases with the number of degrees of freedom. Under the same number of degrees of freedom, the polygonal CS-FEM requires less computational time than the conventional FEM, with a reduction of approximately 30\%--40\%. This demonstrates the computational efficiency of the proposed method.

\FloatBarrier

\FloatBarrier

\FloatBarrier

\section{Conclusions}
\label{sec:5}
This paper presents a polygonal cell-based smoothed finite element method (CS-FEM) for two-dimensional steady-state and transient heat-conduction analysis. In the proposed framework, Wachspress shape functions are used to construct polygonal elements, and the smoothed temperature gradient is evaluated through boundary integration over cell-based smoothing domains. Based on the numerical results, the following conclusions can be drawn.

(1) The proposed polygonal CS-FEM can accurately solve steady-state and transient heat-conduction problems. The numerical results agree well with analytical or reference solutions under different boundary conditions, and the linear patch test confirms that the method satisfies the first-order consistency requirement.

(2) The use of polygonal elements improves the flexibility of mesh generation and enables the method to effectively handle complex geometries. Compared with conventional triangular and quadrilateral elements, polygonal discretizations provide better adaptability to irregular domains and curved boundaries.

(3) The cell-based smoothing operation provides a stable and accurate approximation of the temperature gradient. The convergence studies show that the proposed method achieves favorable convergence behavior and generally provides accuracy comparable to or higher than that of conventional FEM under similar mesh resolutions.

Overall, the proposed polygonal CS-FEM provides an accurate, stable, and flexible numerical framework for two-dimensional heat-conduction analysis. Future work will focus on extending the present method to three-dimensional heat-conduction problems and coupled thermo-mechanical analyses.

\section{Acknowledgements}
The Fundamental Research Funds for the Central Universities (grant NO. B240201147), the Xing Dian Talent Support Program of Yunnan Province (grant NO. XDYC-QNRC-2022-0764), the Yunan Fundamental Research Projects (grant NO. 202401CF070043), and the Key Science and Technology Project of POWERCHINA Ltd (grant NO. DJ-ZDXM-2024-45) provided support for this study.
\bibliographystyle{elsarticle-num} 
\bibliography{cas-refs}



\end{document}